\documentclass[11pt,letterpaper]{article}
\usepackage[T1]{fontenc}
\usepackage[margin=1in]{geometry}
\usepackage{amsmath,amsthm,amssymb,mathtools,booktabs,array}
\usepackage{lmodern}
\usepackage{microtype}
\usepackage[authoryear,round]{natbib}
\usepackage{hyperref}
\hypersetup{
  colorlinks=true,
  linkcolor=blue,
  citecolor=blue,
  urlcolor=blue,
  pdftitle={Sharp Cramer-type moderate deviation bounds for the self-normalized maximum of partial sums},
  pdfauthor={Songhao Liu, Qiman Shao, Hao Shi, Jingyu Xu},
  pdfsubject={},
  pdfkeywords={Cramer-type moderate deviations, random-walk maximum, self-normalized bounds, sharp relative error}
}
\numberwithin{equation}{section}
\theoremstyle{plain}
\newtheorem{theorem}{Theorem}[section]
\newtheorem{proposition}[theorem]{Proposition}
\newtheorem{lemma}[theorem]{Lemma}
\newtheorem{corollary}[theorem]{Corollary}
\theoremstyle{definition}

\newcommand{\E}{\mathbb E}
\newcommand{\PP}{\mathbb P}
\newcommand{\R}{\mathbb R}

\newcommand{\ind}{\mathbf 1}
\newcommand{\Phibar}{\overline\Phi}
\newcommand{\var}{\operatorname{Var}}

\newcommand{\supp}{\operatorname{supp}}
\newcommand{\ii}{\mathrm i}
\newcommand{\dd}{\,\mathrm d}
\newcommand{\ee}{\mathrm e}

\newcommand{\calN}{\mathcal N}
\newcommand{\calC}{\mathcal C}
\newcommand{\calL}{\mathcal L}
\newcommand{\calF}{\mathcal F}

\newcommand{\Cend}{C_{\mathrm{end},n}}

\newcommand{\norm}[1]{\lVert#1\rVert}

\newcommand{\Mref}[1]{\ref{#1}}
\newcommand{\Sref}[1]{\ref{#1}}
\providecommand{\bmhead}[1]{\section*{#1}}
\allowdisplaybreaks[1]

\title{Sharp Cram\'er-type moderate deviation bounds for the self-normalized maximum of partial sums}
\author{Song-Hao Liu$^{1}$ \and Qi-Man Shao$^{2}$ \and Hao Shi$^{3}$ \and Jing-Yu Xu$^{3}$}
\date{}
\begin{document}
\maketitle
\vspace{-1.5em}
\begin{center}
\small
\begin{minipage}{0.98\textwidth}
\centering
$^{1}$School of Mathematical Sciences, Dalian University of Technology,\\
Dalian, Liaoning, China\\[0.4em]
$^{2}$Department of Statistics and Data Science,\\
Shenzhen International Center for Mathematics,\\
Southern University of Science and Technology,\\
Shenzhen, Guangdong, 518055, China\\[0.4em]
$^{3}$Department of Statistics and Data Science,\\
Southern University of Science and Technology,\\
Shenzhen, Guangdong, 518055, China\\[0.7em]
\href{mailto:liusonghao@dlut.edu.cn}{\texttt{liusonghao@dlut.edu.cn}};\quad
\href{mailto:shaoqm@sustech.edu.cn}{\texttt{shaoqm@sustech.edu.cn}}\\
\href{mailto:Hao.Shi@link.cuhk.edu.hk}{\texttt{Hao.Shi@link.cuhk.edu.hk}};\quad
\href{mailto:12131253@mail.sustech.edu.cn}{\texttt{12131253@mail.sustech.edu.cn}}
\end{minipage}
\end{center}

\begin{abstract}
  Let $X_1,X_2,\ldots$ be i.i.d.  random variables with zero mean and variance $\sigma^2\in(0,\infty)$.
In this paper, we establish a sharp Cram\'er-type moderate deviation for the self-normalized maximum of partial sums under a finite $p$-th absolute moment, $2<p\le3$.
In particular, we show that $\PP(M_n\ge xV_n)=2(1-\Phi(x))\{1+O(1)\beta_{n,p}(1+x)^p\}$ for $0\le x\le\beta_{n,p}^{-1/p}$, where $M_n=\max_{1\le k\le n}\sum_{j=1}^kX_j$, $V_n=(\sum_{j=1}^nX_j^2)^{1/2}$, $\beta_{n,p}=\E|X_1|^p/(n^{(p-2)/2}\sigma^p)$, and $\Phi$ is the standard normal distribution function.
This proves a conjecture of \citet{Shao2010}.
More precisely, we bound the relative error by the same truncated-moment quantity as for self-normalized sums.
The proof compares the maximum with the absolute sum while preserving their common denominator, through a reflected-walk representation and Fourier estimates with distribution-adapted smoothing.
\end{abstract}

\noindent\textbf{Keywords:} Cram\'er-type moderate deviations;
self-normalized maximum of partial sums; sharp relative error bound; Lindley recursion.
\medskip

\section{Introduction}\label{sec:introduction}
Let $X,X_1,X_2,\ldots$ be i.i.d.\ real random variables such that
\begin{equation}\label{eq:assumptions}
 \E X=0,\qquad 0<\sigma^2=\E X^2<\infty,\qquad
 \mu_p=\E|X|^p<\infty,\quad 2<p\le3.
\end{equation}
For $1\le k\le n$, set
\begin{equation}\label{eq:statistics}
 S_k^X=\sum_{j=1}^kX_j,\qquad
 M_n=\max_{1\le k\le n}S_k^X,\qquad V_n^2=\sum_{j=1}^nX_j^2.
\end{equation}
Write $\Phi$ and $\phi$ for the standard normal distribution function and density, and put $\Phibar=1-\Phi$. The self-normalized invariance principle of \citet{CsorgoSzyszkowiczWang2003}, together with the Brownian reflection principle, identifies $2\Phibar(x)$ as the reference upper tail for $M_n/V_n$. Our aim is to bound the relative error when the threshold $x\ge0$
is allowed to grow with $n$.

For the terminal sum $S_n^X$, self-normalization replaces the unknown scale $\sigma\sqrt n$ by  $V_n$ and permits Cram\'er--type moderate deviations without exponential-moment assumptions; see \citet{deLaPenaLaiShao2009}. In particular, \citet{Shao1999} showed that a finite third absolute moment suffices for the relative normal approximation of $S_n^X/V_n$ to hold uniformly over $0\le x\le a_n$, where $a_n=o(n^{1/6})$.

The relevant error scale for the terminal self-normalized sum $S_n^X/V_n$ combines a second-moment tail with a truncated third moment. For $x\ge0$, set $t_{n,x}=\sigma\sqrt n/(1+x)$ and define
\begin{equation}\label{eq:truncatedscale}
 R_n(x)=\frac{(1+x)^2}{\sigma^2}
 \left\{\E[X^2\ind_{\{|X|>t_{n,x}\}}]
 +t_{n,x}^{-1}\E[|X|^3\ind_{\{|X|\leq t_{n,x}\}}]\right\}.
\end{equation}
\citet{JingShaoWang2003} established relative error bounds in terms of this truncated-moment quantity for independent, not necessarily identically distributed, observations. In the i.i.d.\ setting, their error quantity is exactly $R_n(x)$, and their bounds yield
\begin{equation}\label{eq:jswbenchmark}
 \left|\frac{\PP(S_n^X\ge xV_n)}{\Phibar(x)}-1\right|
 \le C R_n(x),\qquad
 0\le x\le \frac{\sigma n^{(p-2)/(2p)}}{\mu_p^{1/p}},
\end{equation}
where $C$ is an absolute constant. Subsequent work includes refined tail approximations by \citet{Wang2011} and \citet{SangGe2017}, and an extension to nonlinear statistics by \citet{ShaoZhou2016}.

For maxima normalized by $\sigma\sqrt n$, \citet{Aleshkyavichene1979} obtained Cram\'er-type expansions under exponential-moment assumptions. For self-normalized maxima, \citet{HuShaoWang2009} established a Cram\'er-type relative-tail approximation under a finite fourth moment. \citet{LiuShaoWang2013} extended the approximation to a finite $p$-th moment, $2<p\le3$. In the i.i.d.\ setting, their result gives
\[
 \frac{\PP(M_n\ge xV_n)}{\Phibar(x)}\longrightarrow2
 \quad\text{uniformly for }0\le x\le a_n,
 \qquad a_n=o\bigl(n^{(p-2)/(2p)}\bigr).
\]
The remaining question is whether the maximum admits the same error bound $R_n(x)$ as the self-normalized sum in \eqref{eq:jswbenchmark}.

Write
\begin{equation}\label{eq:beta}
 \beta_{n,p}=\frac{\mu_p}{n^{(p-2)/2}\sigma^p}.
\end{equation}
We prove that, under \eqref{eq:assumptions},
\[
 \left|\frac{\PP(M_n\ge xV_n)}{\Phibar(x)}-2\right|
 \le C R_n(x),\qquad 0\le x\le\beta_{n,p}^{-1/p},
\]
with an absolute constant $C$, independent of the distribution, $n$, $x$ and $p$. In particular, the relative error tends to zero uniformly on $0\le x\le a_n\le\beta_{n,p}^{-1/p}$ whenever $R_n(a_n)\to0$. A bound in terms of the $p$th absolute moment follows immediately.
Since $R_n(x)\le \beta_{n,p}(1+x)^p$, the relative error for the maximum is at most $C\beta_{n,p}(1+x)^p$ throughout $0\le x\le\beta_{n,p}^{-1/p}$. This proves the bound conjectured by \citet[Section~3.3, Conjecture~1]{Shao2010}, whose finite-third-moment case was also stated by \citet[Remark~1]{LiuShaoWang2013}.

For terminal self-normalized sums $S_n^X$, \citet{ChistyakovGotze2004} established sharpness of the relative-error bounds in \citet{JingShaoWang2003}. We show that the scale $R_n$ is also sharp for the maximum $M_n$, even for a fixed distribution. In fact, for each $p\in(2,3]$, we construct a single symmetric distribution satisfying \eqref{eq:assumptions}. The corresponding absolute relative-tail error is asymptotically equivalent to a positive multiple of $R_n$; see Theorem~\ref{thm:sharp}. 
Thus the upper bound cannot in general be replaced by an $o(R_n)$ bound.

{The proof compares the maximum with the absolute self-normalized sum $|S_n^X|/V_n$, whose upper tail has the same normal approximation $2\Phibar(x)$. 
Using the Lindley reflected-walk representation and a telescoping identity, we obtain a comparison between the maximum $M_n/V_n$ and the absolute self-normalized sum $|S_n^X|/V_n$. Fourier estimates with convolution smoothing then control the comparison error between $M_n/V_n$ and $|S_n^X|/V_n$. Finally, a two-sided comparison removes the smoothing.}

\textbf{Notation} Throughout, $C,c$ denote positive absolute constants that may change between occurrences. A subscript on $C$ or $c$ indicates dependence on a fixed parameter. Unless otherwise specified, constants in the upper bounds do not depend on the increment distribution.

\textbf{Organization} Section~\ref{sec:mainresults} states the two main theorems
and their Lyapunov-moment consequence.
Section~\ref{sec:mainproof} proves the main theorems.
The three supporting propositions are proved in
Sections~\ref{sec:preliminary}--\ref{sec:comparison}.
The proofs of the auxiliary lemmas are collected in the Supplement.

\section{Main results}\label{sec:mainresults}
We use the weak event $\{M_n\ge xV_n\}$, including the all-zero sample, rather than impose a quotient convention when $V_n=0$. The same convention is used for $\{S_n^X\ge xV_n\}$; the zero-denominator correction is included in its error bound, as verified in Section~\ref{sec:preliminary}.

\begin{theorem}\label{thm:main}
Let $X,X_1,X_2,\ldots$ be i.i.d.\ real random variables with $\E X=0$, $\sigma^2=\E X^2\in(0,\infty)$ and $\mu_p=\E|X|^p<\infty$ for some $2<p\le3$. For $n\ge1$, write $M_n=\max_{1\le k\le n}\sum_{j=1}^kX_j$, $V_n^2=\sum_{j=1}^nX_j^2$ and $\beta_{n,p}=\mu_p/(n^{(p-2)/2}\sigma^p)$. For $x\ge0$, set $t_{n,x}=\sigma\sqrt n/(1+x)$ and define
\[
 R_n(x)=\frac{(1+x)^2}{\sigma^2}
 \left\{\E[X^2\ind_{\{|X|>t_{n,x}\}}]
 +t_{n,x}^{-1}\E[|X|^3\ind_{\{|X|\le t_{n,x}\}}]\right\}.
\]
There is an absolute constant $C$ such that, for every $0\le x\le\beta_{n,p}^{-1/p}$,
\begin{equation}\label{eq:main}
 \left|\frac{\PP(M_n\ge xV_n)}{\Phibar(x)}-2\right|
 \le C R_n(x).
\end{equation}
The constant is independent of the distribution, $n$, $x$ and $p\in(2,3]$.
\end{theorem}

\begin{theorem}\label{thm:sharp}
For each $p\in(2,3]$, there is a symmetric distribution $P_*$ with countably infinite support such that
\[
 0<\E_{P_*}X^2<\infty,\qquad
 \E_{P_*}|X|^p<\infty,\qquad
 \E_{P_*}|X|^{p+\eta}=\infty\quad(\eta>0),
\]
and there are integers $N_k\uparrow\infty$ and thresholds $x_k\ge0$ for which an infinite i.i.d.\ sample from $P_*$ satisfies
\begin{equation}\label{eq:sharp-limit}
 \frac{\mathcal E_{N_k}(x_k;P_*)}{R_{N_k}(x_k;P_*)}
 \longrightarrow\kappa_p>0.
\end{equation}
Here $R_n(x;P)$ and $\beta_{n,p}(P)$ denote $R_n(x)$ and $\beta_{n,p}$ from Theorem~\ref{thm:main}, computed under the law $P$, and $\mathcal E_n(x;P)=\left|\PP_P(M_n\ge xV_n)/\Phibar(x)-2\right|$.
Moreover, $x_k\beta_{N_k,p}(P_*)^{1/p}\to0$ and $R_{N_k}(x_k;P_*)\to0$.
\end{theorem}

The original Lyapunov-moment bound follows from Theorem~\ref{thm:main}.

\begin{corollary}[Lyapunov bound]\label{cor:lyapunov}
Under the assumptions of Theorem~\ref{thm:main},
\begin{equation}\label{eq:budget-lyapunov}
 R_n(x)\le\beta_{n,p}(1+x)^p,\qquad x\ge0.
\end{equation}
Consequently,
\begin{equation}\label{eq:lyapunov-main}
 \left|\frac{\PP(M_n\ge xV_n)}{\Phibar(x)}-2\right|
 \le C\beta_{n,p}(1+x)^p,\qquad 0\le x\le\beta_{n,p}^{-1/p}.
\end{equation}
\end{corollary}
\section{Proofs of the main theorems}\label{sec:mainproof}

We prove Theorem~\ref{thm:main} by comparing the tail probabilities
$\PP(M_n\ge xV_n)$ and $\PP(|S_n^X|\ge xV_n)$.
Proposition~\ref{prop:preliminary} provides the normal approximation
for the latter probability and settles the cases where $x$ is close
to zero or $R_n(x)$ is not small.
For the remaining thresholds,
Propositions~\ref{prop:reflection} and~\ref{prop:smoothed}
provide the reflection identity and the smoothed comparison.
We remove the smoothing in Step~4.
The supporting propositions are proved in
Sections~\ref{sec:preliminary}--\ref{sec:comparison}.

\subsection{Proof of Theorem~\ref{thm:main}}\label{sec:proofpropositions}
\begin{proof}
We use the assumptions and notation of Theorem~\ref{thm:main}, and write $S_k^X=\sum_{j=1}^kX_j$. For $t>0$, define the functions $C_n$ and $\Cend$,
the reference function $G$, and the error weight $w_n$ by
\begin{align}
 C_n(t)&=\PP(M_n<tV_n),&
 \Cend(t)&=\PP(|S_n^X|<tV_n),\label{eq:cdfs}\\
 G(t)&=2\Phi(t)-1,&w_n(t)&=R_n(t)\Phibar(t).\label{eq:weight}
\end{align}
In particular,
\begin{equation}\label{eq:complement}
 E_n(t):=C_n(t)-G(t)=2\Phibar(t)-\PP(M_n\ge tV_n),\qquad t>0.
\end{equation}
Thus, for $0<t\le\beta_{n,p}^{-1/p}$, the required relative-tail bound is equivalent to $|E_n(t)|\le Cw_n(t)$.

\noindent\textbf{Step 1: preliminary bounds and truncation.}\par\nobreak

Write $r_n=R_n(0)$.
We first record the properties of $R_n$ used in the preliminary estimates.

\begin{lemma}\label{lem:truncatedbudget}
Let $F(z)=z^2\min\{|z|,1\}$ and $\xi=X/(\sigma\sqrt n)$. The function $R_n(t)=n\E F((1+t)\xi)$ is continuous and strictly increasing on $[0,\infty)$, satisfies $R_n(t)\le \beta_{n,p}(1+t)^p$, and obeys $R_n(t)/(1+t)^2\to1$ as $t\to\infty$. Moreover,
\begin{equation}\label{eq:budgetcomparison}
 \tfrac18 R_n(x)\le R_n(t)\le8R_n(x),\qquad
 x,t\ge0,\quad |t-x|\le(1+x)^{-1}.
\end{equation}
If $R_n(t)\le\rho\le1/16$, then
\begin{equation}\label{eq:budgetlower}
 n\ge\frac{(1-\rho)^3(1+t)^6}{R_n(t)^2}\ge(1+t)^6.
\end{equation}
\end{lemma}
The proof is given in Supplemental Material~\Sref{supp:truncatedbudgetproof}.

The following estimates control the maximum, the absolute sum, and the truncation error.

\begin{proposition}\label{prop:preliminary}
There are absolute constants $\rho_*\in(0,1/16]$ and $C>0$
with the following property. For $n\ge1$, define
\begin{equation}\label{eq:comparisoninterval}
 \mathcal I_n
 =\{t\ge0:t\le\beta_{n,p}^{-1/p}\ \text{or}\ R_n(t)\le\rho_*\}.
\end{equation}
Then, for every $t\in\mathcal I_n$, $t>0$, and constant $H>0$,
\begin{align}
  \PP(M_n\ge tV_n)
 &\le C\Phibar(t),\label{eq:coarse}\\
 |\Cend(t)-G(t)|
 &\le Cw_n(t),\label{eq:endbound}\\
 \PP(M_n\ge tV_n,\max_{j\le n}|X_j|>H)
 &\le Cn\PP(|X|>H)\Phibar(t).\label{eq:badlarge}
\end{align}
In addition,
\begin{equation}\label{eq:nearzero}
 |\PP(M_n\ge tV_n)-2\Phibar(t)|
 \le C(t+r_n),\qquad 0\le t\le1.
\end{equation}
\end{proposition}
The proof is given in Section~\ref{sec:preliminary}.

Fix $\rho_*$ as in Proposition~\ref{prop:preliminary}.
We shall prove
\[
 |E_n(t)|\le Cw_n(t),\qquad t\in\mathcal I_n,\quad t>0.
\]
Since $[0,\beta_{n,p}^{-1/p}]\subseteq\mathcal I_n$,
this implies Theorem~\ref{thm:main} for positive thresholds.
The enlarged interval accommodates the neighboring thresholds
used in smoothing and desmoothing.

By \eqref{eq:endbound} and the decomposition
\[
 E_n(t)=\{C_n(t)-\Cend(t)\}+\{\Cend(t)-G(t)\},
\]
we reduce the problem to comparing $C_n$ and $\Cend$.

We choose the constants as follows.
Fix sufficiently small absolute constants $\varepsilon_0\in(0,1]$ and $\eta_0>0$.
Fix an absolute $L\ge1600$ and set
\begin{equation}\label{eq:proofparameters}
 K=8(1+L),\qquad A_0=4KL,\qquad \varepsilon=\min\left\{\frac{\varepsilon_0}{2},
                    \frac1{4(K+1)L}\right\}.
\end{equation}
For $0<u\le\varepsilon_0$, put
\[
 \eta(u)=\min\left\{\frac{u^3}{8},\frac{u^2\eta_0}{40},
                       \frac{\rho_*}{32},\frac1{4A_0}\right\}.
\]
We divide $\mathcal I_n$ into three regions:

\begin{equation}\label{eq:1}
  \begin{aligned}
  \mathcal{D}_1 =& \mathcal I_n\cap \{t\ge0: R_n(t)\ge \eta(\varepsilon)\},\qquad\\
  \mathcal{D}_2 = &\mathcal I_n\cap \{t\ge0: R_n(t)< \eta(\varepsilon), t\le A_0r_n\},\qquad\\
  \mathcal{D}_3 =& \mathcal I_n\cap \{t\ge0: R_n(t)< \eta(\varepsilon), t>A_0r_n\},
  \end{aligned}
\end{equation}

Applying \eqref{eq:nearzero} at $t=0$ and using $\Phibar(0)=1/2$ gives
\[
 \left|\frac{\PP(M_n\ge0)}{\Phibar(0)}-2\right|
 =2|\PP(M_n\ge0)-1|\le2CR_n(0).
\]
This proves the assertion at $x=0$. For the rest of the proof, let $x>0$.

For $\mathcal{D}_1$,
by \eqref{eq:complement}, the triangle inequality, and the  bound \eqref{eq:coarse},
\[
 \begin{aligned}
 |E_n(x)|
 &\le2\Phibar(x)+\PP(M_n\ge xV_n)\le(C+2)\Phibar(x)
 \le\frac{C+2}{\eta(\varepsilon)}w_n(x).
 \end{aligned}
\]

In $\mathcal{D}_2$, $R_n$ is increasing by Lemma~\ref{lem:truncatedbudget}. Thus $r_n=R_n(0) \le R_n(x)<\eta(\varepsilon)$, and hence $x\le A_0r_n<A_0\eta(\varepsilon)\le1/2$. Thus \eqref{eq:nearzero} applies, while $w_n(x)\ge r_n\Phibar(1)$. Together with \eqref{eq:complement}, these bounds give
\[
 |E_n(x)|\le C(x+r_n)\le C(A_0+1)r_n
 \le\frac{C(A_0+1)}{\Phibar(1)}w_n(x).
\]
Thus $|E_n(x)|\le Cw_n(x)$ in both ranges.
It remains to consider $\mathcal D_3$.

To compare the tails in the remaining region, we condition on bounded increments.
Fix $x\in\mathcal D_3$.
With $\varepsilon$ fixed in \eqref{eq:proofparameters}, define
\begin{equation}\label{eq:truncation}
 \xi=\frac{X}{\sigma\sqrt n},\qquad
 \xi_j=\frac{X_j}{\sigma\sqrt n},\qquad
 h_x=\frac{\varepsilon}{1+x},\qquad \calN_x=\bigcap_{j\le n}\{|\xi_j|\le h_x\},\qquad
 p_x=\PP(\calN_x).
\end{equation}
We use the conditional law given $\calN_x$ for this comparison.

For $t>0$, splitting the tail probability over $\calN_x$ and $\calN_x^c$ gives
\begin{equation}
\label{eq:002}
 E_n(t)
 =\{2\Phibar(t)p_x-\PP(M_n\ge tV_n,\calN_x)\} +\{2\Phibar(t)(1-p_x)-\PP(M_n\ge tV_n,\calN_x^c)\}.
\end{equation}

We first bound the terms involving $\calN_x^c$.
The definition of $R_n$ gives
\[
 R_n(x)=n\E\bigl[((1+x)\xi)^2
 \min\{(1+x)|\xi|,1\}\bigr].
\]
The union bound yields
\begin{equation}\label{eq:truncationprob}
 1-p_x=\PP(\calN_x^c)
 \le n\PP(|\xi|>h_x)
 \le\varepsilon^{-3}R_n(x).
\end{equation}
Since $\eta(\varepsilon)\le\varepsilon^3/8$, we have
\begin{equation}\label{eq:conditioningprob}
 p_x\ge3/4.
\end{equation}

Consider $t>0$ with $|t-x|\le(1+x)^{-1}$.
By \eqref{eq:budgetcomparison},
\begin{equation}\label{eq:neighborhood}
 R_n(t)\le8R_n(x)<8\eta(\varepsilon)\le\rho_*,\qquad t\in\mathcal I_n.
\end{equation}
In this range, $|t^2-x^2|\le3$ and
$1/2\le(1+t)/(1+x)\le2$.
Mills' inequalities and \eqref{eq:budgetcomparison} give
$\Phibar(t)\asymp\Phibar(x)$ and $w_n(t)\asymp w_n(x)$.
Apply \eqref{eq:badlarge} with $H=\sigma\sqrt n\,h_x$.
Then
\[
 \PP(M_n\ge tV_n,\calN_x^c)
 \le Cn\PP(|\xi|>h_x)\Phibar(t) \le C\varepsilon^{-3}R_n(x)\Phibar(t).
\]
For the absolute sum, use
\[
 \{|S_n^X|\ge tV_n\}
 \subseteq\{M_n\ge tV_n\}
 \cup\left\{\max_{k\le n}(-S_k^X)\ge tV_n\right\}
\]
and apply the same estimate to $X$ and $-X$.
Together with \eqref{eq:badlarge} and \eqref{eq:truncationprob}, these bounds give
\begin{equation}\label{eq:strategytruncation}
 \PP(M_n\ge tV_n,\calN_x^c)+\PP(|S_n^X|\ge tV_n,\calN_x^c)
 \le Cw_n(x),\qquad 2\Phibar(t)(1-p_x)\le Cw_n(x).
\end{equation}
The second term in the decomposition \eqref{eq:002} therefore satisfies
\[
 \bigl|2\Phibar(t)(1-p_x)-\PP(M_n\ge tV_n,\calN_x^c)\bigr|
 \le Cw_n(x).
\]

For the first term in \eqref{eq:002}, we use the conditional law given $\calN_x$.
Let 
\begin{equation}
\label{eq:003}
  \PP_x^*=\PP(\,\cdot\mid\calN_x)
\end{equation}
and write $\E_x^*$
for expectation under this law.
The event $\calN_x$ remains fixed as $t$ varies.
For $t>0$, define
\begin{equation}\label{eq:conditionalcdfs}
 C_x^*(t)=\PP_x^*(M_n<tV_n).
\end{equation}
The first term in the decomposition \eqref{eq:002} equals
\[
 2\Phibar(t)p_x-\PP(M_n\ge tV_n,\calN_x)
 =p_x\{C_x^*(t)-G(t)\}.
\]
The decomposition of $E_n(t)$ becomes
\begin{equation}\label{eq:baseconditioning}
 E_n(t)=p_x\{C_x^*(t)-G(t)\} +2\Phibar(t)(1-p_x)-\PP(M_n\ge tV_n,\calN_x^c).
\end{equation}

For $t>0$, write
\[
 C_{\mathrm{end},x}^*(t)=\PP_x^*(|S_n^X|<tV_n).
\]
Then
\[
 C_x^*(t)-G(t)
 =\{C_x^*(t)-C_{\mathrm{end},x}^*(t)\}
  +\{C_{\mathrm{end},x}^*(t)-G(t)\}.
\]
We first bound $C_{\mathrm{end},x}^*(t)-G(t)$.
Splitting the absolute-sum tail over $\calN_x$ and $\calN_x^c$ gives
\[
 \Cend(t)-G(t)
 =p_x\{C_{\mathrm{end},x}^*(t)-G(t)\} +2\Phibar(t)(1-p_x)-\PP(|S_n^X|\ge tV_n,\calN_x^c).
\]
For $t>0$ with $|t-x|\le(1+x)^{-1}$, this identity,
\eqref{eq:endbound}, \eqref{eq:neighborhood} and \eqref{eq:strategytruncation} give
\begin{equation}\label{eq:conditionalbasebound}
 |C_{\mathrm{end},x}^*(t)-G(t)|
 \le\frac{|\Cend(t)-G(t)|+2\Phibar(t)(1-p_x)}{p_x} +\frac{\PP(|S_n^X|\ge tV_n,\calN_x^c)}{p_x} \le Cw_n(x).
\end{equation}
At $t=x$, \eqref{eq:baseconditioning} and the bounds above yield
\begin{equation}\label{eq:conditionaltarget}
 |E_n(x)|\le p_x|C_x^*(x)-C_{\mathrm{end},x}^*(x)|+Cw_n(x).
\end{equation}

\noindent\textbf{Step 2: representation by a reflected walk.}\enspace
By \eqref{eq:conditionaltarget}, we need to compare
$C_x^*(x)$ and $C_{\mathrm{end},x}^*(x)$.
Under $\PP_x^*$, put $A_j=\xi_j$.
These increments are i.i.d. with common law
\begin{equation}\label{eq:conditionallaw}
 A\stackrel d=\xi\mid\{|\xi|\le h_x\}.
\end{equation}
The conditional mean need not be zero.
Define
\begin{equation}\label{eq:boundedstatistics}
 S_i=\sum_{j=1}^i A_j,\qquad Q_i=\sum_{j=1}^i A_j^2,\qquad S_0=Q_0=0.
\end{equation}
Thus $S_i=S_i^X/(\sigma\sqrt n)$ and $Q_n=V_n^2/(n\sigma^2)$.
We compare $\max(0,S_1,\ldots,S_n)$ with $|S_n|$,
using the common denominator $\sqrt{Q_n}$.

We use the Lindley reflected walk \citep{Lindley1952} to represent the maximum.
Set
\begin{equation}\label{eq:lindley}
 W_0=0,\qquad W_i=(W_{i-1}+A_i)^+,\qquad U_i=(-W_{i-1}-A_i)^+.
\end{equation}
Here $u^+=\max(u,0)$.

\begin{proposition}\label{prop:reflection}
For $x\in\mathcal D_3$, we have $\var(A)>0$ and
\begin{equation}\label{eq:reversal}
 (W_n,Q_n)\stackrel d=(\max(0,S_1,\ldots,S_n),Q_n)
 \quad\text{under }\PP_x^*.
\end{equation}

Let $A'_1,\ldots,A'_n$ be an independent copy of the conditional sample and put
\[
 S'_\ell=\sum_{j=1}^\ell A'_j,\qquad
 Q'_\ell=\sum_{j=1}^\ell (A'_j)^2,\qquad S'_0=Q'_0=0.
\]
Write $\E'_x$ for expectation over this copy only.
For every measurable test $\psi:\R^2\to[0,1]$, define
\begin{equation}\label{eq:futurekernel}
 \begin{split}
 f_\ell^\psi(w,q)&=\E'_x\psi(w+S'_\ell,q+Q'_\ell),\qquad f_0^\psi=\psi,\\
 \Delta_i(\psi)&=f_{n-i}^\psi(0,Q_i)-f_{n-i}^\psi(-U_i,Q_i).
 \end{split}
\end{equation}
Then
\begin{equation}\label{eq:duhamel}
 \E_x^*\psi(W_n,Q_n)-\E_x^*\psi(S_n,Q_n)
 =\sum_{i=1}^n\E_x^*\Delta_i(\psi).
\end{equation}
\end{proposition}
The proof is given in Section~\ref{sec:reduction}.

For $t>0$, \eqref{eq:reversal} gives
\begin{equation}\label{eq:reflectedcdfs}
 C_x^*(t)=\PP_x^*(Q_n>0,W_n<t\sqrt{Q_n}),\qquad C_{\mathrm{end},x}^*(t)=\PP_x^*(Q_n>0,|S_n|<t\sqrt{Q_n}).
\end{equation}
For this threshold, put
\[
 \psi_t(w,q)=\ind_{\{q>0,\,|w|<t\sqrt q\}}.
\]
Since $W_n\ge0$, both probabilities use this same indicator.
Applying \eqref{eq:duhamel}, we obtain
\[
 C_x^*(t)-C_{\mathrm{end},x}^*(t)
 =\E_x^*\psi_t(W_n,Q_n)-\E_x^*\psi_t(S_n,Q_n) =\sum_{i=1}^n\E_x^*\Delta_i(\psi_t).
\]

\noindent\textbf{Step 3: Fourier comparison after smoothing.}\enspace
We now estimate $\sum_{i=1}^n\E_x^*\Delta_i(\psi_t)$.
For the indicator $\psi_t$ above, $f_\ell^{\psi_t}$ may not be smooth enough
for Taylor's formula.
We first smooth the test function and obtain the required derivative bounds
by Fourier methods.
Fix $\chi\in C_c^\infty((1/16,4))$ with $0\le\chi\le1$ and $\chi=1$ on $[1/8,3]$, where $C_c^\infty((1/16,4))$ denotes the class of infinitely differentiable functions support in $(1/16,4)$.
The factor $\chi(q)$ restricts $q$ to a bounded interval away from zero.
For $y\ge0$ and a random vector $(Z,Y)$ independent of the conditional sample, define
\begin{equation}\label{eq:test}
 H_y(w,q)=\chi(q)\ind_{\{|w|<y\sqrt q\}},\qquad
 \widetilde\psi_y(w,q)=\E_{Z,Y}H_y(w+Z,q+Y).
\end{equation}
Here $H_y(w,q)=0$ for $q\le0$, and $\E_{Z,Y}$ is expectation over $(Z,Y)$ only.
Write
\begin{equation}\label{eq:smoothedcdfs}
 \widetilde C_x^*(y)=\E_x^*\widetilde\psi_y(W_n,Q_n),\qquad
 \widetilde C_{\mathrm{end},x}^*(y)=\E_x^*\widetilde\psi_y(S_n,Q_n).
\end{equation}
Since $0\le\widetilde\psi_y\le1$, \eqref{eq:duhamel} gives
\[
 \widetilde C_x^*(y)-\widetilde C_{\mathrm{end},x}^*(y)
 =\sum_{i=1}^n\E_x^*\Delta_i(\widetilde\psi_y),
\]

Proposition~\ref{prop:reflection} gives $\E_x^*(A-\E_x^* A)^2>0$.
To state the estimate, put $\gamma=\E_x^*|A|^3/\E_x^*(A-\E_x^* A)^2$ and define
\begin{equation}\label{eq:Lambda}
 \Lambda_y=(1+y)n|\E_x^* A|+(1+y)^3n\E_x^*|A|^3,\qquad y\ge0.
\end{equation}
The following proposition chooses the noise law and bounds the smoothed difference.

\begin{proposition}\label{prop:smoothed}
For $x\in\mathcal D_3$, there is a random vector $(Z,Y)$
independent of the conditional sample
such that
\begin{equation}\label{eq:smoothing-size}
 |Z|+(1+x)|Y|\le100h_x,\qquad
 \E_{Z,Y}|Z|+(1+x)\E_{Z,Y}|Y|\le100\gamma.
\end{equation}
The first inequality holds almost surely.
The noise law may depend on $x$, but is held fixed as $y$ varies.
For every
\begin{equation}\label{eq:comparisonthresholds}
 y\ge\gamma,\qquad |y-x|\le\frac1{2(1+x)},
\end{equation}
the test in \eqref{eq:test}, using this same noise law, satisfies
\begin{equation}\label{eq:smoothed}
 \left|\E_x^*\widetilde\psi_y(W_n,Q_n)-\E_x^*\widetilde\psi_y(S_n,Q_n)\right| =\left|\sum_{i=1}^n\E_x^*\Delta_i(\widetilde\psi_y)\right|
 \le C\Lambda_y\Phibar(y)\le C_{\varepsilon}w_n(x).
\end{equation}
\end{proposition}
The proof, including the noise construction, is given in Section~\ref{sec:comparison}.

Use the noise law supplied by Proposition~\ref{prop:smoothed}.
For every $y$ satisfying \eqref{eq:comparisonthresholds}, we obtain
\begin{equation}
\label{eq:003}
|\widetilde C_x^*(y)-\widetilde C_{\mathrm{end},x}^*(y)|
 =\left|\sum_{i=1}^n\E_x^*\Delta_i(\widetilde\psi_y)\right|
 \le Cw_n(x).
\end{equation}

\noindent\textbf{Step 4: removing the smoothing and completing the proof.}\enspace
We first compare $\widetilde C_{\mathrm{end},x}^*(y)$ with $G(y)$.
The following two-sided estimate will also let us return from the smoothed maximum to its original tail.
Extend $C_x^*,C_{\mathrm{end},x}^*,G$ and $E_n$ by zero
to nonpositive arguments.

\begin{lemma}\label{lem:sandwich}
Let $x\in\mathcal D_3$.
Use the noise supplied by
Proposition~\ref{prop:smoothed} and put
\begin{equation}\label{eq:delta}
 \delta=4\{|Z|+(1+x)|Y|\},\qquad b_x=\PP_x^*(Q_n\notin[1/4,2]).
\end{equation}
For $y>0$ with $|y-x|\le[2(1+x)]^{-1}$,
\begin{equation}\label{eq:sandwich}
 \E_{Z,Y} C_x^*(y-\delta)-b_x\le\widetilde C_x^*(y)
 \le\E_{Z,Y} C_x^*(y+\delta)+b_x.
\end{equation}
The same inequalities hold for $C_{\mathrm{end},x}^*,\widetilde C_{\mathrm{end},x}^*$.
Moreover,
\begin{equation}\label{eq:sandwichcost}
 \delta\le\frac{400\varepsilon}{1+x},\qquad
 \E_{Z,Y}\delta\le400\gamma\le1600r_n,\qquad
 b_x+\phi(x)\E_{Z,Y}\delta\le C_{\varepsilon}w_n(x).
\end{equation}
\end{lemma}
The proof is given in Supplemental Material~\Sref{supp:sandwich-proof}.

Let $y>0$ satisfy $|y-x|\le[2(1+x)]^{-1}$.
By \eqref{eq:sandwichcost}, $\delta\le400\varepsilon/(1+x)\le L\varepsilon/(1+x)$.
By \eqref{eq:proofparameters}, $\delta\le[4(K+1)(1+x)]^{-1}$.
Hence $|y\pm\delta-x|\le(1+x)^{-1}$.
By \eqref{eq:sandwich} and
\eqref{eq:conditionalbasebound} from Step~1,
\[
 \E_{Z,Y} G(y-\delta)-Cw_n(x)-b_x
 \le\widetilde C_{\mathrm{end},x}^*(y) \le\E_{Z,Y} G(y+\delta)+Cw_n(x)+b_x.
\]
On $(0,\infty)$, $G'=2\phi$.
For nonnegative $t$ with $|t-x|\le(1+x)^{-1}$, we have $|t^2-x^2|\le3$.
Hence
\begin{equation}\label{eq:localweights}
 |G(t)-G(s)|\le C\phi(x)|t-s|,\qquad
 \text{for }|s-x|,|t-x|\le(1+x)^{-1}.
\end{equation}
The zero extension gives the same bound when an argument is nonpositive.
Using \eqref{eq:sandwichcost}, we obtain
\begin{equation}
\label{eq:004}
  |\widetilde C_{\mathrm{end},x}^*(y)-G(y)|
 \le Cw_n(x)+C\phi(x)\E_{Z,Y}\delta+b_x
 \le Cw_n(x).
\end{equation}

Thus, for every $y$ satisfying \eqref{eq:comparisonthresholds},
\eqref{eq:003} and \eqref{eq:004} give
\begin{equation}\label{eq:smoothedcdfbound}
 |\widetilde C_x^*(y)-G(y)|
 \le |\widetilde C_x^*(y)-\widetilde C_{\mathrm{end},x}^*(y)|
      +|\widetilde C_{\mathrm{end},x}^*(y)-G(y)| \le Cw_n(x).
\end{equation}

We now return to the estimate of $E_n(x)$.
By \eqref{eq:conditioningprob}, $p_x\ge3/4$.
By \eqref{eq:baseconditioning} and \eqref{eq:strategytruncation},
\begin{equation}\label{eq:maxpullback}
 C_x^*(t)-G(t)=\frac{E_n(t)}{p_x}
 +O\!\bigl(w_n(x)\bigr),
\end{equation}
uniformly for $|t-x|\le(1+x)^{-1}$, including the zero extensions.
Next, inspired by the recursive approach of \citet{LiuZhang2023}, put
\begin{equation}\label{eq:worsterror}
 K_n=\sup_{t\in\mathcal I_n,\,t>0}\frac{|E_n(t)|}{w_n(t)}.
\end{equation}
By Lemma~\ref{lem:truncatedbudget}, $\mathcal I_n$ is compact.
Since $R_n(0)>0$, $w_n$ has a positive minimum on $\mathcal I_n$.
Thus $K_n<\infty$, because $|E_n(t)|\le1$.
Every positive $t$ in this neighborhood belongs to $\mathcal I_n$ by Step~1.
For $t>0$ with $|t-x|\le(1+x)^{-1}$, Mills' inequalities give
$\phi(u)/\Phibar(u)\le u+4$ for $u\ge0$. Hence
\[
 \frac{\Phibar(t)}{\Phibar(x)}
 \le\exp\{(x-t)_+(x+4)\}\le\ee^4.
\]
Together with \eqref{eq:budgetcomparison}, this yields
$w_n(t)\le8\ee^4w_n(x)<Lw_n(x)$, since $8\ee^4<1600\le L$.
The definition of $K_n$ and the zero extension therefore give
\begin{equation}\label{eq:localerrorbound}
 |E_n(t)|\le LK_nw_n(x).
\end{equation}

We split $\delta$ according to whether it exceeds a fixed multiple of its mean.
With $K$ fixed in \eqref{eq:proofparameters}, put
\[
 a=K\E_{Z,Y}\delta.
\]
If $\E_{Z,Y}\delta>0$, Markov's inequality gives
$\PP_{Z,Y}(\delta>a)\le\E_{Z,Y}\delta/a=1/K$.
If $\E_{Z,Y}\delta=0$, then $\delta=0$ almost surely and $\PP_{Z,Y}(\delta>a)=0$.

We first check that the shifts used below remain in the required neighborhood of $x$.
By \eqref{eq:sandwichcost} and $L\ge1600$,
\[
 0\le\delta\le\frac{400\varepsilon}{1+x}
 \le\frac{L\varepsilon}{1+x}
 \qquad\text{almost surely}.
\]
Taking expectation over the noise and using $a=K\E_{Z,Y}\delta$, we obtain
\[
 \E_{Z,Y}\delta\le\frac{L\varepsilon}{1+x},
 \qquad a\le\frac{KL\varepsilon}{1+x}.
\]
Since \eqref{eq:proofparameters} gives
$\varepsilon\le[4(K+1)L]^{-1}$, it follows that
\begin{equation}\label{eq:shiftbudget}
 a+\delta\le\frac{KL\varepsilon}{1+x}
             +\frac{L\varepsilon}{1+x}
 =\frac{(K+1)L\varepsilon}{1+x}
 \le\frac1{4(1+x)}
 \qquad\text{almost surely}.
\end{equation}
Also, \eqref{eq:sandwichcost} gives
$\E_{Z,Y}\delta\le1600r_n\le Lr_n$ and $\gamma\le4r_n\le Lr_n$.
Since $x\in\mathcal D_3$ and $A_0=4KL$, we have $x>4KLr_n$. Hence
\begin{equation}\label{eq:shiftlowerbounds}
 a\le KLr_n<x/4,\qquad
 \gamma\le Lr_n<x/(4K),\qquad
 x-a\ge x/2\ge\gamma.
\end{equation}

Choose $s\in\{-1,1\}$ so that $sE_n(x)=|E_n(x)|$, and set
\[
 y=x+sa,\qquad t_\delta=y-s\delta=x+s(a-\delta).
\]
By \eqref{eq:shiftbudget} and \eqref{eq:shiftlowerbounds},
\[
 y\ge x-a\ge\gamma,\qquad
 |y-x|=a\le\frac1{4(1+x)}\le\frac1{2(1+x)},
\]
and, almost surely,
\[
 |t_\delta-x|=|a-\delta|\le a+\delta
 \le\frac1{4(1+x)}\le\frac1{2(1+x)}.
\]
Thus \eqref{eq:sandwich} and \eqref{eq:smoothedcdfbound} apply at $y$.
For this fixed $x$, the quantities $s,a,y$ are deterministic.
If $s=1$, then $t_\delta=y-\delta$; the lower bound in
\eqref{eq:sandwich}, after subtracting $G(y)$, gives
\[
 \widetilde C_x^*(y)-G(y)
 \ge\E_{Z,Y}\bigl[C_x^*(t_\delta)-G(y)\bigr]-b_x.
\]
If $s=-1$, then $t_\delta=y+\delta$. The upper bound in
\eqref{eq:sandwich} gives
$\widetilde C_x^*(y)\le\E_{Z,Y}C_x^*(t_\delta)+b_x$.
Subtracting $G(y)$ and multiplying by $-1$ therefore yields
\[
 -\{\widetilde C_x^*(y)-G(y)\}
 \ge\E_{Z,Y}\bigl[-\{C_x^*(t_\delta)-G(y)\}\bigr]-b_x.
\]
In both cases,
\[
 s\{\widetilde C_x^*(y)-G(y)\}
 \ge\E_{Z,Y}\bigl[s\{C_x^*(t_\delta)-G(y)\}\bigr]-b_x.
\]
On the other hand, \eqref{eq:smoothedcdfbound} gives
\[
 s\{\widetilde C_x^*(y)-G(y)\}
 \le\bigl|\widetilde C_x^*(y)-G(y)\bigr|\le Cw_n(x).
\]
Combining these inequalities, we obtain
\begin{equation}\label{eq:signedexpectationbound}
 Cw_n(x)\ge s\{\widetilde C_x^*(y)-G(y)\}
 \ge\E_{Z,Y}\bigl[s\{C_x^*(t_\delta)-G(y)\}\bigr]-b_x.
\end{equation}
Since $|s|=1$, the $O(w_n(x))$ term in \eqref{eq:maxpullback},
after multiplication by $s$, is bounded below by $-Cw_n(x)$.
We use this observation on each of the two events.

On $\{\delta\le a\}$, $s(t_\delta-x)=a-\delta\ge0$.
Since $C_x^*$ is nondecreasing, $sC_x^*(t_\delta)\ge sC_x^*(x)$.
Applying \eqref{eq:maxpullback} at $x$ and \eqref{eq:localweights}
with $|x-y|=a$, and using $sE_n(x)=|E_n(x)|$, we obtain
\[
 \begin{aligned}
 s\{C_x^*(t_\delta)-G(y)\}
 &\ge s\{C_x^*(x)-G(y)\}\\
 &=s\{C_x^*(x)-G(x)\}+s\{G(x)-G(y)\}\\
 &\ge\frac{sE_n(x)}{p_x}-Cw_n(x)-|G(x)-G(y)|\\
 &\ge\frac{|E_n(x)|}{p_x}-Cw_n(x)-C\phi(x)a.
 \end{aligned}
\]

On $\{\delta>a\}$, \eqref{eq:shiftbudget} gives
$|t_\delta-x|\le a+\delta\le(1+x)^{-1}$.
Thus \eqref{eq:localerrorbound} applies at $t=t_\delta$ and gives
$|E_n(t_\delta)|\le LK_nw_n(x)$.
We now add and subtract $G(t_\delta)$.
By \eqref{eq:maxpullback} at $t_\delta$ and \eqref{eq:localweights}
with $|t_\delta-y|=\delta$,
\[
 \begin{aligned}
 s\{C_x^*(t_\delta)-G(y)\}
 &=s\{C_x^*(t_\delta)-G(t_\delta)\}
   +s\{G(t_\delta)-G(y)\}\\
 &\ge\frac{sE_n(t_\delta)}{p_x}-Cw_n(x)
       -|G(t_\delta)-G(y)|\\
 &\ge-\frac{|E_n(t_\delta)|}{p_x}-Cw_n(x)-C\phi(x)\delta\\
 &\ge-\frac{LK_nw_n(x)}{p_x}-Cw_n(x)-C\phi(x)\delta.
 \end{aligned}
\]
The zero extensions ensure that the same argument applies when $t_\delta\le0$.
Combining the two cases, we have
\begin{equation}\label{eq:twocaselowerbound}
 s\{C_x^*(t_\delta)-G(y)\}\ge
 \begin{cases}
 |E_n(x)|/p_x-Cw_n(x)-C\phi(x)a,&\delta\le a,\\
 -LK_nw_n(x)/p_x-Cw_n(x)-C\phi(x)\delta,&\delta>a.
 \end{cases}
\end{equation}
Taking expectation over the two events therefore gives
\begin{equation}\label{eq:averagedlowerbound}
 \begin{aligned}
 &\E_{Z,Y}\bigl[s\{C_x^*(t_\delta)-G(y)\}\bigr]\\
 &\quad\ge
 \frac{\PP_{Z,Y}(\delta\le a)|E_n(x)|
       -\PP_{Z,Y}(\delta>a)LK_nw_n(x)}{p_x}
 -Cw_n(x)\\
 &\qquad\quad
 -C\phi(x)\left\{
 a\PP_{Z,Y}(\delta\le a)
 +\E_{Z,Y}\bigl[\delta\ind_{\{\delta>a\}}\bigr]
 \right\}.
 \end{aligned}
\end{equation}
Since $a\ge0$ and $\delta\ge0$, we also have
\[
 a\PP_{Z,Y}(\delta\le a)\le a,
 \qquad
 \E_{Z,Y}\bigl[\delta\ind_{\{\delta>a\}}\bigr]
 \le\E_{Z,Y}\delta.
\]
Combining \eqref{eq:averagedlowerbound} with \eqref{eq:signedexpectationbound}, we obtain
\begin{equation}\label{eq:averagederrorbound}
 \frac{\PP_{Z,Y}(\delta\le a)|E_n(x)|
       -\PP_{Z,Y}(\delta>a)LK_nw_n(x)}{p_x} \le Cw_n(x)+C\phi(x)(a+\E_{Z,Y}\delta)+b_x.
\end{equation}
 \eqref{eq:sandwichcost} gives
\begin{equation}\label{eq:averagingcostbound}
 C\phi(x)(a+\E_{Z,Y}\delta)+b_x
 =C(K+1)\phi(x)\E_{Z,Y}\delta+b_x \le Cw_n(x).
\end{equation}
Using \eqref{eq:averagederrorbound} and \eqref{eq:averagingcostbound},
together with $\PP_{Z,Y}(\delta\le a)\ge1-1/K$,
$\PP_{Z,Y}(\delta>a)\le1/K$ and $p_x\le1$, we get
\begin{equation}\label{eq:contraction}
 \frac{|E_n(x)|}{w_n(x)}
 \le C+\frac{L}{K-1}K_n
 \le C+\tfrac14K_n.
\end{equation}
This holds for every $x\in\mathcal D_3$.
For positive $x\in\mathcal D_1\cup\mathcal D_2$,
Step~1 gives $|E_n(x)|/w_n(x)\le C$.
Taking the supremum over positive $x\in\mathcal I_n$, we get
\[
 K_n\le C+K_n/4.
\]
Since $K_n<\infty$, it follows that $K_n\le4C/3$.
Thus $|E_n(x)|\le Cw_n(x)$ for all positive $x\in\mathcal I_n$.
For $0<x\le\beta_{n,p}^{-1/p}$, \eqref{eq:complement} gives
\[
 \left|\frac{\PP(M_n\ge xV_n)}{\Phibar(x)}-2\right|
 =\frac{|E_n(x)|}{\Phibar(x)}
 \le C\frac{w_n(x)}{\Phibar(x)}=CR_n(x),\qquad 0<x\le\beta_{n,p}^{-1/p}.\qedhere
\]
\end{proof}

\subsection{Proof of Theorem~\ref{thm:sharp}}\label{sec:sharpness}

\begin{proof}
\textbf{Case 1: $2<p<3$.}
Fix $\nu=\lceil8/(3-p)\rceil$.
For $k\ge1$, set
\begin{equation}\label{eq:sharp-parameters}
 b_k=\frac{2^{\nu^k+1}}{\sqrt{2\nu^k\log2}},\qquad
 \pi_k=2^{-k-2}b_k^{-p},
\end{equation}
and define
\begin{equation}\label{eq:sharp-law}
 P_*\{1\}=P_*\{-1\}=\frac{1-\sum_{k\ge1}\pi_k}{2},
 \qquad P_*\{b_k\}=P_*\{-b_k\}=\frac{\pi_k}{2}.
\end{equation}
Since $b_k\uparrow\infty$ and $b_1>1$, we have $\sum_k\pi_k<\sum_k2^{-k-2}=1/4$.
Thus \eqref{eq:sharp-law} defines a symmetric probability law with countably infinite support.
Write $\sigma^2=\E_{P_*}X^2$.
Then
\begin{equation}\label{eq:sharp-moments}
 1<\sigma^2\le\E_{P_*}|X|^p
 =1-\sum_k\pi_k+\sum_k2^{-k-2}
 =\tfrac54-\sum_k\pi_k<\tfrac54.
\end{equation}
For every $\eta>0$, we have $\pi_kb_k^{p+\eta}=2^{-k-2}b_k^\eta\to\infty$,
since $\log b_k\asymp\nu^k$.
Hence $\E_{P_*}|X|^{p+\eta}=\infty$.
Take
\begin{equation}\label{eq:sharp-subsequence}
 N_k=2^{2\nu^k},\qquad x_k=\sqrt{\log N_k}.
\end{equation}
Then $b_k=2\sqrt{N_k}/x_k$, and the thresholds satisfy
\[
 x_k\beta_{N_k,p}(P_*)^{1/p}
 =\frac{(\E_{P_*}|X|^p)^{1/p}}{\sigma}
   \frac{\sqrt{\log N_k}}{N_k^{(p-2)/(2p)}}\longrightarrow0.
\]

Put $a=2/\sigma$ and $\delta_k=N_k\pi_k$.
Here $\delta_k$ is the expected number of observations of magnitude $b_k$ in a sample of size $N_k$ from $P_*$.
The following lemma gives the leading terms of $R_{N_k}(x_k;P_*)$ and the relative tail error.
\begin{lemma}\label{lem:conditional-tails}
For the distribution \eqref{eq:sharp-law} and sequences
\eqref{eq:sharp-parameters} and \eqref{eq:sharp-subsequence},
put $\vartheta_p=\ee^{-a^2/2}\cosh(a)$.
Then $\vartheta_p<1$, and
\begin{align}
 R_{N_k}(x_k;P_*)&\sim a^2\delta_k\longrightarrow0,\label{eq:sharp-rn-scale}\\
 \frac{\PP_{P_*}(M_{N_k}\ge x_kV_{N_k})}{\Phibar(x_k)}-2
 &=2(\vartheta_p-1)\delta_k+o(\delta_k).\label{eq:sharp-tail-expansion}
\end{align}
\end{lemma}
The proof is given in Supplemental Material~\Sref{supp:conditional-tails-proof}.

By the definition of $\mathcal E_n$, \eqref{eq:sharp-tail-expansion} and \eqref{eq:sharp-rn-scale} give
\begin{equation}\label{eq:sharp-scale}
 \frac{\mathcal E_{N_k}(x_k;P_*)}{R_{N_k}(x_k;P_*)}
 \longrightarrow\kappa_p:=\frac{2(1-\vartheta_p)}{a^2}>0.
\end{equation}

\noindent\textbf{Case 2: $p=3$.}
The preceding choice $\nu=\lceil8/(3-p)\rceil$ is undefined at $p=3$.
We use a similar discrete distribution and take the threshold to be zero.
For $j\ge1$, set
\begin{equation}\label{eq:sharp-third-law}
 \begin{split}
 b_j^{(3)}&=2^{j^2},\qquad
 \pi_j^{(3)}=2^{-j-2}(b_j^{(3)})^{-3},\\
 P_*\{1\}=P_*\{-1\}&=\frac{1-\sum_{j\ge1}\pi_j^{(3)}}2,\qquad
 P_*\{b_j^{(3)}\}=P_*\{-b_j^{(3)}\}=\frac{\pi_j^{(3)}}2.
 \end{split}
\end{equation}
This defines a symmetric law with countably infinite support.
Write $\sigma^2(P_*)=\E_{P_*}X^2$ and $\mu_3(P_*)=\E_{P_*}|X|^3$.
Then
\[
 1<\sigma^2(P_*)\le\mu_3(P_*)
 =\tfrac54-\sum_{j\ge1}\pi_j^{(3)}<\tfrac54.
\]
For every $\eta>0$, the terms
$\pi_j^{(3)}(b_j^{(3)})^{3+\eta}=2^{-j-2+\eta j^2}$ tend to infinity.
Hence $\E_{P_*}|X|^{3+\eta}=\infty$.

Take $N_k=k$ and $x_k=0$.
At threshold zero,
\[
 \mathcal E_n(0;P_*)
 =2\PP_{P_*}(M_n<0)
 =2\PP_{P_*}(S_1^X<0,\ldots,S_n^X<0).
\]

The next lemma estimates this probability and $R_n(0;P_*)$.

\begin{lemma}\label{lem:strict-survival}
Let $X_1,X_2,\ldots$ be i.i.d.\ with a symmetric distribution $P$.
Assume $\sigma^2=\E_PX^2\in(0,\infty)$ and $\mu_3=\E_P|X|^3<\infty$.
For $S_j^X=\sum_{i=1}^jX_i$, write
$q_n=\PP_P(S_1^X<0,\ldots,S_n^X<0)$ and $a_j=\PP_P(S_j^X=0)$.
Then
\begin{equation}\label{eq:strict-survival-asymptotic}
 \begin{gathered}
 A(P):=\exp\left\{-\frac12\sum_{j\ge1}\frac{a_j}{j}\right\}>0,\\
 q_n\sim\frac{A(P)}{\sqrt{\pi n}},\qquad
 R_n(0;P)\sim\frac{\mu_3}{\sigma^3\sqrt n}\quad(n\to\infty).
 \end{gathered}
\end{equation}
\end{lemma}
The proof is given in Supplemental Material~\Sref{supp:strict-survival-proof}.

Applying Lemma~\ref{lem:strict-survival} with $P=P_*$ gives
\begin{equation}\label{eq:sharp-third-limit}
 \frac{\mathcal E_n(0;P_*)}{R_n(0;P_*)}
 \longrightarrow\kappa_3:=\frac{2\sigma^3(P_*)A(P_*)}{\sqrt\pi\,\mu_3(P_*)}>0.
\end{equation}
The same lemma gives $R_{N_k}(x_k;P_*)\to0$.
Since $x_k=0$, we also have $x_k\beta_{N_k,3}(P_*)^{1/3}=0$.
\end{proof}
\section{Proof of Proposition~\ref{prop:preliminary}}\label{sec:preliminary}

We use the assumptions and notation of Section~\ref{sec:proofpropositions}.
Let $c>0$ be the small absolute constant in
\citet[Theorem~2.1]{ShaoZhou2016}, and fix
$\rho_*=\min\{1/16,c\}$.
The proof of Lemma~\ref{lem:truncatedbudget} does not depend on Proposition~\ref{prop:preliminary}.
We therefore use its notation and results in the following proof.

\noindent\emph{Proof of \eqref{eq:coarse}.}
Recall from Lemma~\ref{lem:truncatedbudget} that $F(z)=z^2\min\{|z|,1\}$.
For $x\in\mathcal I_n$, $x>0$, we first proof that $x^6\leq n$.
If $x\le\beta_{n,p}^{-1/p}$, then $F(z)\le|z|^p$ gives
$n\E F(x\xi)\le\beta_{n,p}x^p\le1$.
Lyapunov's inequality gives
$\beta_{n,p}^{-1/p}\le n^{(p-2)/(2p)}\le n^{1/6}$,
so $n\ge x^6$ in this case.
If $R_n(x)\le\rho_*$, then Lemma~\ref{lem:truncatedbudget} gives
$n\E F(x\xi)\le R_n(x)\le1$,
and \eqref{eq:budgetlower} gives $n\ge(1+x)^6$.
Thus, in either case,
\begin{equation}\label{supp:xn}
 n\E F(x\xi)\le1,\qquad n\ge x^6.
\end{equation}
For $0<x<8$, \eqref{eq:coarse} follows from
$\PP(M_n\ge xV_n)\le1$ and $\Phibar(x)\ge\Phibar(8)$.
It therefore suffices to consider $x\ge8$.

To continue our proof, we first remove large increments. The next two lemmas control the retained sample and compare its crossing event with the original one.

\begin{lemma}\label{supp:boundedrandom}
Fix $D\ge1$. Let $x\ge1$, $k\ge1$, and let
$\widetilde\xi,\widetilde\xi_1,\ldots,\widetilde\xi_k$
be i.i.d.\ real random variables satisfying
\[
 |\widetilde\xi|\le x^{-1},\qquad \nu:=k\E\widetilde\xi^2\in[1/2,2],\qquad k|\E\widetilde\xi|\le D/x,\qquad k\E|\widetilde\xi|^3\le D/x^3.
\]
The support bound holds almost surely.
For this sample, define
$Q=\sum_{i=1}^k\widetilde\xi_i^2$ and
$M=\max_{1\le j\le k}\sum_{i=1}^j\widetilde\xi_i$.
Then
\begin{equation}\label{supp:boundedrandombound}
 \PP(M\ge x\sqrt Q)\le C_D\Phibar(x),
\end{equation}
where $C_D$ depends only on $D$. The weak event includes $Q=0$.
\end{lemma}
The proof is given in Supplemental Material~\Sref{supp:boundedrandomproof}.

The following inequality generalizes the one-observation case in
\citet[Section~2, equation~(2.6)]{LiuShaoWang2013}.

\begin{lemma}\label{supp:deletion}
Let $J\subset\{1,\ldots,n\}$ with $|J|=j<n$, and write
$\{1,\ldots,n\}\setminus J=\{i_1<\cdots<i_{n-j}\}$.
Define
\[
 M^{(-J)}:=\max_{1\le k\le n-j}\sum_{r=1}^kX_{i_r},
 \qquad
 V^{(-J)}:=\left(\sum_{r=1}^{n-j}X_{i_r}^2\right)^{1/2}.
\]
If $x>0$ and $x^2>j$, then
\begin{equation}\label{eq:deletion}
 \{M_n\ge xV_n\}\subseteq
 \{M^{(-J)}\ge\sqrt{x^2-j}\,V^{(-J)}\}.
\end{equation}
Both events include samples whose respective square sums vanish.
\end{lemma}
The proof is given in Supplemental Material~\Sref{supp:deletionproof}.

Delete the observations with $|\xi_i|>1/x$, and define the random index set
$J=\{i\in\{1,\ldots,n\}:|\xi_i|>1/x\}$ and its cardinality $K=|J|$.
The variable $K$ is binomial, and \eqref{supp:xn} and the Markov inequality give
\[
 \E K=n\PP(|\xi|>1/x)\le n\E F(x\xi)\le1.
\]
Since $K$ is binomial with mean at most one,
\[
 \E\ee^{K/2}\le\exp\{\ee^{1/2}-1\},\qquad
 \PP(K>x^2/2)\le(2\ee/x^2)^{x^2/2}\le\ee^{-x^2}
\]
for $x\ge8$.

We split according to whether $K\le x^2/2$.
To prove \eqref{eq:coarse}, it therefore suffices to show that,
for every deterministic set $I\subseteq\{1,\ldots,n\}$ with
$\PP(J=I)>0$ and $j:=|I|\le x^2/2$,
\begin{equation}\label{supp:conditionalcoarse}
 \PP(M_n\ge xV_n\mid J=I)\le C\ee^{j/2}\Phibar(x),
\end{equation}
where $C$ is absolute.
Indeed, averaging this bound over $J$ would give
\[
 \PP(M_n\ge xV_n)
 \le C\Phibar(x)\E\ee^{K/2}+\PP(K>x^2/2) \le C\Phibar(x)+\ee^{-x^2} \le C\Phibar(x).
\]

It remains to prove \eqref{supp:conditionalcoarse}.
Fix such a set $I$, and put $k=n-j$ and $y=\sqrt{x^2-j}$.
Then $k\ge1$ and $x/\sqrt2\le y\le x$.
Conditional on the event $\{J=I\}$, write the remaining increments,
in their original order, as $\widetilde\xi_1,\ldots,\widetilde\xi_k$, 
and they are i.i.d.\ with common law
\[
 \widetilde\xi\stackrel d=\xi\mid\{|\xi|\le1/x\}.
\]
By \eqref{eq:deletion}, on the event $\{J=I\}$,
the inequality $M_n\ge xV_n$ implies
$M^{(-I)}\ge yV^{(-I)}$ for the retained sample.
Since $y\le x$,
\[
 |\widetilde\xi|\le1/x\le1/y.
\]
To check the moment conditions in Lemma \ref{supp:boundedrandom}, put $q_x=\PP(|\xi|\le1/x)$.
By \eqref{supp:xn}, $q_x\ge1-1/n\ge1/2$ and
\[
 n\E[|\xi|^r\mathbf1_{\{|\xi|>1/x\}}]\le x^{-r}\quad(r=1,2),
 \qquad n\E[|\xi|^3\mathbf1_{\{|\xi|\le1/x\}}]\le x^{-3}.
\]
Using $\E\xi=0$ and $n\E\xi^2=1$, and dividing the retained moments by $q_x$, gives
\[
 n|\E\widetilde\xi|\le2/x,\qquad
 1-x^{-2}\le n\E\widetilde\xi^2\le1+2/n,\qquad
 n\E|\widetilde\xi|^3\le2/x^3.
\]
Since $k/n\ge1-1/(2x^4)$ and $x\ge2$, it follows that
\[
 \frac12\le k\E\widetilde\xi^2\le2,\qquad
 k|\E\widetilde\xi|\le2/y,\qquad
 k\E|\widetilde\xi|^3\le2/y^3.
\]
Thus Lemma~\ref{supp:boundedrandom} applies at sample size $k$
and threshold $y$, with $D=2$, and yields
\[
 \PP(M_n\ge xV_n\mid J=I)
 \le C\Phibar(y)\le C\ee^{j/2}\Phibar(x).
\]
The last inequality follows from Mills' inequalities and $y\ge x/\sqrt2$.

This proves \eqref{supp:conditionalcoarse}, and hence \eqref{eq:coarse}.

\noindent\emph{Proof of \eqref{eq:endbound}.}
Let $x\in\mathcal I_n$, $x>0$.
For $R_n(x)\le\rho_*$, applying
\citet[Theorems~2.1--2.2]{ShaoZhou2016} and accounting for $V_n=0$ gives \footnote{The verification of these hypotheses in these theorems and the treatment of the case $V_n=0$ (i.e. \eqref{supp:jsw} and \eqref{eq:zeroatom})
involve lengthy technical details.
We therefore give the details in Supplemental Material~\Sref{supp:reference-verification}.}
\begin{equation}\label{supp:jsw}
 |\PP(S_n^X\ge xV_n)-\Phibar(x)|\le Cw_n(x),\qquad x>0,\quad R_n(x)\le\rho_*.
\end{equation}
In particular, H\"older's inequality and Mills' lower bound give
\begin{equation}\label{eq:zeroatom}
 \PP(V_n=0)\le Cw_n(x),\qquad x\ge0,\quad R_n(x)\le\rho_*.
\end{equation}
Apply \eqref{supp:jsw} to $X$ and $-X$.
The two events intersect exactly when $V_n=0$, so
\[
 \Cend(x)=1-\PP(S_n^X\ge xV_n)-\PP(-S_n^X\ge xV_n)+\PP(V_n=0).
\]
Together with \eqref{eq:zeroatom}, this proves \eqref{eq:endbound} when $R_n(x)\le\rho_*$.

For $R_n(x)>\rho_*$, apply \eqref{eq:coarse} to $X$ and $-X$.
This gives $\PP(|S_n^X|\ge xV_n)\le2C\Phibar(x)$, and hence
\[
 |\Cend(x)-G(x)|
 \le(2+2C)\Phibar(x)
 \le(2+2C)\rho_*^{-1}w_n(x).
\]
Thus \eqref{eq:endbound} holds for every positive $x\in\mathcal I_n$.

\noindent\emph{Proof of \eqref{eq:badlarge}.}
Fix $H>0$ and $x\in\mathcal I_n$, $x>0$.
First suppose $x\ge2$.
By \eqref{supp:xn}, $n>1$.
Put $y=\sqrt{x^2-1}$.
We check that $y\in\mathcal I_{n-1}$ for the same increment law.
If $x\le\beta_{n,p}^{-1/p}$, write
\[
 \beta_{m,p}^{-2/p}=(\sigma^p/\mu_p)^{2/p}m^{(p-2)/p},
 \qquad(\sigma^p/\mu_p)^{2/p}\le1.
\]
The derivative of the right side, as a function of $m\ge1$, is at most one.
Hence
\[
 \beta_{n,p}^{-2/p}-\beta_{n-1,p}^{-2/p}\le1,\qquad
 y\le\beta_{n-1,p}^{-1/p}.
\]
If $R_n(x)\le\rho_*$, \eqref{eq:budgetlower} gives $n\ge(1+x)^2$.
Since $(1+x)^2-(1+y)^2=1+2(x-y)\ge1$,
\[
 \frac{1+y}{\sqrt{n-1}}\le\frac{1+x}{\sqrt n}.
\]
Monotonicity of $a\mapsto F(aX)$ for $a\ge0$ now gives
\[
 R_{n-1}(y)
 =(n-1)\E F\left(\frac{(1+y)X}{\sigma\sqrt{n-1}}\right)
 \le\frac{n-1}{n}R_n(x)\le\rho_*.
\]
Thus $y\in\mathcal I_{n-1}$ in both cases, and \eqref{eq:coarse} applies to the retained sample.

For each index $i$, apply \eqref{eq:deletion} with $J=\{i\}$.
The event on the right depends only on the remaining observations,
so it is independent of $X_i$.
Using \eqref{eq:coarse} at sample size $n-1$, we obtain
\[
 \begin{aligned}
 \PP(M_n\ge xV_n,|X_i|>H)
 &\le\PP(|X|>H)\PP(M_{n-1}\ge yV_{n-1})\\
 &\le C\PP(|X|>H)\Phibar(y)
 \le C\PP(|X|>H)\Phibar(x).
 \end{aligned}
\]
Here $(M_{n-1},V_{n-1})$ denotes the statistics of the remaining i.i.d.\ sample.
The last inequality follows from Mills' inequalities.
A union bound over $i$ gives
\[
 \PP(M_n\ge xV_n,\max_{i\le n}|X_i|>H)
 \le Cn\PP(|X|>H)\Phibar(x).
\]
For $0<x<2$, the same bound follows from the union bound
and $\Phibar(x)\ge\Phibar(2)>0$.
This proves \eqref{eq:badlarge}.

\noindent\emph{Proof of \eqref{eq:nearzero}.}
We use the following bounds for the maximum and the sum of squares.

\begin{lemma}\label{supp:origin}
Let $X,X_1,X_2,\ldots$ be i.i.d.\ with $\E X=0$, $\sigma^2=\E X^2\in(0,\infty)$ and $\mu_p=\E|X|^p<\infty$ for some $2<p\le3$. Write $S_k^X=\sum_{j=1}^kX_j$. For every integer $m\ge1$ and every $a\ge0$,
\begin{equation}\label{supp:lowsurvivalbound}
 \PP\{\max(0,S_1^X,\ldots,S_m^X)\le a\}
 \le C\min\left\{1,\frac{a}{\sigma\sqrt m}+r_m\right\},
\end{equation}
where $r_m=R_m(0)$ is computed from the same increment law at sample size $m$. Also, with $V_n^2=\sum_{j=1}^nX_j^2$ and $r_n=R_n(0)$, every $n\ge1$ satisfies
\begin{equation}\label{supp:squaresbound}
 \PP(|V_n^2-n\sigma^2|>n\sigma^2/2)\le Cr_n.
\end{equation}
The constant $C$ is absolute and uniform in $p\in(2,3]$.
\end{lemma}
The proof is given in Supplemental Material~\Sref{supp:originproof}.

For $0<x\le1$, split according to whether $V_n^2\le2n\sigma^2$
and apply Lemma~\ref{supp:origin}:
\[
 1-\PP(M_n\ge xV_n)
 \le\PP\{\max(0,M_n)\le x\sigma\sqrt{2n}\}
       +\PP(V_n^2>2n\sigma^2) \le C(x+r_n).
\]
Since $1-2\Phibar(x)\le Cx$, this implies
\[
 |\PP(M_n\ge xV_n)-2\Phibar(x)|
 \le \{1-\PP(M_n\ge xV_n)\}+\{1-2\Phibar(x)\}
 \le C(x+r_n).
\]
At $x=0$, the same lemma with $a=0$ gives
\[
 1-\PP(M_n\ge0)=\PP(S_1^X<0,\ldots,S_n^X<0)
 \le\PP\{\max(0,M_n)=0\}\le Cr_n.
\]
This proves \eqref{eq:nearzero} and completes Proposition~\ref{prop:preliminary}.
All constants are absolute and uniform in $p\in(2,3]$.

\section{Proof of Proposition~\ref{prop:reflection}}\label{sec:reduction}

\begin{proof}
Fix $x\in\mathcal D_3$.
Put
\[
 q=\PP(|\xi|>h_x),\qquad
 v=n\E[\xi^2\ind_{\{|\xi|>h_x\}}].
\]
The definition of $R_n(x)$ and $h_x=\varepsilon/(1+x)$ give
\[
 nq\le\varepsilon^{-3}R_n(x)\le\tfrac14,\qquad
 v\le\frac{R_n(x)}{\varepsilon(1+x)^2}\le\tfrac14.
\]
Indeed, on $\{|\xi|>h_x\}$ we have
$(1+x)^2|\xi|^2\min\{(1+x)|\xi|,1\}\ge\varepsilon^3$, and 
the bounds by $1/4$ follow from
$R_n(x)<\eta(\varepsilon)\le\varepsilon^3/8$.
Since $\E\xi=0$ and $n\E\xi^2=1$,
\eqref{eq:conditionallaw} gives
\[
 \E_x^* A=-\frac{\E[\xi\ind_{\{|\xi|>h_x\}}]}{1-q},\qquad
 n\E_x^* A^2=\frac{1-v}{1-q}\ge\tfrac34.
\]
By the Cauchy--Schwarz inequality,
\[
 n(\E_x^* A)^2\le\frac{qv}{(1-q)^2}\le v\le\tfrac14.
\]
It follows that
\[
 n\E_x^*(A-\E_x^* A)^2=n\E_x^* A^2-n(\E_x^* A)^2\ge\tfrac12>0.
\]

The Lindley recursion \citep{Lindley1952} starts with
\[
 W_1=\max(0,A_1),\qquad
 W_2=\max\{0,\max(0,A_1)+A_2\}
     =\max(0,A_2,A_1+A_2).
\]
More generally, with empty sums interpreted as zero, we have
\[
 W_i=\max_{0\le k\le i}\sum_{j=k+1}^i A_j,
 \qquad 0\le i\le n.
\]
Indeed, this identity holds at $i=0$ since $W_0=0$.
If it holds at $i-1$, then
\[
 \begin{aligned}
 W_i
 &=\max\{0,W_{i-1}+A_i\}\
 =\max\left\{0,\max_{0\le k\le i-1}
                 \sum_{j=k+1}^i A_j\right\}
 =\max_{0\le k\le i}\sum_{j=k+1}^i A_j.
 \end{aligned}
\]
The last equality includes the zero term by taking $k=i$.
Thus induction gives, at $i=n$,
\[
 W_n=\max(0,A_n,A_{n-1}+A_n,\ldots,A_1+\cdots+A_n).
\]
Under $\PP_x^*$, the increments $A_1,\ldots,A_n$ are i.i.d.
Reversing their order preserves their joint law
and leaves $Q_n=\sum_{j=1}^n A_j^2$ unchanged.
The sums in the displayed maximum then become $S_1,\ldots,S_n$.
Therefore
\[
 (W_n,Q_n)\stackrel d=(\max(0,S_1,\ldots,S_n),Q_n)
 \quad\text{under }\PP_x^*.
\]

For the telescoping identity, fix a measurable test $\psi:\R^2\to[0,1]$.
Put $\calF_i=\sigma(A_1,\ldots,A_i)$.
By the definition of $f_\ell^\psi$ and independence of the increments,
\[
 \E_x^*\!\left[
 f_{n-i}^\psi(W_{i-1}+A_i,Q_i)\mid\calF_{i-1}\right]
 =f_{n-i+1}^\psi(W_{i-1},Q_{i-1}),
\]
where $Q_i=Q_{i-1}+A_i^2$.
We next compare the positions before and after reflection at step $i$.
If $W_{i-1}+A_i\ge0$, then $W_i=W_{i-1}+A_i$ and $U_i=0$.
Consequently,
\[
 f_{n-i}^\psi(W_i,Q_i)
 -f_{n-i}^\psi(W_{i-1}+A_i,Q_i)=0,
\]
and the definition in \eqref{eq:futurekernel} also gives
\[
 \Delta_i(\psi)
 =f_{n-i}^\psi(0,Q_i)-f_{n-i}^\psi(0,Q_i)=0.
\]
If $W_{i-1}+A_i<0$, then
\[
 W_i=0,\qquad U_i=-(W_{i-1}+A_i)>0.
\]
The position before reflection is therefore $-U_i$, and
\[
 f_{n-i}^\psi(W_i,Q_i)
  -f_{n-i}^\psi(W_{i-1}+A_i,Q_i) =f_{n-i}^\psi(0,Q_i)-f_{n-i}^\psi(-U_i,Q_i)
 =\Delta_i(\psi).
\]
In both cases the square-sum coordinate remains
$Q_i=Q_{i-1}+A_i^2$.
Thus, 
\[
 f_{n-i}^\psi(W_i,Q_i)-f_{n-i}^\psi(W_{i-1}+A_i,Q_i)=\Delta_i(\psi).
\]
Combining the two identities yields
\begin{equation}\label{eq:reflection-step}
 \E_x^*\!\left[
 f_{n-i}^\psi(W_i,Q_i)-f_{n-i+1}^\psi(W_{i-1},Q_{i-1})
 \mid\calF_{i-1}\right] =\E_x^*[\Delta_i(\psi)\mid\calF_{i-1}].
\end{equation}
Take expectations and sum over $i=1,\ldots,n$.
Since $W_0=Q_0=0$ and $f_0^\psi=\psi$,
\[
 \sum_{i=1}^n\E_x^*\Delta_i(\psi)
 =\E_x^*\psi(W_n,Q_n)-f_n^\psi(0,0) =\E_x^*\psi(W_n,Q_n)-\E_x^*\psi(S_n,Q_n).
\]
The last equality uses the independent copy in \eqref{eq:futurekernel}.
This proves \eqref{eq:duhamel}.
\end{proof}

\section{Proof of Proposition~\ref{prop:smoothed}}\label{sec:comparison}

\begin{proof}
Fix $x\in\mathcal D_3$.  Recall that the conditional increment
\[
 A\stackrel d=\xi\mid\{|\xi|\le h_x\},
\]
and $\PP_x^*$ and $\E_x^*$ denote the  conditional probability and expectation defined in \eqref{eq:003}.

\noindent\textbf{Step 1: construction of the noise.} Put
\begin{equation}\label{eq:conditionalmoments}
 \begin{aligned}
  \mu&=\E_x^* A,\qquad m_r=\E_x^*|A|^r,\qquad \tau^2=m_2-\mu^2,\\
  \gamma&=\frac{m_3}{\tau^2},\qquad
  \Lambda_y=(1+y)n|\mu|+(1+y)^3nm_3.
 \end{aligned}
\end{equation}

Let $P_A$ denote the distribution of $A$ under $\PP_x^*$.
Take independent $A,A'$ with law $P_A$.
Since $\E_x^*(A-A')^2=2\tau^2>0$, define the probability measure
\begin{equation}\label{eq:biasedpair}
 \Pi(\dd a,\dd a')=\frac{(a-a')^2}{2\tau^2}P_A(\dd a)P_A(\dd a').
\end{equation}
Draw $(B,B')$ from $\Pi$.
Independently, draw $V$ with density $(1-|v|)\ind_{[-1,1]}(v)$.
Set
\begin{equation}\label{eq:chord}
 (Z_1,Y_1)=\bigl(V(B-B'),\ V(B^2-B'^2)\bigr).
\end{equation}
Fix $J=8$.
Take two independent families
$\{(Z_j^{(1)},Y_j^{(1)})\}_{j\le J}$ and
$\{(Z_j^{(2)},Y_j^{(2)})\}_{j\le J}$,
each consisting of i.i.d.\ copies of $(Z_1,Y_1)$.
Let $V_0$ be independent of both families.
Its density $\rho$ is fixed, even and in $C_c^\infty$, with support in $[-1,1]$.
Put
\begin{equation}\label{eq:fullnoise}
 Z=\sum_{j=1}^J Z_j^{(1)}-\sum_{j=1}^JZ_j^{(2)}+\gamma V_0,
 \qquad Y=\sum_{j=1}^JY_j^{(1)}+\sum_{j=1}^JY_j^{(2)}.
\end{equation}
All auxiliary variables are independent of $A_1,\ldots,A_n$.
The law of $(Z,Y)$ depends on the conditional law fixed at $x$.
We use this same noise law for every threshold $y$.
In noise expectations, the subscript names the noise variables.

Put $D=A-A'$ and define the characteristic functions
\begin{equation}\label{eq:charfunctions}
 \begin{aligned}
 \varphi(u,v)&=\E_x^*\ee^{\ii uA+\ii vA^2},&
 K(u,v)&=\E_{Z_1,Y_1}\ee^{\ii uZ_1+\ii vY_1},\\
 \Xi(u,v)&=\E_{Z,Y}\ee^{\ii uZ+\ii vY},&
 \widehat\rho(t)&=\E_{V_0}\ee^{\ii tV_0}.
 \end{aligned}
\end{equation}

The following lemma records the properties of the noise $(Z,Y)$.

\begin{lemma}\label{lem:noisemoments}
For $x\in\mathcal D_3$, the noise $(Z,Y)$ in \eqref{eq:fullnoise}
has the following properties.

\textup{(i)} The law is even in its first coordinate,
\begin{equation}\label{eq:noisesymmetry}
 (Z,Y)\stackrel d=(-Z,Y),
\end{equation}
and
\begin{align}
 |Z|&\le34h_x,& |Y|&\le64h_x^2,\label{eq:noisesupport}\\
 \E_{Z,Y}|Z|&\le\frac{67}{3}\gamma,& \E_{Z,Y}|Y|&\le\frac{128}{3}h_x\gamma,&
 \E_{Z,Y} Z^2&\le\frac{70}{3}h_x\gamma.\label{eq:noisemoments}
\end{align}

\textup{(ii)} For real $u,v$, $0\le K(u,v)\le1$, and
\begin{equation}\label{eq:Kexact}
 K(u,v)=\frac1{\tau^2}\E_x^*
 \frac{1-\cos\{D(u+v(A+A'))\}}{(u+v(A+A'))^2}.
\end{equation}
When the denominator is zero, the quotient is defined by $D^2/2$.
The characteristic function of $(Z,Y)$ is
\begin{equation}\label{eq:Xi}
 \Xi(u,v)=K(u,v)^J K(-u,v)^J\widehat\rho(\gamma u).
\end{equation}

\textup{(iii)} On the set
\begin{equation}\label{eq:cone}
 \calC_x=\{(u,v):u\ne0,\ |v|\le |u|/(4h_x)\},
\end{equation}
one has
\begin{equation}\label{eq:gap}
 K(u,v)\asymp\frac{1-|\varphi(u,v)|^2}{\tau^2u^2}.
\end{equation}
\end{lemma}
The proof of Lemma~\ref{lem:noisemoments} is given in Supplemental Material~\Sref{supp:noisemoments-proof}.

By part~\textup{(i)} of Lemma~\ref{lem:noisemoments},
using $(1+x)h_x=\varepsilon\le1$, we obtain
\[
 \begin{aligned}
 |Z|+(1+x)|Y|&\le(34+64\varepsilon)h_x\le98h_x,\\
 \E_{Z,Y}|Z|+(1+x)\E_{Z,Y}|Y|&\le\frac{67+128\varepsilon}{3}\gamma\le65\gamma.
 \end{aligned}
\]
Thus both bounds in \eqref{eq:smoothing-size} hold with the coefficient $100$.

We also need the following moment bounds for $A$.

\begin{lemma}\label{lem:conditionalmoments}
For $x\in\mathcal D_3$, the quantities in \eqref{eq:conditionalmoments} satisfy
\begin{align}
 \tfrac34&\le n{m_2}\le\tfrac43,&\tfrac12&\le n\tau^2\le\tfrac43,\label{eq:momentbudget1}\\
 {n|\mu|}&\le\frac{2\varepsilon^{-2}R_n(x)}{1+x},&nm_3&\le\frac{2R_n(x)}{(1+x)^3},&nm_4&\le h_x nm_3,\label{eq:momentbudget2}\\
 \gamma&\asymp nm_3\le2r_n,&\gamma&\le2h_x,&{m_2}&\le\gamma^2.\label{eq:momentbudget3}
\end{align}
Here $r_n=R_n(0)$; in particular, $\gamma\le4r_n$.
For every $y\ge0$ with $|y-x|\le(1+x)^{-1}$,
one also has $\Lambda_y\le20\varepsilon^{-2}R_n(x)$.
\end{lemma}
The proof is given in Supplemental Material~\Sref{supp:conditionalmoments-proof}.

\noindent\textbf{Step 2: Taylor expansion.}\par\nobreak
Fix $y$ satisfying \eqref{eq:comparisonthresholds}.
Recall the smoothed test $\widetilde\psi_y$ in \eqref{eq:test},
with the noise $(Z,Y)$ constructed in Step~1.
As in \eqref{eq:futurekernel}, write
\[
 f_\ell(w,q)=f_\ell^{\widetilde\psi_y}(w,q)
 =\E'_x\widetilde\psi_y(w+S'_\ell,q+Q'_\ell).
\]
Proposition~\ref{prop:reflection} gives
\[
 \E_x^*\widetilde\psi_y(W_n,Q_n)-\E_x^*\widetilde\psi_y(S_n,Q_n)
 =\sum_{i=1}^n\E_x^*\Delta_i(\widetilde\psi_y).
\]
It remains to prove
\begin{equation}\label{eq:correction-target}
 \left|\sum_{i=1}^n\E_x^*\Delta_i(\widetilde\psi_y)\right|
 \le C\Lambda_y\Phibar(y).
\end{equation}

Recall that
\[
 \Delta_i(\widetilde\psi_y)
 =f_{n-i}(0,Q_i)-f_{n-i}(-U_i,Q_i).
\]
The independent term $\gamma V_0$ in $Z$ makes $f_\ell$ smooth
in its first coordinate.
Primes denote differentiation in that coordinate.
Taylor's formula gives
\begin{equation}\label{eq:taylorgrowing}
 |\Delta_i(\widetilde\psi_y)|
 \le U_i|f_{n-i}'(0,Q_i)|
 +\frac{U_i^2}{2}\sup_{-U_i\le w\le0}|f_{n-i}''(w,Q_i)|.
\end{equation}
To sum this bound, we first estimate the moments of $U_i$.

\begin{lemma}\label{lem:regulators}
Let $x\in\mathcal D_3$, and let $y$ satisfy \eqref{eq:comparisonthresholds}.
For the reflected walk $W_i$ and the quantities $U_i$ in
\eqref{eq:lindley}, we have, for $1\le i\le n$,
\begin{equation}\label{eq:regindividual}
 \E_x^*U_i\le\tau i^{-1/2}+|\mu|,\qquad
 \E_x^*U_i^2\le\frac{Cm_3}{\tau\sqrt i}.
\end{equation}
Consequently,
\begin{align}
 \sum_{i\le n}\E_x^*U_i&\le C,&
 \sum_{i\le n}\E_x^*U_i^2&\le Cnm_3,\label{eq:regtotal}\\
 \E_x^*U_i&\le C/n,&
 \E_x^*U_i^2&\le Cm_3\quad(i\ge n/2).\label{eq:reglate}
\end{align}
For $y\ge1$ and fixed sufficiently small $\kappa>0$,
\begin{equation}\label{eq:weighted}
 \sum_{i\le n/2}\E_x^*[U_i\ee^{-\kappa y^2Q_i}]
 \le\frac C{1+y},\qquad
 \sum_{i\le n/2}\E_x^*[U_i^2\ee^{-\kappa y^2Q_i}]
 \le\frac{Cnm_3}{1+y}.
\end{equation}
\end{lemma}
The proof is given in Supplemental Material~\Sref{supp:regulators-proof}.
We estimate the derivatives separately for $\gamma\le y\le2$ and $y\ge2$.

\noindent\textbf{Step 3: $\gamma\le y\le2$.}\par\nobreak
Since $\Phibar(y)\ge\Phibar(2)>0$, it suffices to prove
\[
 \left|\sum_{i=1}^n\E_x^*\Delta_i(\widetilde\psi_y)\right|
 \le C(n|\mu|+nm_3).
\]
For $i\in\{0,\ldots,n\}$, split $f_{n-i}$ into its even and odd parts:
\[
 e_{n-i}(w,q)=\tfrac12\{f_{n-i}(w,q)+f_{n-i}(-w,q)\},\qquad o_{n-i}(w,q)=\tfrac12\{f_{n-i}(w,q)-f_{n-i}(-w,q)\}.
\]
For $u\ge0$,
\[
 f_{n-i}(0,q)-f_{n-i}(-u,q)
 =e_{n-i}(0,q)-e_{n-i}(u,q)+o_{n-i}(u,q).
\]
Since $e_{n-i}'(0,q)=o_{n-i}(0,q)=0$, Taylor's formula gives, for $1\le i\le n$,
\begin{equation}\label{eq:evenoddlocal}
 |\Delta_i(\widetilde\psi_y)|
 \le U_i\sup_{0\le z\le U_i}|o_{n-i}'(z,Q_i)|
 +\frac{U_i^2}{2}\sup_{0\le z\le U_i}|e_{n-i}''(z,Q_i)|.
\end{equation}
The following lemma gives the derivative estimates used in this bound.

\begin{lemma}\label{lem:realkernel}
Let $x\in\mathcal D_3$, and let $y$ satisfy
\eqref{eq:comparisonthresholds} with $y\le2$.
For $1\le i\le n$, put $T=(n-i)m_2$.
Use the even and odd parts $e_{n-i},o_{n-i}$ above.
There is a fixed sufficiently large absolute constant $L_2$ such that
the following estimates hold.

\textup{(i)} For $i\le n/2$,
\[
 \norm{o_{n-i}'}_\infty\le C(n|\mu|+nm_3),\qquad
 \norm{e_{n-i}''}_\infty\le C.
\]

\textup{(ii)} For $i>n/2$ and $T\ge\max\{L_2\gamma^2,y^2\}$,
\begin{equation}\label{eq:realpoint}
 \norm{o_{n-i}'}_\infty
 \le Cy\{n|\mu|T^{-1/2}+nm_3T^{-3/2}\},\qquad
 \norm{e_{n-i}''}_\infty\le CyT^{-3/2}.
\end{equation}

\textup{(iii)} The sums over $i>n/2$ with $L_2\gamma^2\le T<\max\{L_2\gamma^2,y^2\}$ satisfy
\begin{align}
 m_2\sum_{\substack{n/2<i\le n\\L_2\gamma^2\le(n-i)m_2<\max\{L_2\gamma^2,y^2\}}}
 \sup_{q,\,|w|\le y/16}|o_{n-i}'(w,q)|
 &\le C(n|\mu|+nm_3),\label{eq:realsumodd}\\
 m_2\sum_{\substack{n/2<i\le n\\L_2\gamma^2\le(n-i)m_2<\max\{L_2\gamma^2,y^2\}}}
 \sup_{q,\,|w|\le y/16}|e_{n-i}''(w,q)|
 &\le C.\label{eq:realsumeven}
\end{align}
\end{lemma}
The proof is given in Supplemental Material~\Sref{supp:real}.

For $i\le n/2$, use part~\textup{(i)} of Lemma~\ref{lem:realkernel}
in \eqref{eq:evenoddlocal} and sum with \eqref{eq:regtotal}.
This gives
\[
 \sum_{i\le n/2}\E_x^*|\Delta_i(\widetilde\psi_y)|
 \le C(n|\mu|+nm_3)\sum_{i\le n/2}\E_x^*U_i
       +C\sum_{i\le n/2}\E_x^*U_i^2 \le C(n|\mu|+nm_3).
\]

Next let $i>n/2$ and $T\ge L_2\gamma^2$.
For $T\ge \max\{L_2\gamma^2,y^2\}$, apply part~\textup{(ii)} of Lemma~\ref{lem:realkernel}
and \eqref{eq:reglate}, to \eqref{eq:evenoddlocal}, 
\[
 \E_x^*|\Delta_i(\widetilde\psi_y)|
 \le\frac{Cy}{n}\{n|\mu|T^{-1/2}+nm_3T^{-3/2}\}.
\]
By Lemma~\ref{lem:conditionalmoments},
$n^{-1}\asymp m_2$ and $m_2\le\gamma^2\le y^2$.
Comparing the sums with integrals, we obtain
\[
 \begin{aligned}
 m_2\sum_{\substack{n/2<i<n\\y^2\le(n-i)m_2}}
 ((n-i)m_2)^{-1/2}
 &\le\frac{m_2}{y}+\int_0^{4/3}t^{-1/2}\dd t\le C,\\
 m_2\sum_{\substack{n/2<i<n\\y^2\le(n-i)m_2}}
 ((n-i)m_2)^{-3/2}
 &\le\frac{m_2}{y^3}+\int_{y^2}^{\infty}t^{-3/2}\dd t\le\frac Cy.
 \end{aligned}
\]
The terms $m_2/y$ and $m_2/y^3$ bound the summand with the smallest $T$.
Thus the sum of $\E_x^*|\Delta_i(\widetilde\psi_y)|$ over these indices is at most
\[
 \begin{aligned}
 &Cy m_2\sum_{\substack{n/2<i<n\\y^2\le(n-i)m_2}}
 \{n|\mu|((n-i)m_2)^{-1/2}+nm_3((n-i)m_2)^{-3/2}\}\\
 &\qquad\le C(yn|\mu|+nm_3)\le C(n|\mu|+nm_3),
 \end{aligned}
\]
since $y\le2$.

For $L_2\gamma^2\le T<\max\{L_2\gamma^2,y^2\}$, split according to $U_i\le y/16$.
On this event, \eqref{eq:evenoddlocal} only uses derivatives at $|w|\le y/16$.
Taking expectations and using \eqref{eq:reglate}, we get
\[
 \begin{aligned}
 &\sum_{\substack{i>n/2\\L_2\gamma^2\le(n-i)m_2<\max\{L_2\gamma^2,y^2\}}}
 \E_x^*\bigl[|\Delta_i(\widetilde\psi_y)|\ind_{\{U_i\le y/16\}}\bigr]\\
 &\quad\le\frac Cn\sum_{\substack{i>n/2\\L_2\gamma^2\le(n-i)m_2<\max\{L_2\gamma^2,y^2\}}}
 \left\{\sup_{q,\,|w|\le y/16}|o_{n-i}'(w,q)|
        +nm_3\sup_{q,\,|w|\le y/16}|e_{n-i}''(w,q)|\right\}\\
 &\quad\le\frac{C(n|\mu|+nm_3)}{nm_2}
           +\frac{Cm_3}{m_2}
 \le C(n|\mu|+nm_3).
 \end{aligned}
\]
The second inequality uses \eqref{eq:realsumodd}--\eqref{eq:realsumeven}.

On $U_i>y/16$, use $|\Delta_i(\widetilde\psi_y)|\le1$,
since $0\le f_{n-i}\le1$.
Since $0<(n-i)m_2<y^2$, this range has at most $y^2/m_2\le Cny^2$ indices.
Markov's inequality and \eqref{eq:reglate} give
\begin{equation}\label{eq:largepush}
 \begin{aligned}
 &\sum_{\substack{i>n/2\\L_2\gamma^2\le(n-i)m_2<\max\{L_2\gamma^2,y^2\}}}
 \E_x^*\bigl[|\Delta_i(\widetilde\psi_y)|\ind_{\{U_i>y/16\}}\bigr]\\
 &\quad\le\sum_{\substack{i>n/2\\L_2\gamma^2\le(n-i)m_2<\max\{L_2\gamma^2,y^2\}}}
 \PP_x^*(U_i>y/16)\\
 &\quad\le\frac{256}{y^2}
 \sum_{\substack{i>n/2\\L_2\gamma^2\le(n-i)m_2<\max\{L_2\gamma^2,y^2\}}}\E_x^*U_i^2
 \le\frac C{y^2}(ny^2)m_3=Cnm_3.
 \end{aligned}
\end{equation}

It remains to treat $i>n/2$ with $T<L_2\gamma^2$.
For these indices, we use the independent noise $\gamma V_0$.
Its density is $\rho_\gamma(t)=\gamma^{-1}\rho(t/\gamma)$.
A change of variable in the integral defining $f_{n-i}$
allows us to differentiate $\rho_\gamma$.
Since $0\le H_y\le1$, for $a=1,2$,
\begin{equation}\label{eq:roughderivatives}
 \norm{f_{n-i}^{(a)}}_\infty
 \le\int_\R|\rho_\gamma^{(a)}(t)|\dd t
 =\gamma^{-a}\int_\R|\rho^{(a)}(t)|\dd t
 \le C\gamma^{-a}.
\end{equation}
By Lemma~\ref{lem:conditionalmoments}, $m_2\le\gamma^2$ and $nm_2\asymp1$.
There are therefore at most $Cn\gamma^2$ such indices, including $i=n$.
Insert \eqref{eq:roughderivatives} and \eqref{eq:reglate}
into \eqref{eq:taylorgrowing} to obtain
\[
 \sum_{\substack{i>n/2\\(n-i)m_2<L_2\gamma^2}}
 \E_x^*|\Delta_i(\widetilde\psi_y)|
 \le Cn\gamma^2\left\{\frac1{n\gamma}+\frac{m_3}{\gamma^2}\right\}
 \le Cnm_3,
\]
where the last inequality uses $\gamma\asymp nm_3$.
Adding these bounds proves
\[
 \left|\sum_{i=1}^n\E_x^*\Delta_i(\widetilde\psi_y)\right|
 \le\sum_{i=1}^n\E_x^*|\Delta_i(\widetilde\psi_y)|
 \le C(n|\mu|+nm_3)
 \le C\Lambda_y\Phibar(y).
\]

\noindent\textbf{Step 4: $y\ge2$.}\par\nobreak
Set $k=\lfloor n/2\rfloor$ and $E=\{Q_k<1/16\}$.
Split the sum as
\[
 \sum_{i=1}^n\E_x^*\Delta_i(\widetilde\psi_y)
 =\sum_{i\le k}\E_x^*\Delta_i(\widetilde\psi_y)
   +\E_x^*\left[\ind_E\sum_{i>k}\Delta_i(\widetilde\psi_y)\right] +\E_x^*\left[\ind_{E^c}\sum_{i>k}\Delta_i(\widetilde\psi_y)\right].
\]
We first estimate
$\sum_{i\le k}\E_x^*\Delta_i(\widetilde\psi_y)$.

We use the following derivative estimates.

\begin{lemma}\label{lem:growingkernel}
Let $x\in\mathcal D_3$, and let $y$ satisfy
\eqref{eq:comparisonthresholds} with $y\ge2$.
Write $h=\varepsilon_0/(1+y)$ and $\Lambda=\Lambda_y$.
For $1\le i\le n$, put $T=(n-i)m_2$.
There is a fixed absolute $t_0\in(0,3/8)$ such that
the following estimates hold.

\textup{(i)} For $i\le n/2$ and $q\ge0$,
\begin{align}
 |f_{n-i}'(0,q)|&\le C\Lambda\ee^{-y^2(q+T)/(2T)},\label{eq:macroodd}\\
 \sup_{|w|\le h}|f_{n-i}''(w,q)|
 &\le C(1+y)\ee^{-y^2(q+T)/(2T)}.\label{eq:macroeven}
\end{align}

\textup{(ii)} For $i>n/2$, $T\ge t_0$ and $q\ge1/16$,
\eqref{eq:macroodd} and \eqref{eq:macroeven} also hold.

\textup{(iii)} For $i>n/2$, $L_2\gamma^2\le T<t_0$ and $q\ge1/16$,
\begin{equation}\label{eq:shortcomplex}
 |f_{n-i}'(0,q)|\le C\Lambda\Phibar(y)T^2,\qquad
 \sup_{|w|\le h}|f_{n-i}''(w,q)|\le C\Phibar(y)T^2.
\end{equation}
\end{lemma}
The proof is given in Supplemental Material~\Sref{supp:complex}.

For $i\le k$, we have $0\le U_i\le h_x\le h$.
Apply part~\textup{(i)} of Lemma~\ref{lem:growingkernel}
to \eqref{eq:taylorgrowing}, 
\[
 |\Delta_i(\widetilde\psi_y)|
 \le C\ee^{-y^2/2}\ee^{-\kappa y^2Q_i}
       \{\Lambda U_i+(1+y)U_i^2\}.
\]
Taking expectations and summing with \eqref{eq:weighted} gives
\[
 \begin{aligned}
 \sum_{i\le k}\E_x^*|\Delta_i(\widetilde\psi_y)|
 &\le C\ee^{-y^2/2}\left\{
 \Lambda\sum_{i\le k}\E_x^*[U_i\ee^{-\kappa y^2Q_i}]
 +(1+y)\sum_{i\le k}\E_x^*[U_i^2\ee^{-\kappa y^2Q_i}]
 \right\}\\
 &\le C\ee^{-y^2/2}\left\{\frac{\Lambda}{1+y}+nm_3\right\}
 \le C\Lambda\Phibar(y),
 \end{aligned}
\]
where in last inequality uses $nm_3\le\Lambda/(1+y)^3$
and Mills' inequalities $\Phibar(y)\asymp\ee^{-y^2/2}/(1+y)$.

Next we estimate
$\E_x^*\left[\ind_E\sum_{i>k}\Delta_i(\widetilde\psi_y)\right]$.

\begin{lemma}\label{lem:reflection}
Let $\calF_i=\sigma(A_1,\ldots,A_i)$.
For a measurable test $\psi:\R^2\to[0,1]$,
let $f_\ell^\psi$ and $\Delta_i(\psi)$ be as in \eqref{eq:futurekernel}.
For $0\le k\le n$ and $E\in\calF_k$,
\begin{equation}\label{eq:signedblock}
 \left|\E_x^*\left[\ind_E\sum_{i=k+1}^n\Delta_i(\psi)\right]\right|
 \le\PP_x^*(E).
\end{equation}
\end{lemma}
The proof is given in Supplemental Material~\Sref{supp:reflection-proof}.
The next lemma bounds $\PP_x^*(E)$.

\begin{lemma}\label{lem:clock}
Let $x\in\mathcal D_3$, and let $y$ satisfy
\eqref{eq:comparisonthresholds} with $y\ge2$.
For $Q_i=\sum_{j=1}^iA_j^2$,
\begin{equation}\label{eq:clock}
 \PP_x^*(Q_{\lfloor n/2\rfloor}<1/16)
 +\PP_x^*(Q_n\notin[1/4,2])
 \le C\Lambda_y\Phibar(y).
\end{equation}
\end{lemma}
The proof is given in Supplemental Material~\Sref{supp:clock-proof}.
Apply Lemma~\ref{lem:reflection} with $\psi=\widetilde\psi_y$.
The terms telescope inside the expectation.
Since $0\le\widetilde\psi_y,f_{n-k}\le1$, \eqref{eq:clock} gives
\[
 \begin{aligned}
 &\left|\E_x^*\left[\ind_E\sum_{i>k}\Delta_i(\widetilde\psi_y)\right]\right|\\
 &\quad=\left|\E_x^*\left[\ind_E
 \{\widetilde\psi_y(W_n,Q_n)-f_{n-k}(W_k,Q_k)\}\right]\right|\\
 &\quad\le\E_x^*\left[\ind_E
 |\widetilde\psi_y(W_n,Q_n)-f_{n-k}(W_k,Q_k)|\right]\\
 &\quad\le\PP_x^*(E)=\PP_x^*(Q_k<1/16)
 \le C\Lambda\Phibar(y).
 \end{aligned}
\]

Finally, we estimate
$\E_x^*\left[\ind_{E^c}\sum_{i>k}\Delta_i(\widetilde\psi_y)\right]$.
For $i>k$ and $T\ge t_0$, apply part~\textup{(ii)} of Lemma~\ref{lem:growingkernel}.
On $E^c$, $Q_i\ge1/16$, and $T\le4/3$ gives
$\ee^{-y^2Q_i/(2T)}\le\ee^{-cy^2}$.
Taking expectations in \eqref{eq:taylorgrowing} and using \eqref{eq:reglate}, we obtain
\[
 \E_x^*[\ind_{E^c}|\Delta_i(\widetilde\psi_y)|]
 \le\frac Cn\ee^{-y^2/2}\ee^{-cy^2}\{\Lambda+(1+y)nm_3\}.
\]
There are at most $n$ such indices. Hence
\[
 \sum_{\substack{k<i\le n\\(n-i)m_2\ge t_0}}
 \E_x^*[\ind_{E^c}|\Delta_i(\widetilde\psi_y)|]
 \le C\ee^{-y^2/2}\ee^{-cy^2}\{\Lambda+(1+y)nm_3\} \le C\Lambda\Phibar(y).
\]

For $i>k$ and $L_2\gamma^2\le T<t_0$, use part~\textup{(iii)}
of Lemma~\ref{lem:growingkernel} in \eqref{eq:taylorgrowing}.
Taking expectations and applying \eqref{eq:reglate} gives
\[
 \E_x^*[\ind_{E^c}|\Delta_i(\widetilde\psi_y)|]
 \le C\Phibar(y)T^2
 \{\Lambda\E_x^*U_i+\E_x^*U_i^2\} \le\frac Cn\Phibar(y)T^2(\Lambda+nm_3).
\]
Since $T<t_0<1$ and $nm_3\le\Lambda$, summing over these indices gives
\[
 \sum_{\substack{k<i\le n\\L_2\gamma^2\le(n-i)m_2<t_0}}
 \E_x^*[\ind_{E^c}|\Delta_i(\widetilde\psi_y)|] \le C\Phibar(y)\frac{\Lambda+nm_3}{n}
 \sum_{\substack{k<i\le n\\L_2\gamma^2\le(n-i)m_2<t_0}}((n-i)m_2)^2
 \le C\Lambda\Phibar(y).
\]

It remains to treat $i>k$ with $T<L_2\gamma^2$.
For these indices, we use the scalar noise $\gamma V_0$.
For $q>5$, the cutoff $\chi$ makes $f_{n-i}$ and its derivatives zero.
Take $1/16\le q\le5$ and $|w|\le h$.
On $\supp\chi$, the boundary of $H_y$ lies at least $y/4$ from zero.
The fixed choice of $\varepsilon_0$ gives $|w|+|Z|\le Ch\le y/8$.
Thus a derivative of $\widetilde\psi_y(w+S_{n-i},q+Q_{n-i})$
can be nonzero only if $|S_{n-i}|\ge y/8$.
As in \eqref{eq:roughderivatives}, differentiate the density of $\gamma V_0$ to obtain
\[
 \sup_{|w|\le h}|f_{n-i}^{(a)}(w,q)|
 \le C\gamma^{-a}\PP_x^*(|S_{n-i}|\ge y/8),\qquad a=1,2.
\]
To bound this probability, put $\lambda_0=32(1+y)$.
Since $\lambda_0h=32\varepsilon_0\le1$, Taylor's formula gives
\[
 (n-i)\log\E_x^*\ee^{\pm\lambda_0A}
 \le(n-i)\lambda_0|\mu|+C\lambda_0^2T.
\]
The first term is at most $32\Lambda$.
For the second, Lemma~\ref{lem:conditionalmoments} gives $\gamma\le2h$,
so $\lambda_0^2T\le4L_2\cdot32^2\varepsilon_0^2$.
Both terms are bounded by absolute constants.
Chernoff's inequality therefore gives
\[
 \PP_x^*(|S_{n-i}|\ge y/8)\le C\ee^{-4y^2}.
\]
For $i=n$, the derivatives vanish by the same support argument.
Insert these derivative bounds and \eqref{eq:reglate} into \eqref{eq:taylorgrowing}.
This yields
\[
 \E_x^*[\ind_{E^c}|\Delta_i(\widetilde\psi_y)|]
 \le C\ee^{-4y^2}\left\{\frac1{n\gamma}+\frac{m_3}{\gamma^2}\right\}.
\]
As in Step~3, there are at most $Cn\gamma^2$ indices in this range.
Their sum is therefore bounded by
\[
 \sum_{\substack{k<i\le n\\(n-i)m_2<L_2\gamma^2}}
 \E_x^*[\ind_{E^c}|\Delta_i(\widetilde\psi_y)|]
 \le Cn\gamma^2\ee^{-4y^2}
 \left\{\frac1{n\gamma}+\frac{m_3}{\gamma^2}\right\} \le Cnm_3\ee^{-4y^2}\le C\Lambda\Phibar(y).
\]
Combining the three ranges, we obtain
\[
 \left|\E_x^*\left[\ind_{E^c}\sum_{i>k}\Delta_i(\widetilde\psi_y)\right]\right|
 \le\sum_{i>k}\E_x^*[\ind_{E^c}|\Delta_i(\widetilde\psi_y)|]
 \le C\Lambda\Phibar(y).
\]

Adding the three terms in the decomposition proves
\eqref{eq:correction-target} for $y\ge2$.
Step~3 proves it for $\gamma\le y\le2$.
Together with \eqref{eq:duhamel}, this gives the first inequality in \eqref{eq:smoothed}.
For the second, Lemma~\ref{lem:conditionalmoments} gives $\Lambda_y\le CR_n(x)$,
and Mills' inequalities give $\Phibar(y)\asymp\Phibar(x)$.
Hence
\[
 C\Lambda_y\Phibar(y)\le CR_n(x)\Phibar(x)=Cw_n(x).
\]
\end{proof}

\section*{Acknowledgements}
Liu S.H. was partially supported by the Fundamental Research Funds for the Central Universities DUT25RC(3)133.
Shao Q.M. was partially supported by National Nature Science Foundation of China NSFC 12031005 and Shenzhen Outstanding Talents Training Fund, China.

\bmhead{AI Statement}
During the preparation of this manuscript, the authors used ChatGPT for language editing,
LaTeX formatting, checking references, notation and cross-references, and assistance
with proof organization and mathematical checks.
The authors reviewed and revised the manuscript, verified all mathematical arguments,
and take full responsibility for its content.

\section*{Supplementary material}
{{\bfseries Supplement to “Sharp Cram\'er-type moderate deviation bounds for the self-normalized maximum of partial sums”}

The Supplement contains the auxiliary estimates and proofs used in the main article.}

\clearpage
\begin{center}
{\Large\bfseries Supplement to “Sharp Cram\'er-type moderate deviation bounds for the self-normalized maximum of partial sums”}
\end{center}
\medskip
This  Supplementary material contains the auxiliary technical estimates and proofs used in the main article.
\medskip
\setcounter{section}{0}
\renewcommand{\thesection}{S.\arabic{section}}
\renewcommand{\theHsection}{supp.\arabic{section}}
\renewcommand{\theequation}{\thesection.\arabic{equation}}
\renewcommand{\thetheorem}{\thesection.\arabic{theorem}}

\section{Preliminary estimates}\label{supp:prelim}
We prove the auxiliary estimates used in Section~\Mref{sec:preliminary}.

\subsection{Proof of Lemma~\Mref{lem:truncatedbudget}}\label{supp:truncatedbudgetproof}
\begin{proof}
Put $\xi=X/(\sigma\sqrt n)$ and recall that
\[
 F(z)=z^2\min\{|z|,1\},\qquad
 R_n(t)=n\E F((1+t)\xi),\qquad
 n\E\xi^2=1,\quad n\E|\xi|^p=\beta_{n,p}.
\]
For $0\le t\le T<\infty$, we have
$F((1+t)\xi)\le(1+T)^2\xi^2$.
Dominated convergence gives continuity of $R_n$.
The integrand is strictly increasing in $t$ on $\{\xi\ne0\}$.
This event has positive probability, so $R_n$ is strictly increasing.
The bound $F(z)\le|z|^p$ gives $R_n(t)\le\beta_{n,p}(1+t)^p$.
Also,
\[
 \frac{R_n(t)}{(1+t)^2}
 =n\E\bigl[\xi^2\min\{(1+t)|\xi|,1\}\bigr]\longrightarrow1
\]
by dominated convergence.

For $a\ge1$, the definition of $F$ gives
$a^2F(z)\le F(az)\le a^3F(z)$.
Replacing $a$ by $a^{-1}$ gives the bounds for $0<a\le1$.
Thus
\begin{equation}\label{supp:Fscaling}
 \min\{a^2,a^3\}F(z)\le F(az)\le\max\{a^2,a^3\}F(z),\qquad a>0.
\end{equation}
If $x,t\ge0$ and $|t-x|\le(1+x)^{-1}$, then
$(1+t)/(1+x)\in[1/2,2]$.
Apply \eqref{supp:Fscaling} with $a=(1+t)/(1+x)$ and $z=(1+x)\xi$.
Taking expectations proves (\Mref{eq:budgetcomparison}).

To prove (\Mref{eq:budgetlower}), put
$a=1+t$ and $D=\{|\xi|\le a^{-1}\}$.
Suppose $R_n(t)\le\rho\le1/16$.
Splitting the expectation over $D$ and $D^c$ gives
\begin{equation}\label{supp:truncateddecomposition}
 R_n(t)
 =a^3n\E[|\xi|^3\ind_D]+a^2n\E[\xi^2\ind_{D^c}].
\end{equation}
Since $n\E\xi^2=1$,
\begin{equation}\label{supp:retainedvariance}
 n\E[\xi^2\ind_D]
 \ge1-\frac{R_n(t)}{a^2}\ge1-\rho.
\end{equation}
The variable $\xi\ind_D$ has a finite third absolute moment.
Lyapunov's inequality and \eqref{supp:retainedvariance} yield
\begin{equation}\label{supp:retainedthirdmoment}
 n\E[|\xi|^3\ind_D]
 \ge\frac{\{n\E[\xi^2\ind_D]\}^{3/2}}{\sqrt n}
 \ge\frac{(1-\rho)^{3/2}}{\sqrt n}.
\end{equation}
Substituting into \eqref{supp:truncateddecomposition} and rearranging gives
\[
 n\ge\frac{(1-\rho)^3(1+t)^6}{R_n(t)^2}.
\]
Here $R_n(t)>0$ because $\PP(\xi\ne0)>0$.
Finally, $(1-\rho)^3/R_n(t)^2\ge(15/16)^3/(1/16)^2>1$.
This proves (\Mref{eq:budgetlower}).
\end{proof}

\subsection{Proof of (\Mref{supp:jsw}) and (\Mref{eq:zeroatom})}
\label{supp:reference-verification}

We verify (\Mref{supp:jsw}), including the correction for $V_n=0$.
Mills' inequalities give
\begin{equation}\label{supp:mills}
 c\frac{\ee^{-z^2/2}}{1+z}\le\Phibar(z)\le C\frac{\ee^{-z^2/2}}{1+z},\qquad z\ge0.
\end{equation}
Let $t>0$ with $R_n(t)\le\rho_*$.
Put $\xi=X/(\sigma\sqrt n)$ and $\xi_j=X_j/(\sigma\sqrt n)$.
Then $\E\xi_j=0$ and $\sum_{j=1}^n\E\xi_j^2=1$.
Apply \citet[Theorems~2.1 and~2.2]{ShaoZhou2016} with both nonlinear remainders equal to zero.
For $t\ge1$, write
\[
 \delta_j(t)=\E F(t\xi_j),\qquad
 \Delta(t)=\sum_{j=1}^n\delta_j(t)\le R_n(t).
\]
Their Theorem~2.1 requires $\max_j\delta_j(t)\le1$ and $\Delta(t)\le ct^2$.
Both follow from the choice $\rho_*=\min\{1/16,c\}$ in Section~\Mref{sec:preliminary}.
Since $\Delta(t)\le\rho_*$, the exponential bounds give relative error at most
$C\Delta(t)\le CR_n(t)$.
For $0<t\le1$, their Theorem~2.2 gives absolute error at most $CR_n(t)$.
This is at most $Cw_n(t)$ because $\Phibar(t)\ge\Phibar(1)$.
One-sided limits give the same distribution-function bounds for strict events.

The event $V_n=0$ contributes an additional probability.
Put $q_0=\PP(\xi\ne0)$ and $D=\{|\xi|\le(1+t)^{-1}\}$.
H\"older's inequality and \eqref{supp:retainedvariance} give
\[
 nq_0\ge
 \frac{\{n\E[\xi^2\ind_D]\}^3}{\{n\E[|\xi|^3\ind_D]\}^2}
 \ge\frac{(1-\rho_*)^3(1+t)^6}{R_n(t)^2}.
\]
Consequently,
\[
 \PP(V_n=0)=(1-q_0)^n
 \le\exp\left\{-\frac{(1+t)^6}{2R_n(t)^2}\right\}.
\]
This estimate also holds at $t=0$.
For $a=1+t$ and $r=R_n(t)\le\rho_*\le1/16$,
\[
 \exp\{-a^6/(2r^2)\}
 \le\exp\{-1/(4r^2)\}\exp\{-a^2/(4\rho_*^2)\}
 \le Cr\frac{\ee^{-t^2/2}}{1+t}
 \le Cr\Phibar(t).
\]
Thus (\Mref{eq:zeroatom}) holds.
Adding this probability gives (\Mref{supp:jsw}) for the weak event in the main article.
The same correction, continuity of $R_n$ and one-sided limits give the bounds for strict tail events.

\subsection{Proof of Lemma~\Mref{supp:origin}}\label{supp:originproof}
For the estimate near zero, we use the following bound on the maximum of a random walk.
Let $B,B_1,B_2,\ldots$ be i.i.d.\ with mean $\mu_B$, variance $\tau_B^2>0$
and $m_{3,c}=\E|B-\mu_B|^3<\infty$.
Fix an absolute constant $C\ge1$ for the Berry--Esseen inequality \citep{Petrov1975}:
\begin{equation}\label{supp:BE}
 \sup_z\left|\PP\left(\frac{\sum_{j=1}^m(B_j-\mu_B)}{\tau_B\sqrt m}\le z\right)-\Phi(z)\right|
 \le C\frac{m_{3,c}}{\tau_B^3\sqrt m}.
\end{equation}
The same bound holds for strict distribution functions.
Subtracting two bounds controls interval probabilities.
Integration by parts gives the bound for monotone functions of total variation at most one.

Put $S_k^B=\sum_{j=1}^kB_j$ and $M_m^B=\max(0,S_1^B,\ldots,S_m^B)$.
Spitzer's identity \citep{Spitzer1956} is
\begin{equation}\label{supp:spitzer}
 \sum_{m\ge0}z^m\E\ee^{-sM_m^B}
 =\exp\left\{\sum_{k\ge1}\frac{z^k}{k}\E\ee^{-s(S_k^B)^+}\right\},
 \qquad s>0,\quad0<z<1.
\end{equation}
See \citet[Proposition~2]{JanssenVanLeeuwaarden2018}.
For a general increment law, put
$B_N=N^{-1}\lfloor N\max(-N,\min(B,N))\rfloor$.
Apply that proposition to $NB_N$ with Laplace parameter $s/N$.
As $N\to\infty$, each finite path and its maximum converge.
The Laplace functions are bounded by one, so dominated convergence gives \eqref{supp:spitzer}.

\begin{lemma}\label{supp:survival}
Under the assumptions of \eqref{supp:BE}, for $a\ge0$ and $m\ge0$,
\begin{equation}\label{supp:survivalbound}
 \PP(M_m^B\le a)
 \le C\min\left\{1,\frac{a+m_{3,c}/\tau_B^2}{\tau_B\sqrt{m+1}}\right\}
 \exp\left\{C\frac{|\mu_B|\sqrt m}{\tau_B}\right\}.
\end{equation}
\end{lemma}
\begin{proof}
Fix $s>0$ and put $R=m_{3,c}/\tau_B^3+(s\tau_B)^{-1}\ge1$.
For $Z\sim N(0,k\tau_B^2)$, Mills' upper bound gives
\[
 \E\ee^{-sZ^+}
 =\tfrac12+\ee^{s^2k\tau_B^2/2}\Phibar(s\tau_B\sqrt k)
 \le\tfrac12+C\min\{1,(s\tau_B\sqrt k)^{-1}\}.
\]
The function $z\mapsto\ee^{-sz^+}$ has total variation one.
Apply \eqref{supp:BE} and bound the effect of the mean $k\mu_B$ on the normal distribution function.
Together with the bound by one, this gives
\[
 \E\ee^{-s(S_k^B)^+}
 \le\tfrac12+\min\{\tfrac12,CRk^{-1/2}\}
       +C|\mu_B|\sqrt k/\tau_B.
\]
Split the sum at $k=R^2$:
\[
 \sum_{k\ge1}\frac{z^k}{k}\min\{\tfrac12,CRk^{-1/2}\}
 \le\tfrac12\sum_{k\le R^2}\frac1k+CR\sum_{k>R^2}k^{-3/2}
 \le\log R+C.
\]
For $m\ge1$, take $z=\ee^{-1/m}$.
Integral comparison gives $\sum_{k\ge1}z^kk^{-1/2}\le C\sqrt m$
and $(1-z)^{-1/2}\le C\sqrt m$.
Substituting these bounds into \eqref{supp:spitzer} yields
\[
 \sum_{j\ge0}z^j\E\ee^{-sM_j^B}
 \le CR\sqrt m\exp\{C|\mu_B|\sqrt m/\tau_B\}.
\]
Since $M_j^B$ is nondecreasing in $j$, the terms $0\le j\le m$ on the left sum to at least
$\ee^{-1}(m+1)\E\ee^{-sM_m^B}$.
Hence
\[
 \E\ee^{-sM_m^B}
 \le CRm^{-1/2}\exp\{C|\mu_B|\sqrt m/\tau_B\}.
\]
For $a>0$, use $\PP(M_m^B\le a)\le\ee^{sa}\E\ee^{-sM_m^B}$
with $s=a^{-1}$.
The bound by one and $\sqrt{m+1}\le\sqrt{2m}$ give \eqref{supp:survivalbound}.
For $a=0$, let $s\to\infty$ in the Laplace bound.
For $m=0$, the conclusion follows from $m_{3,c}/\tau_B^3\ge1$.
\end{proof}

\begin{proof}
Put $\zeta=X/(\sigma\sqrt m)$ and $\zeta_j=X_j/(\sigma\sqrt m)$.
Then $m\E\zeta^2=1$ and
\[
 r_m=m\E[\zeta^2\ind_{\{|\zeta|>1\}}]
       +m\E[|\zeta|^3\ind_{\{|\zeta|\le1\}}].
\]
If $r_m$ is bounded below by a fixed positive constant, both conclusions follow from the bound by one.
Assume $r_m$ is sufficiently small.
Put $q=\PP(|\zeta|>1)$ and $B\stackrel d=\zeta\mid\{|\zeta|\le1\}$.
The definition of $r_m$ and $\E\zeta=0$ give
\[
 mq\le r_m,\qquad m|\E B|\le2r_m,\qquad
 \tfrac34\le m\E B^2\le\tfrac43,\qquad m\E|B|^3\le2r_m.
\]
Since $m(\E B)^2\le4r_m^2/m$, we have $m\var(B)\asymp1$.
Also, $\E|B-\E B|^3\le8\E|B|^3$.
For independent copies $B_1,\ldots,B_m$, put
\[
 M_B^+=\max\left(0,B_1,B_1+B_2,\ldots,\sum_{j=1}^mB_j\right),
 \qquad Q_B=\sum_{j=1}^mB_j^2.
\]
In \eqref{supp:survivalbound}, the denominator $\sqrt{(m+1)\var(B)}$ is comparable to one.
The other terms in that bound satisfy
\[
 \frac{\E|B-\E B|^3}{\var(B)}=O(r_m),\qquad
 \frac{|\E B|\sqrt m}{\sqrt{\var(B)}}=O(r_m).
\]
Thus
\[
 \PP(M_B^+\le v)\le C\min\{1,v+r_m\},\qquad v\ge0.
\]
Conditional on $|\zeta_j|\le1$ for every $j$, the normalized sample has this product law.
The complementary event has probability at most $mq\le r_m$.
Take $v=a/(\sigma\sqrt m)$ and add this probability.
Together with the bound by one, this proves (\Mref{supp:lowsurvivalbound}).

For the square sum, $|B|\le1$ gives
\[
 \var(Q_B)\le m\E B^4\le m\E|B|^3\le2r_m,
 \qquad \tfrac34\le\E Q_B\le\tfrac43.
\]
The mean is at least $1/6$ from the complement of $[1/2,3/2]$.
Chebyshev's inequality therefore gives
$\PP(|Q_B-1|>1/2)\le Cr_m$.
The conditional law of $V_m^2/(m\sigma^2)$ is that of $Q_B$.
Adding the probability $mq\le r_m$ and setting $m=n$ proves (\Mref{supp:squaresbound}).
\end{proof}

\subsection{Proof of Lemma~\Mref{supp:boundedrandom}}\label{supp:boundedrandomproof}
For the bound with a random denominator, we first compare the maximum with the final sum.
\begin{lemma}\label{supp:completion}
Let $B,B_1,\ldots,B_k$ be i.i.d., where $k\ge1$, with mean $\mu_B$,
variance $\tau_B^2>0$ and $m_{3,c}=\E|B-\mu_B|^3<\infty$.
Put
\[
 a=2k|\mu_B|+16Cm_{3,c}/\tau_B^2,
\]
where $C$ is the constant fixed in \eqref{supp:BE}.
For every $b\ge0$,
\begin{equation}\label{supp:completionbound}
 \PP\left(\max_{j\le k}\sum_{i\le j}B_i\ge b\right)
 \le4\PP\left(\sum_{i\le k}B_i\ge b-a\right).
\end{equation}
\end{lemma}
\begin{proof}
For $0\le j\le k$, we claim that $\PP(\sum_{i\le j}B_i\ge-a)\ge1/4$.
The claim is immediate for $j=0$.
For $j\ge1$, we have $j|\mu_B|\le a/2$.
If $\tau_B\sqrt j\le a/4$, Chebyshev's inequality bounds the complementary probability
by $4j\tau_B^2/a^2\le1/4$.
Otherwise, the error in \eqref{supp:BE} is at most
$4Cm_{3,c}/(a\tau_B^2)\le1/4$.
The normal upper tail is at least $1/2$, which proves the claim.

Let $T=\inf\{1\le j\le k:\sum_{i=1}^jB_i\ge b\}$,
with $T=\infty$ if the set is empty.
On $\{T=j\}$, a remaining sum of at least $-a$ gives a final sum of at least $b-a$.
Independence of the remaining increments gives
\[
 \PP\left(\sum_{i=1}^kB_i\ge b-a\right)
 \ge\sum_{j=1}^k\PP(T=j)\PP\left(\sum_{i=j+1}^kB_i\ge-a\right)
 \ge\tfrac14\PP(T\le k).
\]
This proves \eqref{supp:completionbound}.
\end{proof}

We estimate the final-sum probability by exponential tilting.
\begin{lemma}\label{supp:boundedendpoint}
Fix $D\ge1$. Let $x\ge1$, $k\ge1$, and let $B,B_1,\ldots,B_k$ be i.i.d.\ satisfying
\begin{equation}\label{supp:boundedassumptions}
 |B|\le x^{-1},\quad \nu=k\E B^2\in[1/4,4],\quad
 k|\E B|\le D/x,\quad k\E|B|^3\le D/x^3.
\end{equation}
For $0\le b\le3x$,
\begin{equation}\label{supp:boundedendpointbound}
 \PP\left(\max_{j\le k}\sum_{i\le j}B_i\ge b\right)
 \le C_D\Phibar(b/\sqrt\nu).
\end{equation}
\end{lemma}
\begin{proof}
For bounded $x$, the conclusion follows from $b/\sqrt\nu\le6x$ and the bound by one.
Assume $x$ is sufficiently large in terms of $D$.
Put $\mu_B=\E B$, $Y=B-\mu_B$ and $\tau_B^2=\var(B)$.
Then
\[
 |Y|\le2/x,\qquad k\E|Y|^3\le8D/x^3,\qquad
 k\mu_B^2\le D/x^2,\qquad k\tau_B^2\in[1/8,4].
\]
The value of $a$ in Lemma~\ref{supp:completion} satisfies $a\le C_D/x$.
That lemma reduces the proof to bounding
$\PP(\sum_iY_i\ge t)$, where $t=b-a-k\mu_B$.
Put $z=t/(\sqrt k\,\tau_B)$.
If $z\le1$, then $b/\sqrt\nu$ is bounded in terms of $D$, so the conclusion follows from the bound by one.
Assume $z>1$ and set $\lambda=t/(k\tau_B^2)$.
We have $c\le\lambda\le Cx$ and $|\lambda Y|\le C$.
Taylor's inequality gives
\begin{equation}\label{supp:boundedmgf}
 k\log\E\ee^{\lambda Y}
 \le\lambda^2k\tau_B^2/2+C\lambda^3k\E|Y|^3
 \le\lambda^2k\tau_B^2/2+C_D.
\end{equation}

Define $\dd\PP_\lambda/\dd\PP=\ee^{\lambda Y}/\E\ee^{\lambda Y}$
and use its product law for the sample.
The density is bounded above and below by positive absolute constants.
The independent-copy formula then gives $\var_\lambda(Y)\asymp\var(Y)$.
Also, $\E_\lambda|Y-\E_\lambda Y|^3\le C\E|Y|^3$.
Hence the tilted sum has variance comparable to one and Berry--Esseen error at most $C_Dx^{-3}$.
For every interval $I$ of length $a_0\ge0$, \eqref{supp:BE} yields
\begin{equation}\label{supp:tiltedinterval}
 \PP_\lambda\left(\sum_{i=1}^kY_i\in I\right)
 \le C_D(a_0+x^{-3}).
\end{equation}

Changing measure gives
\[
 \PP\left(\sum_iY_i\ge t\right)
 =\ee^{k\log\E\ee^{\lambda Y}-\lambda t}
  \E_\lambda\left[\ee^{-\lambda(\sum_iY_i-t)};\sum_iY_i\ge t\right].
\]
Partition $[t,\infty)$ into intervals
$[t+j/\lambda,t+(j+1)/\lambda)$, $j\ge0$.
By \eqref{supp:tiltedinterval},
\[
 \E_\lambda\left[\ee^{-\lambda(\sum_iY_i-t)};\sum_iY_i\ge t\right]
 \le\sum_{j\ge0}\ee^{-j}C_D(\lambda^{-1}+x^{-3})
 \le C_D/(1+z).
\]
Since $\lambda t-\lambda^2k\tau_B^2/2=z^2/2$, \eqref{supp:boundedmgf} gives
$\PP(\sum_iY_i\ge t)\le C_D\ee^{-z^2/2}/(1+z)$.
Finally,
\[
 |t-b|\le C_D/x,\qquad |k\tau_B^2-\nu|\le C_D/x^2,\qquad
 |z^2-b^2/\nu|\le C_D.
\]
Thus $1+z\asymp_D1+b/\sqrt\nu$.
Mills' inequalities and Lemma~\ref{supp:completion} give \eqref{supp:boundedendpointbound}.
\end{proof}

\begin{proof}
Use $\nu=k\E\widetilde\xi^2\in[1/2,2]$,
$Q=\sum_{i=1}^k\widetilde\xi_i^2$ and $M=\max_{j\le k}\sum_{i\le j}\widetilde\xi_i$
from the statement.
For bounded $x$, the conclusion follows from the bound by one.
Assume $x$ is sufficiently large in terms of $D$.
The increment law satisfies \eqref{supp:boundedassumptions}, and
$k\E\widetilde\xi^4\le x^{-1}k\E|\widetilde\xi|^3\le D/x^4$.
Split the crossing event according to $Q\ge\nu$ and $Q<\nu$.
On $Q\ge\nu$, Lemma~\ref{supp:boundedendpoint}, with $b=x\sqrt\nu\le3x$, gives
\[
 \PP(M\ge x\sqrt Q,\ Q\ge\nu)
 \le\PP(M\ge x\sqrt\nu)\le C_D\Phibar(x).
\]

For $Q<\nu$, put $s=-4x^2/\nu$ and $\mathcal Z_s=\E\ee^{s\widetilde\xi^2}$.
Define $\dd\PP_s/\dd\PP=\ee^{s\widetilde\xi^2}/\mathcal Z_s$
and use its product law for the sample.
Write $\E_s$ for expectation under this law.
Since $|s|\widetilde\xi^2\le8$, we have $\mathcal Z_s\ge\ee^{-8}$.
The identity
\begin{equation}\label{supp:squarecov}
 \E_s f(\widetilde\xi)-\E f(\widetilde\xi)
 =\frac{\E[(f(\widetilde\xi)-\E f(\widetilde\xi))(\ee^{s\widetilde\xi^2}-1)]}{\mathcal Z_s}
\end{equation}
and $|\ee^{s\widetilde\xi^2}-1|\le|s|\widetilde\xi^2$ give
\begin{equation}\label{supp:squaretiltmoments}
 k|\E_s\widetilde\xi-\E\widetilde\xi|\le C_D/x,\qquad
 |\nu_s-\nu|\le C_D/x^2,\qquad k\E_s|\widetilde\xi|^3\le C_D/x^3,
\end{equation}
where $\nu_s=k\E_s\widetilde\xi^2$.
For the mixed terms, use
$|\E\widetilde\xi|\E\widetilde\xi^2\le\E|\widetilde\xi|^3$
and $(\E\widetilde\xi^2)^2\le\E\widetilde\xi^4$.
The covariance of $\widetilde\xi^2$ and $\ee^{s\widetilde\xi^2}$ is nonpositive, so $\nu_s\le\nu$.
The bound \eqref{supp:squaretiltmoments} also gives $\nu_s\ge1/4$.
Thus \eqref{supp:boundedassumptions} holds under $\PP_s$ with a constant depending only on $D$.
Using $\ee^u\le1+u+u^2/2$ for $u\le0$, we obtain
\begin{equation}\label{supp:squaremgf}
 k\log\mathcal Z_s\le s\nu+\tfrac12s^2k\E\widetilde\xi^4
 \le-4x^2+32D.
\end{equation}

Put $J_0=\lfloor x^2\rfloor$.
Partition the range $Q<\nu$ at $u_j=\nu(1-j/x^2)$, $0\le j\le J_0$.
For $0\le j<J_0$, set $y_j=\sqrt{x^2-j-1}$ and
$\mathcal A_j=\{M\ge x\sqrt Q,\ u_{j+1}<Q\le u_j\}$.
On $\mathcal A_j$, $M\ge\sqrt\nu\,y_j$.
Changing measure and using \eqref{supp:squaremgf} give
\[
 \PP(\mathcal A_j)
 =\mathcal Z_s^k\E_s[\ee^{-sQ};\mathcal A_j] \le\exp\{k\log\mathcal Z_s-su_j\}
       \PP_s(M\ge\sqrt\nu\,y_j) \le C_D\ee^{-4j}\Phibar(y_j).
\]
For the last inequality, apply Lemma~\ref{supp:boundedendpoint} under $\PP_s$.
Its threshold is at most $\sqrt\nu\,x\le3x$,
and $\nu_s\le\nu$ gives $\sqrt{\nu/\nu_s}\,y_j\ge y_j$.
To sum these bounds, use Mills' inequalities:
\[
 \frac{\Phibar(y_j)}{\Phibar(x)}
 \le C\frac{1+x}{1+y_j}\ee^{(j+1)/2}.
\]
For $j+1\le x^2/2$, we have $y_j\ge x/\sqrt2$.
The sum of the normalized bounds is at most $C_D\sum_{j\ge0}\ee^{-7j/2}$.
For the remaining indices, it is at most
$C_D(1+x)\sum_{j\ge x^2/2-1}\ee^{-7j/2}\le C_D$.
The lowest layer satisfies $0\le Q\le u_{J_0}\le\nu/x^2$.
The same change of measure gives
\[
 \PP(Q\le u_{J_0})
 \le\mathcal Z_s^k\ee^{-su_{J_0}}
 \le\ee^{32D-4J_0}\le C_D\Phibar(x).
\]
This includes $Q=0$.
Adding the bounds for $Q\ge\nu$, the intermediate layers and the lowest layer proves
(\Mref{supp:boundedrandombound}).
\end{proof}

\subsection{Proof of Lemma~\Mref{supp:deletion}}\label{supp:deletionproof}
\begin{proof}
Write $v^2$ and $z^2$ for the retained and deleted square sums.
On $\{M_n\ge xV_n\}$, choose a prefix that crosses $x\sqrt{v^2+z^2}$.
Cauchy--Schwarz bounds its deleted part above by $\sqrt j\,z$.
The retained part therefore has sum at least
\[
 x\sqrt{v^2+z^2}-\sqrt j\,z\ge\sqrt{x^2-j}\,v,
\]
because $\sqrt{x^2-j}\,v+\sqrt j\,z\le x\sqrt{v^2+z^2}$.
If $v>0$, this sum is positive, so the retained prefix is nonempty.
If $v=0$, all retained observations vanish.
At least one remains because $j<n$, so the retained weak crossing event holds.
This proves (\Mref{eq:deletion}).
\end{proof}

\section{Proofs of the sharpness asymptotics and auxiliary estimates}\label{supp:sharpness-lemmas}
We first give two estimates used in the proof of Lemma~\Mref{lem:conditional-tails}.
\begin{lemma}\label{lem:scale-separation}
Fix $2<p<3$.
Let $b_k,\pi_k,N_k,x_k$ be as in
(\Mref{eq:sharp-parameters}) and (\Mref{eq:sharp-subsequence}).
Put $\delta_k=N_k\pi_k$, $H_1=1$ and $H_k=b_{k-1}$ for $k\ge2$.
Then $\delta_k\to0$, $H_k^2/N_k\to0$, and
\begin{equation}\label{eq:scale-errors}
 \frac{(1+x_k)^3H_k}{\sqrt{N_k}}=o(\delta_k).
\end{equation}
Moreover,
\begin{equation}\label{eq:scale-tails}
 N_k\sum_{j>k}\pi_j=o\{\delta_k\Phibar(x_k)\},\qquad
 \sum_{j\ge k}\pi_jb_j^2\sim\pi_kb_k^2.
\end{equation}
\end{lemma}
\begin{proof}
Put $n=N_k$, $x=x_k$, $H=H_k$ and $\delta=\delta_k$.
For $k\ge2$, we have $N_{k-1}=n^{1/\nu}$ and $x_{k-1}=x/\sqrt\nu$.
Hence
\begin{equation}\label{eq:sharp-separated-parameters}
 H=\frac{2\sqrt\nu\,n^{1/(2\nu)}}x,\qquad
 \delta=2^{-p-k-2}n^{-(p-2)/2}x^p.
\end{equation}
Since $x^2=\log n$ and $p>2$, the second formula gives $\delta\to0$.
The first gives $H^2/n=4\nu n^{1/\nu-1}/x^2\to0$.
Also,
\[
 \frac{(1+x)^3H}{\delta\sqrt n}
 \le C_p2^{k+2}n^{-(3-p)/2+1/(2\nu)}x^{2-p}.
\]
The exponent of $n$ is negative because $\nu\ge8/(3-p)$.
Here $k=O(\log\log n)$, so $2^k=n^{o(1)}$.
The negative power of $n$ therefore makes the right-hand side tend to zero.
This proves \eqref{eq:scale-errors}.

To prove \eqref{eq:scale-tails}, first note that
\begin{equation}\label{eq:rare-probability-ratio}
 \frac{\pi_{k+1}}{\pi_k}
 =\frac12\nu^{p/2}N_k^{-p(\nu-1)/2}\longrightarrow0.
\end{equation}
These ratios decrease. In particular,
\[
 \pi_{k+1}\le\sum_{j>k}\pi_j
 \le\frac{\pi_{k+1}}{1-\pi_{k+2}/\pi_{k+1}}\sim\pi_{k+1}.
\]
Using $\Phibar(x_k)\asymp N_k^{-1/2}/x_k$, we obtain
\[
 \frac{N_k\sum_{j>k}\pi_j}{\delta_k\Phibar(x_k)}
 \asymp\frac{\pi_{k+1}}{\pi_k\Phibar(x_k)}
 \le C_p x_kN_k^{-\{p(\nu-1)-1\}/2}\longrightarrow0.
\]
For the second sum,
\[
 \frac{\pi_{k+1}b_{k+1}^2}{\pi_kb_k^2}
 =\frac12\nu^{(p-2)/2}N_k^{-(p-2)(\nu-1)/2}\longrightarrow0.
\]
These ratios also decrease.
Thus $\sum_{j\ge k}\pi_jb_j^2\sim\pi_kb_k^2$, which completes the proof.
\end{proof}

The next lemma gives tail estimates for weighted sums of independent signs and their partial-sum maximum.
\begin{lemma}\label{lem:weighted-signs}
Let $m\ge1$ and let $\epsilon_1,\ldots,\epsilon_m$ be independent Rademacher variables.
Let $0<a_i\le H$ be deterministic.
For $1\le j\le m$, put
\[
 T_j=\sum_{i=1}^j a_i\epsilon_i,\qquad
 Q_j=\sum_{i=1}^j a_i^2,\qquad T_0=Q_0=0,\qquad M_j^T=\max_{0\le\ell\le j}T_\ell,\qquad Q=Q_m.
\]
If $z\ge1$ and $z^4H^2/Q$ is sufficiently small, then
\begin{equation}\label{eq:weighted-sign-tails}
 \begin{split}
 \frac{\PP(T_m\ge z\sqrt Q)}{\Phibar(z)}
 &=1+O\left(\frac{(1+z)H}{\sqrt Q}+\frac{z^4H^2}{Q}\right),\\
 \frac{\PP(M_m^T\ge z\sqrt Q)}{2\Phibar(z)}
 &=1+O\left(\frac{(1+z)H}{\sqrt Q}+\frac{z^4H^2}{Q}\right).
 \end{split}
\end{equation}
The first estimate also holds with a strict inequality.
For every $D>0$ and $1\le j\le m$,
\begin{equation}\label{eq:weighted-sign-coarse}
 \PP(M_j^T\ge D)\le2\exp\{-D^2/(2Q_j)\}.
\end{equation}
All constants are absolute.
\end{lemma}
\begin{proof}
We first estimate $\PP(T_m\ge h)$, where $h=z\sqrt Q$.
Put
\[
 K(t)=\sum_{i=1}^m\log\cosh(ta_i).
\]
The assumption on $z^4H^2/Q$ implies that $zH/\sqrt Q$ is small.
Since $\tanh u\le u$ for $u\ge0$, $K'(z/\sqrt Q)\le h$.
Also, $\tanh u=u+O(u^3)$ and $\sum_i a_i^4\le H^2Q$ give
\[
 K'(2z/\sqrt Q)\ge2h-Cz^3H^2/\sqrt Q\ge h.
\]
Thus $K''>0$ gives a unique $t\in[z/\sqrt Q,2z/\sqrt Q]$ with $K'(t)=h$.
On this interval,
\[
 K'(t)=tQ+O(t^3H^2Q),\qquad
 K(t)=t^2Q/2+O(t^4H^2Q),\qquad K''(t)=Q+O(t^2H^2Q).
\]
Substitute $K'(t)=h$ to obtain, uniformly in the weights,
\begin{align}
 t&=\frac z{\sqrt Q}+O\left(\frac{z^3H^2}{Q^{3/2}}\right),&
 Q_t:=K''(t)&=Q\{1+O(z^2H^2/Q)\},\label{eq:weighted-saddle}\\
 I:=th-K(t)&=\frac{z^2}{2}+O(z^4H^2/Q),&
 k_t:=t\sqrt{Q_t}&=z+O(z^3H^2/Q).\label{eq:weighted-rate}
\end{align}
Define $\dd\PP_t=\exp\{tT_m-K(t)\}\dd\PP$ and write $\E_t$ for its expectation.
Under this law, the summands $Y_i=a_i\epsilon_i$ remain independent.
Their total mean is $h$ and their total variance is $Q_t$.
For $U=(T_m-h)/\sqrt{Q_t}$, the Berry--Esseen inequality \citep{Petrov1975} gives
\[
 \sup_v\left|\PP_t(U\le v)-\Phi(v)\right|
 \le C\frac{\sum_i\E_t|Y_i-\E_tY_i|^3}{Q_t^{3/2}}
 \le\frac{CH}{\sqrt Q}.
\]
Here $\sum_i\E_t|Y_i-\E_tY_i|^3\le C\sum_i a_i^3\le CHQ$ and $Q_t\asymp Q$.
The function $v\mapsto\ee^{-k_tv}\ind_{\{v\ge0\}}$ has total variation two.
Changing measure therefore gives
\begin{equation}\label{eq:rade-tilt}
 \PP(T_m\ge h)
 =\ee^{-I}\E_t[\ee^{-k_tU}\ind_{\{U\ge0\}}] =\ee^{-I}\{J(k_t)+O(H/\sqrt Q)\},
\end{equation}
where, for $k\ge0$,
\[
 J(k)=\int_0^\infty\ee^{-kv}\phi(v)\dd v
 =\ee^{k^2/2}\Phibar(k)\asymp(1+k)^{-1}.
\]
Differentiating under the integral gives
\[
 |J'(k)|=\int_0^\infty v\ee^{-kv}\phi(v)\dd v\le C(1+k)^{-2},\qquad
 |(\log J)'(k)|\le\frac C{1+k}.
\]
Replacing $\{U\ge0\}$ by $\{U>0\}$ gives the same error.
Since $J(k_t)\asymp(1+z)^{-1}$, the relative error in \eqref{eq:rade-tilt}
is $O((1+z)H/\sqrt Q)$.
Equation \eqref{eq:weighted-rate} and the bound on $(\log J)'$ give
\[
 |\log J(k_t)-\log J(z)|\le C\frac{|k_t-z|}{1+z}
 \le C\frac{z^2H^2}{Q}.
\]
Therefore
\[
 \frac{\PP(T_m\ge z\sqrt Q)}{\Phibar(z)}
 =\ee^{z^2/2-I}\frac{J(k_t)}{J(z)}
   \left\{1+O\left(\frac{(1+z)H}{\sqrt Q}\right)\right\}
 =1+O\left(\frac{(1+z)H}{\sqrt Q}+\frac{z^4H^2}{Q}\right).
\]
This proves the first estimate in \eqref{eq:weighted-sign-tails}, including its strict version.

The bound for $M_m^T$ follows by reflection.
For any $h>0$, let $\tau$ be the first crossing of $h$.
On $\{\tau\le m\}$, the remaining signed sum is conditionally symmetric.
It is nonnegative with probability at least $1/2$.
Since $T_\tau\ge h$, we obtain
\begin{equation}\label{eq:weighted-reflection-upper}
 \PP(M_m^T\ge h)\le2\PP(T_m\ge h).
\end{equation}
For the reverse bound, take a path with $T_m>h+2H$.
Reflect its signs after the first crossing of $h$.
Since $T_\tau<h+H$, the reflected final value is $2T_\tau-T_m<h$.
The crossing time and the preceding path are unchanged.
A second reflection recovers the original path, so the map is injective.
Each sign vector has probability $2^{-m}$; hence the map preserves probability.
Its image is disjoint from $\{T_m\ge h\}$, giving
\begin{equation}\label{eq:weighted-reflection-lower}
 \PP(M_m^T\ge h)\ge\PP(T_m\ge h)+\PP(T_m>h+2H).
\end{equation}
Apply the first estimate in \eqref{eq:weighted-sign-tails} at $z$ and $z+2H/\sqrt Q$.
The shifted threshold satisfies
$(z+2H/\sqrt Q)^4H^2/Q\le16z^4H^2/Q$.
Also,
\[
 \frac{\Phibar(z+2H/\sqrt Q)}{\Phibar(z)}
 =1+O\left(\frac{(1+z)H}{\sqrt Q}\right).
\]
Substitution into \eqref{eq:weighted-reflection-upper}--\eqref{eq:weighted-reflection-lower}
proves the second estimate in \eqref{eq:weighted-sign-tails}.

Finally, $\log\cosh u\le u^2/2$ gives $\E\ee^{tT_j}\le\ee^{t^2Q_j/2}$.
Chernoff's inequality yields
\[
 \PP(T_j\ge D)\le\exp\{-D^2/(2Q_j)\},\qquad D>0.
\]
Applying \eqref{eq:weighted-reflection-upper} to the first $j$ terms proves \eqref{eq:weighted-sign-coarse}.
For $j=0$, the crossing probability is zero.
\end{proof}

\subsection{Proof of Lemma~\Mref{lem:conditional-tails}}\label{supp:conditional-tails-proof}
\begin{proof}
Let $H_k$ be as in Lemma~\ref{lem:scale-separation}, and write
\[
 n=N_k,\qquad x=x_k,\qquad b=b_k,\qquad H=H_k,\qquad \delta=\delta_k.
\]
Recall $\sigma^2=\E_{P_*}X^2<5/4$ and $a=2/\sigma>1$.
We have $b=2\sqrt n/x$ and $x^2=\log n$.
We first prove (\Mref{eq:sharp-rn-scale}).
Put $t=\sigma\sqrt n/(1+x)$.
By \eqref{eq:scale-errors}, $H/t\to0$, and
\[
 \frac bt=\frac2\sigma\frac{1+x}{x}\longrightarrow a>1.
\]
Thus $H<t<b$ for all sufficiently large $k$.
The two terms in $R_n(x;P_*)$ satisfy
\[
 \begin{aligned}
 \frac{(1+x)^2}{\sigma^2}\E_{P_*}[X^2;|X|>t]
 &=\frac{(1+x)^2}{\sigma^2}\sum_{j\ge k}\pi_jb_j^2
 \sim a^2\delta,\\
 \frac{(1+x)^3}{\sigma^3\sqrt n}\E_{P_*}[|X|^3;|X|\le t]
 &\le\frac{(1+x)^3H}{\sigma\sqrt n}
 =o(\delta).
 \end{aligned}
\]
Here we used \eqref{eq:scale-tails} for the first term.
For the second, $|X|^3\le HX^2$ on $|X|\le t$, and \eqref{eq:scale-errors} applies.
Adding the two terms proves (\Mref{eq:sharp-rn-scale}), since $\delta\to0$ by Lemma~\ref{lem:scale-separation}.

For (\Mref{eq:sharp-tail-expansion}), condition the law on magnitudes at most $b$:
\[
 P^{[k]}=P_*(\,\cdot\mid |X|\le b).
\]
Given $\max_{i\le n}|X_i|\le b$, the sample has the product law of $P^{[k]}$.
The two event probabilities can differ by at most
$\PP_{P_*}(\max_{i\le n}|X_i|>b)\le n\sum_{j>k}\pi_j$.
Equation \eqref{eq:scale-tails} therefore gives
\begin{equation}\label{eq:sharp-law-truncation}
 \left|\PP_{P_*}(M_n\ge xV_n)-\PP_{P^{[k]}}(M_n\ge xV_n)\right|
 \le n\sum_{j>k}\pi_j=o\{\delta\Phibar(x)\}.
\end{equation}
Under this law, define
\begin{equation}\label{eq:rare-count}
 K=\#\{i:|X_i|=b\},\qquad
 q_j=\PP_{P^{[k]}}(M_n\ge xV_n\mid K=j).
\end{equation}
Then $K\sim\operatorname{Bin}(n,\widetilde\pi_k)$, where
\[
 \widetilde\pi_k=\frac{\pi_k}{1-\sum_{\ell>k}\pi_\ell},\qquad
 \widetilde\delta=n\widetilde\pi_k\sim\delta.
\]
Split the crossing probability as
\begin{equation}\label{eq:rare-decomposition}
 \begin{aligned}
 \PP_{P^{[k]}}(M_n\ge xV_n)
 &=\PP_{P^{[k]}}(K=0)q_0+\PP_{P^{[k]}}(K=1)q_1\\
 &\quad+\PP_{P^{[k]}}(M_n\ge xV_n,\ K\ge2).
 \end{aligned}
\end{equation}
Since $\widetilde\delta\sim\delta\to0$, expanding
$(1-\widetilde\pi_k)^n$ and $n\widetilde\pi_k(1-\widetilde\pi_k)^{n-1}$ gives
\begin{equation}\label{eq:rare-weights}
 \PP_{P^{[k]}}(K=0)=1-\delta+o(\delta),\qquad
 \PP_{P^{[k]}}(K=1)=\delta+o(\delta).
\end{equation}
Let $B_k$ have law $P_*(\,\cdot\mid |X|<b)$.
Write $\PP_{B_k},\E_{B_k}$ for probabilities and expectations of independent copies of $B_k$.
Put $v_k^2=\E_{B_k}B_k^2$. Since $1\le|B_k|\le H$,
\[
 1\le v_k^2=\frac{\E_{P_*}[X^2;|X|<b]}{\PP_{P_*}(|X|<b)}
 \le\frac{\sigma^2}{\PP_{P_*}(|X|<b)}\le2.
\]
The numerator tends to $\sigma^2$ and the denominator to one, so $v_k^2\to\sigma^2$.
For the weighted-sign estimates at thresholds of order $x$, put
\[
 e_k=\frac{(1+x)H}{\sqrt n}+\frac{x^4H^2}{n}.
\]
Equation \eqref{eq:scale-errors} gives $(1+x)H/\sqrt n=o(\delta)$.
Also, $x^4H^2/n\le\{(1+x)^3H/\sqrt n\}^2=o(\delta^2)$.
Hence $e_k=o(\delta)$.

We begin with $q_0$.
Given $K=0$, the observations are independent copies of $B_k$.
Condition on their magnitudes.
The signs are independent and symmetric, and the weights are bounded by $H$.
The sum of squares $Q$ is at least $n$.
Lemma~\ref{lem:weighted-signs} therefore applies at $z=x$.
Its maximal estimate in \eqref{eq:weighted-sign-tails} has error at most $Ce_k$.
Averaging over the magnitudes gives
\begin{equation}\label{eq:no-rare}
 \frac{q_0}{\Phibar(x)}=2+O(e_k)=2+o(\delta).
\end{equation}

Next consider $q_1$.
Put $m=n-1$.
Given $K=1$, the observation with magnitude $b$ has a uniform insertion position
and an independent symmetric sign $s\in\{-1,1\}$.
The other $m$ observations are independent copies of $B_k$ and are independent of the position and sign.
Write their magnitudes as $a_1,\ldots,a_m$ and their signs as $\epsilon_1,\ldots,\epsilon_m$.
The signs are independent Rademacher variables and are independent of the magnitudes.
Put
\[
 Q=\sum_{i=1}^ma_i^2,\qquad
 T_0=0,\qquad T_j=\sum_{i\le j}a_i\epsilon_i,\qquad L_Q=x\sqrt{Q+b^2}.
\]
For large $k$, $L_Q-b>0$, so we may include time zero in the maxima below.
If $r\in\{0,\ldots,m\}$ observations precede the insertion, the crossing event is
\[
 \mathcal A_{r,s}=\left\{\max_{0\le j\le r}T_j\ge L_Q\right\}
 \cup\left\{\max_{r\le j\le m}T_j\ge L_Q-sb\right\}.
\]
Averaging over the insertion position and sign gives
\begin{equation}\label{eq:one-rare-average}
 q_1=\frac1{2n}\sum_{r=0}^m\sum_{s=\pm1}\PP_{B_k}(\mathcal A_{r,s}).
\end{equation}
Put $\mathcal B_s=\{\max_{0\le j\le m}T_j\ge L_Q-sb\}$.
For $s=1$, $\mathcal A_{r,1}\subseteq\mathcal B_1$.
The difference requires a crossing of $L_Q-b$ before insertion.
For $s=-1$, $\mathcal B_{-1}\subseteq\mathcal A_{r,-1}$.
The difference requires a crossing of $L_Q$ before insertion.
Thus, for either sign,
\begin{equation}\label{eq:insertion-difference}
 \mathcal A_{r,s}\mathbin{\triangle}\mathcal B_s
 \subseteq\left\{\max_{0\le j\le r}T_j\ge L_Q-b\right\}.
\end{equation}
To bound this probability, condition on $a_1,\ldots,a_m$.
Put $D=L_Q-b$ and $Q_r=\sum_{i\le r}a_i^2$.
Since $Q\ge m$ and $b=2\sqrt n/x$, we have $b/\sqrt Q\le3/x$ for large $k$.
Thus $D\ge\sqrt Q(x-3/x)$, and
\[
 D^2\ge Q(x^2-C),\qquad D^2\ge Qx^2/2.
\]
For $r\ge1$,
\[
 \frac{D^2}{2Q_r}
 =\frac{D^2}{2Q}+\frac{D^2(Q-Q_r)}{2QQ_r}
 \ge\frac{x^2}{2}-C+cx^2\frac{m-r}{Q}.
\]
Here $Q-Q_r\ge m-r$ and $Q_r\le Q$.
Apply \eqref{eq:weighted-sign-coarse}.
The sum over insertion positions satisfies
\[
 \sum_{r=0}^m\ee^{-cx^2(m-r)/Q}\le1+\frac{Q}{cx^2}.
\]
Hence
\[
 \frac1n\sum_{r=0}^m
 \PP_{B_k}\left(\max_{0\le j\le r}T_j\ge D\mid a_1,\ldots,a_m\right)
 \le\frac{C\ee^{-x^2/2}}n\left(1+\frac Q{x^2}\right).
\]
The empty prefix contributes zero.
Take $\E_{B_k}$ and use $\E_{B_k}Q=mv_k^2\le2n$ and $x^2=\log n\le n$.
Mills' inequalities give $\ee^{-x^2/2}\asymp x\Phibar(x)$, so
\begin{equation}\label{eq:insertion-cost}
 \frac1n\sum_{r=0}^m\PP_{B_k}\left(\max_{0\le j\le r}T_j\ge D\right)
 \le\frac{C\Phibar(x)}x=o(\Phibar(x)).
\end{equation}
Combining \eqref{eq:one-rare-average}, \eqref{eq:insertion-difference} and \eqref{eq:insertion-cost}, we obtain
\begin{equation}\label{eq:one-rare-reduction}
 q_1=\frac12\sum_{s=\pm1}\PP_{B_k}(\mathcal B_s)+o\{\Phibar(x)\}.
\end{equation}
It remains to estimate $\PP_{B_k}(\mathcal B_s)$.
Put $A_k=xb/\sqrt Q=2\sqrt{n/Q}$.
Independence and $\E_{B_k}B_k^4\le H^2v_k^2\le2H^2$ give
\[
 \operatorname{Var}_{B_k}(Q/n)\le\frac{2H^2}{n}\longrightarrow0,\qquad
 \E_{B_k}(Q/n)=\frac{n-1}{n}v_k^2\longrightarrow\sigma^2.
\]
Chebyshev's inequality gives $Q/n\to\sigma^2$ in $\PP_{B_k}$-probability.
Hence $A_k\to a=2/\sigma$ in the same sense.
Since $Q\ge n-1$, we have $A_k\le3$ for large $k$, uniformly in the magnitudes.
The standardized levels satisfy
\begin{equation}\label{eq:one-rare-threshold}
 z_s=\frac{L_Q-sb}{\sqrt Q}
 =x\sqrt{1+A_k^2/x^2}-sA_k/x
 =x+\frac{A_k^2/2-sA_k}{x}+O(x^{-3}),
\end{equation}
with uniform remainder.
For large $k$, $z_s\ge1$.
Since $Q\ge n-1$ and $z_s=x+O(x^{-1})$, Lemma~\ref{lem:weighted-signs} applies uniformly.
Equation \eqref{eq:one-rare-threshold} also gives
$z_s^2-x^2=A_k^2-2sA_k+O(x^{-2})$ uniformly.
Mills' expansion and the bound on $e_k$ therefore give
\[
 \frac{\PP_{B_k}(\mathcal B_s\mid a_1,\ldots,a_m)}{2\Phibar(x)}
 =\frac{\Phibar(z_s)}{\Phibar(x)}\{1+O(e_k)\}
 =\exp\{sA_k-A_k^2/2\}\{1+o(1)\}.
\]
The remainder is uniform in the magnitudes.
Since $A_k\le3$ and $A_k\to a$ in probability,
the expectations of $\exp\{sA_k-A_k^2/2\}$ tend to $\exp\{sa-a^2/2\}$.
Averaging over the magnitudes gives
\[
 \frac{\PP_{B_k}(\mathcal B_s)}{2\Phibar(x)}
 =\exp\{sa-a^2/2\}+o(1).
\]
Substitute this into \eqref{eq:one-rare-reduction} to obtain
\begin{equation}\label{eq:one-rare}
 \frac{q_1}{\Phibar(x)}=2\vartheta_p+o(1),\qquad
 \vartheta_p=\ee^{-a^2/2}\cosh a.
\end{equation}
For $a>0$, $\log\cosh a=\int_0^a\tanh u\,\dd u<a^2/2$.
Hence $\vartheta_p<1$.

It remains to bound the last term in \eqref{eq:rare-decomposition}.
Put $h=\lfloor x^2/2\rfloor$ and split it as
\begin{equation}\label{eq:multiple-count-split}
 \PP_{P^{[k]}}(M_n\ge xV_n,\ K\ge2)
 =\sum_{j=2}^h\PP_{P^{[k]}}(K=j)q_j+\sum_{j=h+1}^n\PP_{P^{[k]}}(K=j)q_j.
\end{equation}
For $2\le j\le h$, condition on $K=j$ and the positions of the observations with magnitude $b$.
Delete those observations and keep the order of the others.
The remaining $n-j$ observations are independent copies of $B_k$.
Since $x^2>j$, (\Mref{eq:deletion}) gives
\[
 q_j\le\PP_{B_k}\left(M_{n-j}^{B_k}\ge\sqrt{x^2-j}\,V_{n-j}^{B_k}\right).
\]
Here $M_{n-j}^{B_k}$ and $V_{n-j}^{B_k}$ refer to the retained sample.
Condition on its magnitudes.
We have $n-j\sim n$, the threshold lies between $x/\sqrt2$ and $x$,
and $x^4H^2/(n-j)\to0$ uniformly.
Lemma~\ref{lem:weighted-signs} therefore gives
\[
 q_j\le C\Phibar(\sqrt{x^2-j})
 \le C\ee^{j/2}\Phibar(x).
\]
The last inequality follows from \eqref{supp:mills} and $\sqrt{x^2-j}\ge x/\sqrt2$.
Since $\PP_{P^{[k]}}(K=j)\le\widetilde\delta^j/j!$,
\begin{equation}\label{eq:moderate-count}
 \sum_{j=2}^h\PP_{P^{[k]}}(K=j)q_j
 \le C\Phibar(x)\sum_{j\ge2}\frac{(\ee^{1/2}\widetilde\delta)^j}{j!}
 =O(\delta^2\Phibar(x)).
\end{equation}
For the second sum in \eqref{eq:multiple-count-split}, use $q_j\le1$.
Markov's inequality applied to $\binom K{h+1}$ gives
\begin{align*}
 \sum_{j=h+1}^n\PP_{P^{[k]}}(K=j)q_j
 &\le\PP_{P^{[k]}}(K\ge h+1)\\
 &\le\E_{P^{[k]}}\binom K{h+1}
 =\binom n{h+1}\widetilde\pi_k^{h+1}
 \le\left(\frac{\ee\widetilde\delta}{h+1}\right)^{h+1}.
\end{align*}
Since $n=2^{2\nu^k}$, $k=O(\log\log n)$.
By \eqref{eq:sharp-separated-parameters},
\[
 h+1\ge\frac{\log n}{2},\qquad
 \log\widetilde\delta
 =-\frac{p-2}{2}\log n+O(\log\log n),
\]
where the implied constant may depend on $p$.
Thus $\PP_{P^{[k]}}(K\ge h+1)\le\exp\{-c_p(\log n)^2\}=o\{\delta\Phibar(x)\}$.
Combining the two sums in \eqref{eq:multiple-count-split} gives
\begin{equation}\label{eq:multi-rare}
 \PP_{P^{[k]}}(M_n\ge xV_n,\ K\ge2)
 =o\{\delta\Phibar(x)\}.
\end{equation}
Substituting \eqref{eq:no-rare}, \eqref{eq:one-rare}, \eqref{eq:multi-rare} and \eqref{eq:rare-weights}
into \eqref{eq:rare-decomposition} gives
\[
 \begin{aligned}
 \frac{\PP_{P^{[k]}}(M_n\ge xV_n)}{\Phibar(x)}
 &=(1-\delta+o(\delta))(2+o(\delta))\\
 &\quad+(\delta+o(\delta))(2\vartheta_p+o(1))+o(\delta)\\
 &=2+2(\vartheta_p-1)\delta+o(\delta).
 \end{aligned}
\]
By \eqref{eq:sharp-law-truncation}, the same expansion holds under $P_*$:
\[
 \frac{\PP_{P_*}(M_{N_k}\ge x_kV_{N_k})}{\Phibar(x_k)}-2
 =2(\vartheta_p-1)\delta_k+o(\delta_k).
\]
This proves (\Mref{eq:sharp-tail-expansion}).
\end{proof}

\subsection{Proof of Lemma~\Mref{lem:strict-survival}}\label{supp:strict-survival-proof}
\begin{proof}
Put $q_0=1$ and fix $0<z<1$.
Apply \eqref{supp:spitzer} to $X_i+\epsilon$, with $\epsilon>0$, and let the Laplace parameter tend to infinity.
As $\epsilon\downarrow0$,
\[
 \bigcap_{j=1}^n\{S_j^X+j\epsilon\le0\}
 \uparrow\bigcap_{j=1}^n\{S_j^X<0\}.
\]
The bounds $\sum_{n\ge0}z^n<\infty$ and $\sum_{j\ge1}z^j/j<\infty$
justify both limits by dominated convergence.
Symmetry gives $\PP_P(S_j^X<0)=(1-a_j)/2$, so
\begin{equation}\label{eq:strict-survival-series}
 \sum_{n\ge0}q_nz^n
 =\exp\left\{\sum_{j\ge1}\frac{z^j}{j}\PP_P(S_j^X<0)\right\}
 =(1-z)^{-1/2}H(z),
\end{equation}
where
\[
 d_j=\frac{a_j}{2j},\qquad
 H(z)=\exp\left\{-\sum_{j\ge1}d_jz^j\right\}.
\]
Apply \eqref{supp:BE} to the distribution function and its left limit at zero.
We obtain
\[
 a_j\le\frac{2C\mu_3}{\sigma^3\sqrt j},\qquad
 0\le d_j\le Bj^{-3/2},
\]
for a constant $B$ depending on $P$.
Thus $\sum_jd_j<\infty$ and
\[
 H(1)=\exp\left\{-\sum_{j\ge1}d_j\right\}=A(P)>0.
\]

Write $H(z)=\sum_{n\ge0}h_nz^n$.
To recover $q_n$ from \eqref{eq:strict-survival-series}, we first bound $h_n$.
Put $d_0=f_0=0$ and $f_j=j^{-3/2}$ for $j\ge1$, and write $*$ for convolution.
Splitting the sum at $j=n/2$ gives
\[
 (f*f)_n=\sum_{j=1}^{n-1}j^{-3/2}(n-j)^{-3/2}
 \le Kn^{-3/2},\qquad
 K=2^{5/2}\sum_{j\ge1}j^{-3/2}<\infty.
\]
Since $d\le Bf$ and $f*f\le Kf$, induction gives
$(d^{*r})_n\le B^rK^{r-1}n^{-3/2}$ for $r,n\ge1$.
Expanding the exponential defining $H$, we have
\[
 h_n=\sum_{r=1}^n\frac{(-1)^r}{r!}(d^{*r})_n\qquad(n\ge1).
\]
Consequently,
\begin{equation}\label{eq:survival-coefficient-bound}
 h_0=1,\qquad
 |h_n|\le\sum_{r\ge1}\frac{(d^{*r})_n}{r!}
 \le C_Pn^{-3/2},\qquad n\ge1.
\end{equation}
In particular, $\sum_n|h_n|<\infty$ and $\sum_nh_n=H(1)=A(P)$.

The coefficients of $(1-z)^{-1/2}$ are
\[
 v_n=4^{-n}\binom{2n}{n}.
\]
By \eqref{eq:strict-survival-series},
\[
 q_n=\sum_{j=0}^nh_jv_{n-j}.
\]
Stirling's formula gives $v_n\sim(\pi n)^{-1/2}$ and $v_n\asymp(n+1)^{-1/2}$.
For $j\le n/2$, the ratios $v_{n-j}/v_n$ are uniformly bounded
and tend to one for each fixed $j$.
Since $\sum_j|h_j|<\infty$, dominated convergence gives
\[
 \sum_{j\le n/2}h_j\frac{v_{n-j}}{v_n}\longrightarrow A(P).
\]
For the remaining terms, \eqref{eq:survival-coefficient-bound} gives
\[
 \sum_{j>n/2}|h_j|v_{n-j}
 \le C_Pn^{-3/2}\sum_{r=0}^{\lfloor n/2\rfloor}v_r
 \le C_Pn^{-1}=o(v_n).
\]
Therefore $q_n/v_n\to A(P)$, and
\[
 q_n\sim\frac{A(P)}{\sqrt{\pi n}}.
\]

Finally, the definition of $R_n$ gives
\[
 \sqrt n\,R_n(0;P)
 =\frac1{\sigma^3}\E_P\!\left[X^2\min\{|X|,\sigma\sqrt n\}\right]
 \longrightarrow\frac{\mu_3}{\sigma^3}
\]
by monotone convergence.
This proves (\Mref{eq:strict-survival-asymptotic}).
\end{proof}

\section{Auxiliary estimates for the Fourier comparison}\label{supp:reduction-lemmas}\label{supp:comparison-lemmas}
Recall the conditional law of $A$ and the noise $(Z,Y)$
from Section~\Mref{sec:comparison}.
The notation in (\Mref{eq:conditionalmoments}) is
\[
 \begin{aligned}
 \mu&=\E_x^* A,\qquad m_r=\E_x^*|A|^r\quad(r=2,3,4),\qquad
 \tau^2=m_2-\mu^2,\\
 \gamma&=m_3/\tau^2,\qquad
 \Lambda_y=(1+y)n|\mu|+(1+y)^3 nm_3.
 \end{aligned}
\]
The characteristic functions $\varphi,K,\Xi,\widehat\rho$
are defined in (\Mref{eq:charfunctions}).

\subsection{Proof of Lemma~\Mref{lem:noisemoments}}\label{supp:noisemoments-proof}
\begin{proof}
\textup{(i)} Interchange the two families in (\Mref{eq:fullnoise}) and replace
$V_0$ by $-V_0$. The joint law is unchanged and $(Z,Y)$ becomes
$(-Z,Y)$. This proves (\Mref{eq:noisesymmetry}).

The conditional law gives $|A|\le h_x$. Hence the construction of $(Z_1,Y_1)$ gives
\[
 |Z_1|\le2h_x,\qquad |Y_1|\le2h_x|Z_1|\le4h_x^2.
\]
The proof of Proposition~\Mref{prop:reflection} gives
$nm_2\ge3/4$ and $n\mu^2\le1/4$.
Hence $\mu^2\le m_2/3$ and $\tau^2\ge m_2/2$.
Since $m_3\le h_xm_2$, we have $\gamma\le2h_x$.
Using $J=8$ and $|V_0|\le1$ in (\Mref{eq:fullnoise}), we obtain
\[
 |Z|\le32h_x+\gamma\le34h_x,\qquad |Y|\le64h_x^2.
\]

Let $A,A'$ be independent with law $P_A$ and put $D=A-A'$.
Since $|a-b|^3\le4(|a|^3+|b|^3)$ and $|D|\le2h_x$,
\[
 \E_x^*|D|^3\le8m_3,\qquad
 \E_x^*D^4\le2h_x\E_x^*|D|^3\le16h_xm_3.
\]
The density of $V$ gives $\E_V|V|=1/3$ and $\E_VV^2=1/6$.
Thus (\Mref{eq:biasedpair}) gives
\[
 \E_{Z_1,Y_1}|Z_1|=\frac{\E_x^*|D|^3}{6\tau^2}\le\frac43\gamma,
 \qquad
 \E_{Z_1,Y_1}Z_1^2=\frac{\E_x^*D^4}{12\tau^2}\le\frac43h_x\gamma.
\]
Also,
\[
 \E_{Z_1,Y_1}|Y_1|\le2h_x\E_{Z_1,Y_1}|Z_1|
 \le\frac83h_x\gamma.
\]
Summing the $2J$ contributions in (\Mref{eq:fullnoise}) now yields
\[
 \E_{Z,Y}|Z|\le\frac{64}{3}\gamma+\gamma=\frac{67}{3}\gamma,
 \qquad \E_{Z,Y}|Y|\le\frac{128}{3}h_x\gamma.
\]
Since $V$ is independent of $(B,B')$ and $\E_VV=0$, $\E_{Z_1,Y_1}Z_1=0$.
Also $\E_{V_0}V_0=0$ by symmetry.
Independence, $\E_{V_0}V_0^2\le1$ and $\gamma^2\le2h_x\gamma$ now give
\[
 \E_{Z,Y}Z^2=2J\E_{Z_1,Y_1}Z_1^2+\gamma^2\E_{V_0}V_0^2
 \le\frac{64}{3}h_x\gamma+\gamma^2\le\frac{70}{3}h_x\gamma.
\]
This proves (\Mref{eq:noisemoments}).

\textup{(ii)} Direct integration of the density of $V$ gives
\[
 \E_V\ee^{\ii tV}=\int_{-1}^1(1-|r|)\ee^{\ii tr}\,\dd r
 =2\int_0^1(1-r)\cos(tr)\,\dd r
 =\frac{2(1-\cos t)}{t^2},
\]
with value one at $t=0$.
By (\Mref{eq:biasedpair}) and the definition of $(Z_1,Y_1)$,
\[
 K(u,v)=\frac1{2\tau^2}\E_x^*
 \left[D^2\E_V\ee^{\ii VD(u+v(A+A'))}\right].
\]
Substituting the preceding formula gives
\[
 K(u,v)=\frac1{\tau^2}\E_x^*
 \frac{1-\cos\{D(u+v(A+A'))\}}{(u+v(A+A'))^2}.
\]
When $u+v(A+A')=0$, the quotient has continuous value $D^2/2$.
This proves (\Mref{eq:Kexact}).
Since $\E_V\ee^{\ii tV}\ge0$, averaging under $\Pi$ gives $K\ge0$.
Also $|K|\le1$, since $K$ is a characteristic function.
Independence of the two families and $V_0$ in (\Mref{eq:fullnoise}) gives
\[
 \Xi(u,v)=K(u,v)^J K(-u,v)^J\widehat\rho(\gamma u),
\]
which is (\Mref{eq:Xi}).

\textup{(iii)} For the independent pair $A,A'$, we have
\[
 |\varphi(u,v)|^2
 =\E_x^*\ee^{\ii D(u+v(A+A'))}
 =\E_x^*\cos\{D(u+v(A+A'))\}.
\]
The last equality follows by interchanging $A,A'$.
On the cone (\Mref{eq:cone}), $|A+A'|\le2h_x$ gives
\[
 |v(A+A')|\le\frac{|u|}{2},\qquad
 \frac{|u|}{2}\le|u+v(A+A')|\le\frac{3|u|}{2}.
\]
Consequently,
\[
 \frac4{9u^2}\le\frac1{(u+v(A+A'))^2}\le\frac4{u^2}.
\]
Multiply by the nonnegative numerator in (\Mref{eq:Kexact}) and take expectations.
The identity for $|\varphi(u,v)|^2$ then gives
\[
 \frac49\frac{1-|\varphi(u,v)|^2}{\tau^2u^2}
 \le K(u,v)\le4\frac{1-|\varphi(u,v)|^2}{\tau^2u^2}.
\]
This proves (\Mref{eq:gap}).
\end{proof}

\subsection{Proof of Lemma~\Mref{lem:conditionalmoments}}\label{supp:conditionalmoments-proof}
\begin{proof}
Here $\PP$ and $\E$ refer to the original law of $\xi=X/(\sigma\sqrt n)$.
We first prove (\Mref{eq:momentbudget1}).
Put $a=1+x$, $q=\PP(|\xi|>h_x)$ and
$v_{\mathrm{tail}}=n\E[\xi^2\ind_{\{|\xi|>h_x\}}]$.
For $F(z)=z^2\min\{|z|,1\}$, (\Mref{eq:truncatedscale}) gives
$R_n(x)=n\E F(a\xi)$.
Since $0<\varepsilon\le1$, for $|z|>\varepsilon$ we have
\[
 F(z)\ge\varepsilon^3,\qquad F(z)\ge\varepsilon z^2.
\]
On $|\xi|>h_x=\varepsilon/a$, these inequalities yield
\[
 nq\le\varepsilon^{-3}R_n(x),\qquad
 v_{\mathrm{tail}}\le\frac{R_n(x)}{\varepsilon a^2}.
\]
Since $x\in\mathcal D_3$, $R_n(x)<\eta(\varepsilon)\le\varepsilon^3/8$.
With $a\ge1$ and $\varepsilon\le1$, the preceding bounds give
$q\le nq\le1/4$ and $v_{\mathrm{tail}}\le1/4$.
Since $\E\xi=0$ and $n\E\xi^2=1$, conditioning gives
\[
 \mu=-\frac{\E[\xi\ind_{\{|\xi|>h_x\}}]}{1-q},\qquad
 nm_2=\frac{1-v_{\mathrm{tail}}}{1-q}.
\]
Thus
\[
 \frac34\le1-v_{\mathrm{tail}}\le nm_2\le\frac1{1-q}\le\frac43.
\]
Cauchy--Schwarz gives
\[
 n\mu^2\le\frac{qv_{\mathrm{tail}}}{(1-q)^2}.
\]
As $q/(1-q)^2\le(1/4)/(3/4)^2=4/9\le1$, we obtain
$n\mu^2\le v_{\mathrm{tail}}\le1/4$.
It follows that
\[
 \frac12\le n\tau^2=nm_2-n\mu^2\le\frac43,
\]
which proves (\Mref{eq:momentbudget1}).

To estimate $n|\mu|$, use $F(z)\ge\varepsilon^2|z|$ for $|z|>\varepsilon$.
Thus
\[
 n\E[|\xi|\ind_{\{|\xi|>h_x\}}]\le\frac{R_n(x)}{\varepsilon^2a}.
\]
The formula for $\mu$ and $1/(1-q)\le4/3\le2$ now give
\[
 n|\mu|\le\frac{n\E[|\xi|\ind_{\{|\xi|>h_x\}}]}{1-q}
 \le\frac{2\varepsilon^{-2}R_n(x)}a.
\]
On $|\xi|\le h_x$, we have $|a\xi|\le\varepsilon\le1$, so
$F(a\xi)=a^3|\xi|^3$. Consequently,
\[
 nm_3=\frac{n\E[|\xi|^3\ind_{\{|\xi|\le h_x\}}]}{1-q}
 \le\frac{2R_n(x)}{a^3}.
\]
The bound $|A|\le h_x$ also gives $nm_4\le h_xnm_3$.
This proves (\Mref{eq:momentbudget2}).

Since $h_x\le1$, we have
\[
 nm_3\le2n\E[|\xi|^3\ind_{\{|\xi|\le1\}}]\le2r_n.
\]
Using $1/2\le n\tau^2\le4/3$, we obtain
\[
 \frac34nm_3\le\gamma=\frac{nm_3}{n\tau^2}\le2nm_3\le4r_n.
\]
Moreover, $n\mu^2\le1/4\le nm_2/3$, so $m_2/2\le\tau^2\le m_2$.
Lyapunov's inequality and $|A|\le h_x$ give
$m_2^{3/2}\le m_3\le h_xm_2$.
Therefore,
\[
 \sqrt{m_2}\le\frac{m_3}{m_2}\le\gamma\le2h_x,
 \qquad m_2\le\gamma^2.
\]
This proves (\Mref{eq:momentbudget3}).

For $y\ge0$ with $|y-x|\le(1+x)^{-1}$,
we have $1+y\le a+a^{-1}\le2a$.
The bounds for $n|\mu|$ and $nm_3$ then give
\[
 \begin{aligned}
 \Lambda_y
 &\le\left\{2\varepsilon^{-2}\frac{1+y}{a}
       +2\left(\frac{1+y}{a}\right)^3\right\}R_n(x)\\
 &\le(4\varepsilon^{-2}+16)R_n(x)
 \le20\varepsilon^{-2}R_n(x).
 \end{aligned}
\]
\end{proof}

\subsection{Proof of Lemma~\Mref{lem:regulators}}\label{supp:regulators-proof}
\begin{proof}
We first estimate $\E_x^*U_i$ in (\Mref{eq:regindividual}).
For a path $a_1,\ldots,a_i$, let $u_i^{(j)}$ be the last regulator
increment of the rotation $(a_{j+1},\ldots,a_i,a_1,\ldots,a_j)$,
where $0\le j<i$.
We claim that
\[
 \sum_{j=0}^{i-1}u_i^{(j)}=\left(-\sum_{j=1}^ia_j\right)^+.
\]
Put $t=\sum_{j=1}^ia_j$.
If $t\ge0$, the last partial sum of each rotation is $t$.
The preceding minimum includes zero, so the last step creates no new negative minimum.
Thus every $u_i^{(j)}$ is zero.
For $t<0$, put $s_0=0$, $s_k=\sum_{j=1}^ka_j$, and extend by $s_{k+i}=s_k+t$.
The partial sums of rotation $j$ are $s_{j+k}-s_j$, for $0\le k\le i$.
Write
\[
 \underline{s}_j=\min_{j\le k<j+i}s_k,
 \qquad u_i^{(j)}=(\underline{s}_j-s_{j+i})^+.
\]
The next window drops $s_j$ and adds $s_{j+i}<s_j$.
Hence $\underline{s}_{j+1}=\min\{\underline{s}_j,s_{j+i}\}$ and
$u_i^{(j)}=\underline{s}_j-\underline{s}_{j+1}$.
Also,
\[
 \underline{s}_i=\min_{0\le k<i}(s_k+t)=\underline{s}_0+t.
\]
Summing gives $\sum_{j=0}^{i-1}u_i^{(j)}=-t$, proving the claim.

The i.i.d.\ law is invariant under cyclic rotations.
Taking expectations therefore gives
\[
 i\E_x^*U_i=\E_x^*S_i^-
 \le\E_x^*|S_i-i\mu|+i|\mu|\le\tau\sqrt i+i|\mu|.
\]
Here $v^-=\max(-v,0)$.
Dividing by $i$ proves the first bound in (\Mref{eq:regindividual}).

For the second moment, first let $i=1$.
Since $\tau^2\le m_2$, Lyapunov's inequality gives
$\tau m_2\le m_2^{3/2}\le m_3$.
Thus $\E_x^*U_1^2\le m_2\le m_3/\tau$.
For $i\ge2$, the definition of $U_i$ gives
\[
 U_i^2\le(A_i^-)^2\ind_{\{W_{i-1}<A_i^-\}}.
\]
By reversal, $W_{i-1}$ has the law of
$\max(0,S_1,\ldots,S_{i-1})$.
We apply Lemma~\ref{supp:survival} to this maximum.
Lemma~\Mref{lem:conditionalmoments} gives $n\tau^2\asymp1$ and $m_2\asymp\tau^2$.
The centered third moment satisfies
\[
 \E_x^*|A-\mu|^3\le4(m_3+|\mu|^3)\le8m_3.
\]
The same lemma and the choice of $\eta(\varepsilon)$ give
$\Lambda_y\le20\varepsilon^{-2}R_n(x)<\eta_0/2$.
For $i\le n$, the drift factor is bounded, since
\[
 \frac{|\mu|\sqrt i}{\tau}
 \le\frac{n|\mu|}{\tau\sqrt n}\le Cn|\mu|\le C\Lambda_y\le C.
\]
Lemma~\ref{supp:survival} therefore yields, for $a\ge0$,
\[
 \PP_x^*(W_{i-1}\le a)
 \le C\min\left\{1,\frac{a+\gamma}{\tau\sqrt i}\right\}.
\]
Using independence of $W_{i-1}$ and $A_i$, we obtain
\[
 \begin{aligned}
 \E_x^*U_i^2
 &\le\int_{\mathbb R}(a^-)^2\PP_x^*(W_{i-1}<a^-)P_A(\dd a)\\
 &\le\frac C{\tau\sqrt i}\E_x^*[(A^-)^2(A^-+\gamma)]\\
 &\le\frac C{\tau\sqrt i}(m_3+\gamma m_2)
 \le\frac{Cm_3}{\tau\sqrt i}.
 \end{aligned}
\]
Here $\gamma m_2\le2m_3$, since $\tau^2\ge m_2/2$.
This proves the second bound in (\Mref{eq:regindividual}).

Since $\sum_{i\le n}i^{-1/2}\le2\sqrt n$, the individual bounds give
\[
 \begin{aligned}
 \sum_{i\le n}\E_x^*U_i&\le2\tau\sqrt n+n|\mu|\le C,\\
 \sum_{i\le n}\E_x^*U_i^2
 &\le\frac{Cm_3\sqrt n}{\tau}
 =\frac{Cnm_3}{\sqrt{n\tau^2}}\le Cnm_3.
 \end{aligned}
\]
This proves (\Mref{eq:regtotal}).
For $i\ge n/2$, the same bounds give
\[
 \begin{aligned}
 \E_x^*U_i
 &\le\frac{\sqrt{2n\tau^2}+n|\mu|}{n}\le\frac Cn,\\
 \E_x^*U_i^2
 &\le\frac{Cm_3}{\sqrt{n\tau^2}}\le Cm_3.
 \end{aligned}
\]
Thus (\Mref{eq:reglate}) holds.

We now prove (\Mref{eq:weighted}) for $y\ge1$.
Put $s=\kappa y^2$ and $\mathcal Z_s=\E_x^*\ee^{-sA^2}$.
Define the tilted law by
\[
 \PP_s(A\in\dd a)=\frac{\ee^{-sa^2}}{\mathcal Z_s}P_A(\dd a).
\]
Use its product law for the sample, and write $\E_s$ for its expectation.
Let $\mu_s,\tau_s^2$ be the mean and variance under $\PP_s$.
The threshold condition and $\varepsilon\le\varepsilon_0/2$ give
\[
 (1+y)h_x\le2\varepsilon\le\varepsilon_0,
 \qquad sh_x^2\le\kappa\varepsilon_0^2\le1.
\]
Since $|A|\le h_x$, we have $\ee^{-sh_x^2}\le\mathcal Z_s\le1$.
On the support of $P_A$, it follows that
\[
 \ee^{-sh_x^2}\le\frac{\ee^{-sa^2}}{\mathcal Z_s}\le\ee^{sh_x^2}.
\]
For independent copies $A,A'$ under $\PP_s$, the variance identity gives
\[
 \tau_s^2=\frac12\E_s(A-A')^2
 \asymp\frac12\E_x^*(A-A')^2=\tau^2.
\]
The same density comparison gives
$\E_sA^2\asymp m_2$ and $\E_s|A|^3\le Cm_3$.

For the mean, subtract $\mu$ before taking expectations:
\[
 \mu_s-\mu
 =\frac{\E_x^*[(A-\mu)(\ee^{-sA^2}-1)]}{\mathcal Z_s}.
\]
Cauchy--Schwarz and Lyapunov's inequality give
$|\mu|m_2\le m_2^{3/2}\le m_3$.
Using $|\ee^{-sA^2}-1|\le sA^2$, we obtain
\[
 |\mu_s-\mu|\le Cs(m_3+|\mu|m_2)\le Csm_3.
\]
Consequently, $n\tau_s^2\asymp1$ and
\[
 n|\mu_s|\le n|\mu|+Cs nm_3\le C\Lambda_y.
\]
Also $\E_s|A-\mu_s|^3\le8\E_s|A|^3\le Cm_3$.
Thus the drift factor in Lemma~\ref{supp:survival} is bounded under $\PP_s$.
The cyclic argument above and the same second-moment argument now give
\[
 \E_sU_i\le C\tau i^{-1/2}+|\mu_s|,
 \qquad \E_sU_i^2\le\frac{Cm_3}{\tau\sqrt i}.
\]

The variables $U_i,Q_i$ depend only on the first $i$ increments.
The product change of measure therefore gives, for $r=1,2$,
\[
 \E_x^*[U_i^r\ee^{-sQ_i}]=\mathcal Z_s^i\E_sU_i^r.
\]
For $0\le sA^2\le1$, we have $\ee^{-sA^2}\le1-sA^2/2$.
Hence
\[
 \mathcal Z_s\le1-sm_2/2\le\ee^{-sm_2/2},
 \qquad \mathcal Z_s^i\le\ee^{-csi/n},
\]
where the last bound uses $nm_2\ge3/4$.
Inserting the bounds under $\PP_s$, we obtain
\[
 \begin{aligned}
 \sum_{i=1}^n\E_x^*[U_i\ee^{-sQ_i}]
 &\le C\tau\sum_{i=1}^ni^{-1/2}\ee^{-csi/n}\\
 &\quad+C(|\mu|+sm_3)\sum_{i=1}^n\ee^{-csi/n},\\
 \sum_{i=1}^n\E_x^*[U_i^2\ee^{-sQ_i}]
 &\le\frac{Cm_3}{\tau}\sum_{i=1}^ni^{-1/2}\ee^{-csi/n}.
 \end{aligned}
\]
The two sequences in these sums are decreasing.
Comparison with their integrals gives
\[
 \begin{aligned}
 \sum_{i=1}^ni^{-1/2}\ee^{-csi/n}
 &\le\int_0^nt^{-1/2}\ee^{-cst/n}\,\dd t\\
 &=\sqrt n\int_0^1r^{-1/2}\ee^{-csr}\,\dd r
 \le\frac{C\sqrt n}{1+\sqrt s},\\
 \sum_{i=1}^n\ee^{-csi/n}
 &\le\int_0^n\ee^{-cst/n}\,\dd t\\
 &=\frac{n(1-\ee^{-cs})}{cs}\le\frac{Cn}{1+s}.
 \end{aligned}
\]
Since $n\tau^2\asymp1$, we conclude that
\[
 \begin{aligned}
 \sum_{i=1}^n\E_x^*[U_i\ee^{-sQ_i}]
 &\le\frac C{1+\sqrt s}
 +C\frac{n|\mu|+s nm_3}{1+s},\\
 \sum_{i=1}^n\E_x^*[U_i^2\ee^{-sQ_i}]
 &\le\frac{Cnm_3}{1+\sqrt s}.
 \end{aligned}
\]
By the definition of $\Lambda_y$,
$n|\mu|\le\Lambda_y/(1+y)$ and $nm_3\le\Lambda_y/(1+y)^3$.
Since $s=\kappa y^2$, we have $1+s\asymp(1+y)^2$ and
\[
 \frac{n|\mu|+s nm_3}{1+s}\le\frac{C\Lambda_y}{(1+y)^3}.
\]
Also $1+\sqrt s\asymp1+y$ and $\Lambda_y<\eta_0/2$.
Restricting the nonnegative sums to $i\le n/2$ proves (\Mref{eq:weighted}).
\end{proof}

\subsection{Proof of Lemma~\Mref{lem:realkernel}}\label{supp:real}
For $i\in\{0,\ldots,n\}$, put $f_{n-i}=f_{n-i}^{\widetilde\psi_y}$
and $T=(n-i)m_2$.
By Lemma~\Mref{lem:conditionalmoments}, $T\le4/3$.
Use the Fourier transform
$\widehat H(u,v)=\iint H(w,q)\ee^{-\ii(uw+vq)}\dd w\dd q$.
Independence of the noise and the future increments gives
\begin{equation}\label{supp:transform}
 \widehat f_{n-i}(u,v)=\widehat H_y(u,v)\Xi(u,v)\varphi(u,v)^{n-i},\qquad
 \widehat H_y(u,v)=\int\chi(r)\frac{2\sin(uy\sqrt r)}u\ee^{-\ii vr}\dd r.
\end{equation}
At $u=0$, use the continuous extension.
For the odd and even parts in (\Mref{eq:evenoddlocal}), put
\[
 D_{n-i}=\tfrac12(\varphi(u,v)^{n-i}-\varphi(-u,v)^{n-i})\Xi,\qquad
 E_{n-i}=\tfrac12(\varphi(u,v)^{n-i}+\varphi(-u,v)^{n-i})\Xi.
\]
Evenness of $\widehat H_y$ and $\Xi$ in $u$ gives
\begin{align}
 o_{n-i}'(w,q)&=\frac{2\ii}{(2\pi)^2}\iint\ee^{\ii vq}
 \left[\int\chi(r)\sin(uy\sqrt r)\cos(uw)\ee^{-\ii vr}\dd r\right]D_{n-i}(u,v)\dd u\dd v,\label{supp:oddinversion}\\
 e_{n-i}''(w,q)&=-\frac2{(2\pi)^2}\iint\ee^{\ii vq}
 \left[\int\chi(r)\sin(uy\sqrt r)\cos(uw)\ee^{-\ii vr}\dd r\right]uE_{n-i}(u,v)\dd u\dd v.\label{supp:eveninversion}
\end{align}
The smooth scalar noise and the compact support of $\chi$
justify inversion and differentiation in these formulas.

To bound these integrals, we first estimate $\varphi$.

\begin{lemma}\label{supp:lowfreq}
There are absolute $\kappa,c>0$ such that
\begin{equation}\label{supp:lowfreqbound}
 |u|\le \kappa/\gamma,\quad |v|\le \kappa|u|/\gamma
 \quad\Longrightarrow\quad
 |\varphi(\pm u,v)|\le\ee^{-c\tau^2u^2}.
\end{equation}
\end{lemma}
The proof is given in Section~\Sref{supp:auxiliary-proofs}.

The threshold condition and $\varepsilon\le\varepsilon_0/2$ give
$h_x\le h:=\varepsilon_0/(1+y)$, so $|A|\le h$.
We use the following bounds for the difference of the two characteristic functions.
\begin{lemma}\label{supp:oddcommon}
Let $z=u-\ii\lambda$ and $\zeta=v+\ii s$,
where $u,v\in\R$ and $\lambda,s\ge0$.
Put $d(z,\zeta)=\varphi(z,\zeta)-\varphi(-z,\zeta)$.
Then
\begin{align}
 \text{(a)}\quad |d(z,\zeta)|&\le C\ee^{C{\lambda h}}\{ |\mu||z|+m_3[|z|^3+|z|(|v|+s)]\},\label{supp:oddcommon0}\\
 \text{(b)}\quad |\partial_u d(u,v)|&\le C\{|\mu|+m_3(u^2+|v|)\},\label{supp:oddcommon1}\\
 \text{(c)}\quad |\partial_u^2d(u,v)|&\le Cm_3|u|.\label{supp:oddcommon2}
\end{align}
\end{lemma}
The proof is given in Section~\Sref{supp:auxiliary-proofs}.

For $n-i\ge1$, the difference of the two powers is
\begin{equation}\label{supp:powerdifference}
 a^{n-i}-b^{n-i}=(a-b)\sum_{j=0}^{n-i-1}a^{n-i-1-j}b^j.
\end{equation}
The following estimates control the integrals in $v$.
\begin{lemma}\label{supp:marknorms}
Let $a_*>0$ and let $a$ be smooth and supported in a fixed compact interval
$[a_0,a_1]\subset(0,\infty)$. Put
\[
 B_u(v)=\int a(r)\ee^{\ii uy\sqrt r-\ii vr}\dd r,
 \quad B^{\sin}_{u,w}(v)=\int a(r)\sin\{u(y\sqrt r+w)\}\ee^{-\ii vr}\dd r.
\]
\textup{(a)} For $0<y\le2$ and $a$ in a fixed bounded smooth family,
\begin{align}
 \int |B_u(v)|\dd v&\le C(1+\sqrt{y|u|}),\label{supp:markexp}\\
 \int |v||B_u(v)|\dd v&\le C(1+y|u|)(1+\sqrt{y|u|}).\nonumber
\end{align}
For $|w|\le y/16$, put $F_y(u)=\min\{y|u|,1+\sqrt{y|u|}\}$. Then
\begin{align}
 \int |B^{\sin}_{u,w}(v)|\dd v&\le CF_y(u),\label{supp:marksine}\\
 \int |v||B^{\sin}_{u,w}(v)|\dd v&\le C(1+y|u|)F_y(u).\nonumber
\end{align}
If $a_*y$ is sufficiently small, then for a fixed integer $N_0\ge3$ and $j=0,1$,
\begin{equation}\label{supp:markoutside}
 \int_{|v|>c|u|/a_*}|v|^j|B^{\sin}_{u,w}(v)|\dd v
 \le C_{N_0}\min(y|u|,1)(1+|u|/a_*)^{j+1-N_0}.
\end{equation}
The exponential version omits the minimum factor.

\textup{(b)} Let $0<\Delta\le1$, with $N$ in a fixed compact subinterval of $(0,\infty)$, and assume
\begin{equation}\label{supp:localampcondition}
 |a^{(k)}(r)|\le C_k\Delta^{-k}\ee^{-c|r-N|/\Delta}
\end{equation}
for each needed derivative order. For $y\ge1$,
\begin{align}
 \int |B_u(v)|\dd v&\le C\{1+\Delta\sqrt{y|u|}\},\label{supp:localizedmark}\\
 \int |v||B_u(v)|\dd v&\le C(\Delta^{-1}+y|u|)\{1+\Delta\sqrt{y|u|}\}.\label{supp:localizedweighted}
\end{align}
If $a_*y$ is sufficiently small, then for a fixed integer $K_0\ge2$,
\begin{equation}\label{supp:localizedoutside}
 \int_{|v|>c|u|/a_*}|B_u(v)|\dd v
 \le C_{K_0}(1+\Delta|u|/a_*)^{1-K_0}.
\end{equation}
The same bounds hold with the opposite sign of the spatial phase.
\end{lemma}
The proof is given in Section~\Sref{supp:auxiliary-proofs}.

For integers $k\ge0$, $j\in\{0,1\}$ and $N_0\ge k+5$, two direct consequences are
\begin{align}
 \int_0^\infty u^k\int_{|v|>cu/a_*}|v|^j|B^{\sin}_{u,w}(v)|\dd v\dd u
 &\le C_{N_0}y a_*^{k+2},\label{supp:outsideintegral}\\
 \int_{c/\gamma}^\infty u^k(1+u/a_*)^{1-N_0}\dd u
 &\le C_{N_0}a_*^{N_0-1}\gamma^{N_0-k-2}.\label{supp:highoutsideintegral}
\end{align}
The first follows from $\min(yu,1)\le yu$ and the substitution $u=a_*z$.
For the second, use $(1+u/a_*)^{1-N_0}\le a_*^{N_0-1}u^{1-N_0}$ and
\[
 \int_{c/\gamma}^{\infty}u^{k+1-N_0}\dd u
 =\frac{(c/\gamma)^{k+2-N_0}}{N_0-k-2}.
\]
We also use
\begin{equation}\label{supp:gaussianintegral}
 \int_\R |u|^k\ee^{-cTu^2}\dd u=C_{k,c}T^{-(k+1)/2},\qquad k>-1.
\end{equation}
The fixed $L_2$ satisfies $L_2\ge2J+4$ and $\kappa\sqrt{L_2}\ge1$.
Since $m_2\le\gamma^2$, $T\ge L_2\gamma^2$ implies $n-i\ge L_2$.

\begin{proof}[Proof of Lemma~\Mref{lem:realkernel}~\textup{(i)}]
We first bound the Fourier integrals for $\gamma\le y\le2$ and
$L_2\gamma^2\le T\le4/3$.
In (\ref{supp:oddinversion}) and (\ref{supp:eveninversion}),
$\cos(uw)$ is independent of $r,v$ and has modulus at most one.
We therefore use (\ref{supp:marksine}) with $w=0$.
Write
\[
 \calL=\{|u|\le\kappa/\gamma,\ |v|\le\kappa|u|/\gamma\},
 \qquad \calC=\{(u,v):u\ne0,\ |v|\le|u|/(4h)\}.
\]
Since $h_x\le h$, we have $\calC\subseteq\calC_x$.

On $\calL$, (\ref{supp:powerdifference}), (\ref{supp:oddcommon0})
and Lemma~\ref{supp:lowfreq} give
\[
 \begin{aligned}
 |\varphi(u,v)^{n-i}-\varphi(-u,v)^{n-i}|
 &\le C\ee^{-cTu^2}\bigl\{(n-i)|\mu||u|\\
 &\hspace{25mm}+(n-i)m_3(|u|^3+|uv|)\bigr\}.
 \end{aligned}
\]
Here $m_2\asymp\tau^2$ by Lemma~\Mref{lem:conditionalmoments}.
Using $F_y(u)\le y|u|$ in (\ref{supp:marksine}) and integrating in $v$
reduces the odd integral to
\[
 Cy\int_\R\ee^{-cTu^2}
 \bigl\{(n-i)|\mu|u^2+(n-i)m_3[u^4+u^2(1+y|u|)]\bigr\}\dd u.
\]
Equation (\ref{supp:gaussianintegral}) bounds this by
\[
 Cy\bigl\{(n-i)|\mu|T^{-3/2}+(n-i)m_3T^{-5/2}
 +(n-i)m_3T^{-3/2}+y(n-i)m_3T^{-2}\bigr\}.
\]
For the even integral, $|E_{n-i}|\le C\ee^{-cTu^2}$.
The extra factor $u$ gives $CyT^{-3/2}$.

For $|u|\le\kappa/\gamma$ outside $\calL$, use $|D_{n-i}|,|E_{n-i}|\le1$.
Since $\gamma\le2h_x\le2h$ and $hy\le\varepsilon_0$,
the fixed choice of $\varepsilon_0$ allows (\ref{supp:outsideintegral}) with $a_*\asymp\gamma$.
The odd and even contributions are at most $Cy\gamma^2$ and $Cy\gamma^3$.

For $|u|>\kappa/\gamma$ inside $\calC$, we use the factors
$K(u,v)^J$ and $K(-u,v)^J$ in $\Xi$.
For fixed $a\ge0$ and integer $j$, if $n-i\ge2j$, then
\begin{equation}\label{supp:gapopt}
 \sup_{0\le z\le1}z^{n-i-j}(1-z^2)^a\le C_a(n-i)^{-a}.
\end{equation}
This follows by setting $t=1-z^2$
and maximizing $\ee^{-(n-i) t/4}t^a$.

By (\Mref{eq:gap}), (\Mref{eq:Xi}) and (\ref{supp:gapopt}),
\[
 \begin{aligned}
 |\varphi(\pm u,v)^{n-i}\Xi(u,v)|
 &\le\frac{C|\varphi(\pm u,v)|^{n-i}(1-|\varphi(\pm u,v)|^2)^J}
 {\tau^{2J}|u|^{2J}}\\
 &\le CT^{-J}|u|^{-2J}.
 \end{aligned}
\]
For either sign, we used the corresponding factor $K(\pm u,v)^J$;
the other factors have modulus at most one.
Integrating with $F_y(u)\le1+\sqrt{y|u|}$ gives the odd bound
\begin{equation}\label{supp:highlong}
 CT^{-J}\{\gamma^{2J-1}+\sqrt y\,\gamma^{2J-3/2}\}.
\end{equation}
The additional $u$ in the even integral reduces both powers of $\gamma$ by one.

For $|u|>\kappa/\gamma$ outside $\calC$, use $|D_{n-i}|,|E_{n-i}|\le1$.
As $4hy\le4\varepsilon_0$, (\ref{supp:highoutsideintegral}) applies with $a_*=4h$ and $N_0=32$.
The odd and even bounds are $Ch^{31}\gamma^{30}$ and $Ch^{31}\gamma^{29}$.

Now let $1\le i\le n/2$.
Lemma~\Mref{lem:conditionalmoments} gives $3/8\le T\le4/3$ and
\[
 \gamma\le2nm_3\le2\Lambda_y
 \le40\varepsilon^{-2}R_n(x)<\eta_0.
\]
The fixed choice of $\eta_0$ ensures $L_2\gamma^2\le3/8$.
Thus all four frequency bounds apply.
Since $y\le2$ and $T\ge3/8$, the bounds on $\calL$ are at most
$C(n|\mu|+nm_3)$ and $C$.
In the other regions, use $\gamma\asymp nm_3$ and boundedness of $h,\gamma$.
For example, (\ref{supp:highlong}) is at most $C\gamma\le Cnm_3$,
and its even counterpart is at most $C$.
Adding the four contributions gives
\[
 \norm{o_{n-i}'}_\infty\le C(n|\mu|+nm_3),\qquad
 \norm{e_{n-i}''}_\infty\le C.
\]
This proves part~\textup{(i)}.
\end{proof}

\begin{proof}[Proof of Lemma~\Mref{lem:realkernel}~\textup{(ii)}]
Let $i>n/2$ and $T\ge\max\{L_2\gamma^2,y^2\}$.
Apply the four frequency bounds proved above.
Use $(n-i)|\mu|\le CTn|\mu|$, $(n-i)m_3\le CTnm_3$
and $y/\sqrt T\le1$ in the low-frequency bounds.
They become
\[
 Cy\{n|\mu|T^{-1/2}+nm_3T^{-3/2}\},\qquad CyT^{-3/2}.
\]
To compare \eqref{supp:highlong} with the odd bound, use $nm_3\asymp\gamma$:
\[
 \frac{T^{-J}(\gamma^{2J-1}+\sqrt y\,\gamma^{2J-3/2})}
      {y\gamma T^{-3/2}}
 \le C\{(\gamma/y)^{2J-2}+(\gamma/y)^{2J-5/2}\}\le C.
\]
The even ratio satisfies the same bound.
Here $T\ge y^2$ and $\gamma\le y$.
The other two frequency regions are also bounded by the required expressions:
\[
 Cy\gamma^2+Ch^{31}\gamma^{30}\le Cy\gamma T^{-3/2},\qquad
 Cy\gamma^3+Ch^{31}\gamma^{29}\le CyT^{-3/2}.
\]
We used $T\le4/3$ and boundedness of $h,\gamma$.
Adding the four contributions proves (\Mref{eq:realpoint}).

\end{proof}

For part~\textup{(iii)}, suppose $\gamma\le y\le2$, $L_2\gamma^2\le T<y^2$ and $|w|\le y/16$.
In \eqref{supp:oddinversion} and \eqref{supp:eveninversion}, split
\[
 \sin(uy\sqrt r)\cos(uw)
 =\tfrac12\sin\{u(y\sqrt r+w)\}+\tfrac12\sin\{u(y\sqrt r-w)\}.
\]
On $\supp\chi$, $y\sqrt r\pm w\ge3y/16$.
For either $a=y\sqrt r\pm w$, integrate twice by parts on the whole real $u$ line:
\begin{equation}\label{supp:ibpu}
 \int_\R\sin(au)F(u,v)\dd u
 =-a^{-2}\int_\R\sin(au)\partial_u^2F(u,v)\dd u.
\end{equation}
The smooth scalar noise gives rapid decay in $u$, so the boundary terms vanish.
The amplitude becomes $\chi(r)/(\sqrt r\pm w/y)^2$,
which belongs to a fixed bounded smooth family.
Factor out $y^{-2}$ and then split the integrals into five regions:
\[
 \begin{array}{c|l}
 1&\calL\cap\calC\\
 2&\calL\setminus\calC\\
 3&\{|u|\le \kappa/\gamma\}\setminus\calL\\
 4&\{|u|>\kappa/\gamma\}\cap\calC\\
 5&\{|u|>\kappa/\gamma\}\setminus\calC.
 \end{array}
\]
We need bounds for $D_{n-i}^{\prime\prime}$ and $(uE_{n-i})^{\prime\prime}$.

For the noise factor $\Xi$, we use the following derivative bounds.
\begin{lemma}\label{supp:highderivatives}
\textup{(a)} On $\calC$,
\begin{equation}\label{supp:Kderivatives}
 |K_u|\le C/|u|,\quad |K_{uu}|\le C/u^2,\quad
 |\Xi_u|\le C/|u|,\quad |\Xi_{uu}|\le C/u^2.
\end{equation}
\textup{(b)} On the whole real plane,
\begin{equation}\label{supp:Xiglobal}
 |\Xi_u|\le C\gamma,\qquad |\Xi_{uu}|\le Ch\gamma.
\end{equation}

\textup{(c)} For $J=8$ and the fixed sufficiently large absolute $L_2$,
if $T\ge L_2\gamma^2$, $|u|\ge \kappa/\gamma$
and $(u,v)\in\calC$, then
\begin{equation}\label{supp:highderivativebound}
 \left|\partial_u^a\{\varphi(\pm u,v)^{n-i}\Xi(u,v)\}\right|
 \le C T^{a/2-6}|u|^{-12},\qquad a=0,1,2.
\end{equation}
\textup{(d)} On the whole real plane,
\begin{equation}\label{supp:powerglobal}
 |\partial_u^a\varphi(\pm u,v)^{n-i}|\le CT^{a/2},\qquad a=0,1,2.
\end{equation}
\end{lemma}
The proof is given in Section~\Sref{supp:auxiliary-proofs}.

Put $R(T,y)=1+\sqrt y\,T^{-1/4}$.
The following lemma bounds the five contributions after (\ref{supp:ibpu}),
including the factor $y^{-2}$.

\begin{lemma}\label{supp:shortregions}
Suppose $L_2\gamma^2\le T<y^2$ and $|w|\le y/16$.
After \eqref{supp:ibpu}, the contributions to $o_{n-i}'$ and $e_{n-i}''$
are bounded, respectively, as follows.
For Region~1,
\begin{equation}\label{supp:region1bounds}
 Cy^{-2}(nm_3+T{n|\mu|})R(T,y),\qquad Cy^{-2}R(T,y).
\end{equation}
For Region~2, they are
\begin{equation}\label{supp:region2bounds}
 Cy^{-2}T({n|\mu|}+nm_3),\qquad Cy^{-2}(\sqrt T+\gamma).
\end{equation}
For Region~3, they are
\begin{equation}\label{supp:region3bounds}
 Cy^{-2}\gamma(T+h\gamma),\qquad
 Cy^{-2}\{\gamma(\sqrt T+\gamma)+\gamma^2(T+h\gamma)\}.
\end{equation}
For Region~4, using $J-2=6$, they are
\begin{align}
 Cy^{-2}T^{-5}\{\gamma^{11}+\sqrt y\,\gamma^{21/2}\},\label{supp:region4odd}\\
 Cy^{-2}T^{-5}\{\gamma^{10}+\sqrt y\,\gamma^{19/2}\}.\label{supp:region4even}
\end{align}
For Region~5, they are
\begin{align}
 Cy^{-2}h^{31}\gamma^{30}(T+h\gamma),\label{supp:region5odd}\\
 Cy^{-2}h^{31}\{\gamma^{30}(\sqrt T+\gamma)+\gamma^{29}(T+h\gamma)\}.\label{supp:region5even}
\end{align}
\end{lemma}
The proof is given in Section~\Sref{supp:auxiliary-proofs}.

\begin{proof}[Proof of Lemma~\Mref{lem:realkernel}~\textup{(iii)}]
We prove (\Mref{eq:realsumodd})--(\Mref{eq:realsumeven}).
For $i>n/2$ with $L_2\gamma^2\le T<y^2$,
use Lemma~\ref{supp:shortregions}.
We first integrate its bounds in $T$.
For Region~1,
\begin{equation}\label{supp:Rintegrals}
 y^{-2}\int_0^{y^2}R(T,y)\dd T=\tfrac73,\qquad
 y^{-2}\int_0^{y^2}TR(T,y)\dd T=\tfrac{15}{14}y^2.
\end{equation}
Thus the odd and even bounds from Region~1 integrate to $C(n|\mu|+nm_3)$ and $C$.
Region~2 gives $Cy^2(n|\mu|+nm_3)/2$ and $C(2y/3+\gamma)$.
For Region~3, the integrals are
\[
 C\gamma(y^2/2+h\gamma),\qquad
 C\{2\gamma y/3+\gamma^2+\gamma^2y^2/2+h\gamma^3\}.
\]
They are bounded by $C nm_3$ and $C$, respectively.
For Region~4, use
\[
 \int_{L_2\gamma^2}^{y^2}T^{-5}\dd T\le C\gamma^{-8}.
\]
The resulting bounds are
\[
 C\gamma\{(\gamma/y)^2+(\gamma/y)^{3/2}\}\le C\gamma,
 \qquad C\{(\gamma/y)^2+(\gamma/y)^{3/2}\}\le C,
\]
because $y\ge\gamma$. Region~5 integrates to
\[
 Ch^{31}\gamma^{30}(y^2/2+h\gamma),\qquad
 Ch^{31}\{\gamma^{30}(2y/3+\gamma)+\gamma^{29}(y^2/2+h\gamma)\}.
\]
The odd integral is at most $C nm_3$, and the even integral is at most $C$.

We now bound the sums by these integrals.
The grid $T_i=(n-i)m_2$ has spacing $m_2\le\gamma^2$.
For $n/2<i\le n$ with $L_2\gamma^2\le T_i<y^2$ and $L_2\ge2$,
we have $[T_i-m_2,T_i]\subset[T_i/2,T_i]$.
Each power $a$ in the preceding bounds therefore satisfies
\[
 m_2 T_i^a\le2^{|a|}\int_{T_i-m_2}^{T_i}t^a\dd t.
\]
The intervals are disjoint and contained in $[L_2\gamma^2/2,y^2]$.
The integral bounds remain valid with lower cutoff $L_2\gamma^2/2$.
For Region~4, this changes only the absolute constant in the integral of $T^{-5}$.
Apply the comparison to each term of the five regional bounds to obtain
(\Mref{eq:realsumodd}) and (\Mref{eq:realsumeven}).
\end{proof}
\subsection{Proof of Lemma~\Mref{lem:growingkernel}}\label{supp:complex}
Here $y\ge2$. Write $h=\varepsilon_0/(1+y)$ and $\Lambda=\Lambda_y$.
The threshold condition gives $h_x\le h$.
For $1\le i\le n$, put $T=(n-i)m_2$ and $N=q+T$.
Lemma~\Mref{lem:conditionalmoments} gives $T\le4/3$. We first derive the Fourier bounds for
$T\ge L_2\gamma^2$, $q\ge0$ and $N\in[1/16,7]$.
Let $\lambda\ge y/\sqrt N$ and $s\in\{0,\lambda y/\sqrt N\}$.
Then $0\le s\le\lambda^2$.
Define
\begin{equation}\label{supp:jointtilt}
 \mathcal Z_{\pm,s}=\E_x^*\ee^{\pm\lambda A-sA^2},\qquad
 \PP_{\pm,s}(A\in\dd a)
 =\frac{\ee^{\pm\lambda a-sa^2}}{\mathcal Z_{\pm,s}}P_A(\dd a).
\end{equation}
Under either product law, the future increments remain independent.
Write $\E_{\pm,s}$ for expectation under the corresponding product law.

The next lemma controls the moments and normalizing constants under the tilted laws.
\begin{lemma}\label{supp:tilt}
Write $\tau_{\pm,s}^2=\operatorname{Var}_{\pm,s}(A)$.
\textup{(a)} There is a fixed absolute $L_3$ such that, for sufficiently small $yh$,
\begin{align}
 c\ee^{-C{\lambda h}}\tau^2\le\tau_{\pm,s}^2&\le C\ee^{C{\lambda h}}\tau^2,&
 \frac{\E_{\pm,s}|A|^3}{\tau_{\pm,s}^2}&\le\Gamma:=L_3\ee^{L_3{\lambda h}}\gamma,\label{supp:tiltedmoments}\\
 |\E_{\pm,s}A|&\le C\ee^{C{\lambda h}}(|\mu|+\lambda{m_2}+sm_3).&&\label{supp:tiltedmean}
\end{align}
\textup{(b)} The normalizing constants satisfy
\begin{equation}\label{supp:normalizer}
 \log\mathcal Z_{\pm,s}\le\pm\lambda\mu+(\lambda^2/2-s){m_2}+C\lambda^3m_3\ee^{C{\lambda h}}.
\end{equation}
\textup{(c)} If, in addition,
\begin{equation}\label{supp:tiltedcondition}
 (y+\lambda)\Gamma\le c.
\end{equation}
then the tilted raw second moments and variances are comparable.
With $\varphi_{\pm,s}(u,v)=\E_{\pm,s}\ee^{\ii uA+\ii vA^2}$, we have
\begin{equation}\label{supp:tiltedlow}
 |\varphi_{\pm,s}(u,v)|^{n-i}\le\exp\{-c\ee^{-C{\lambda h}}Tu^2\},\quad
 |u|\le c/\Gamma,\quad |v|\le c|u|/\Gamma.
\end{equation}
\end{lemma}
The proof is given in Section~\Sref{supp:auxiliary-proofs}.

To shift the derivative integrals, split the boundary amplitude as follows.
Fix a smooth $\theta$ with values in $[0,1]$,
zero on $(-\infty,-1]$ and one on $[1,\infty)$.
With $s=\lambda y/\sqrt N$, split the amplitude at $r=N$ by setting
\begin{align}
 \Delta&=(1+\lambda y)^{-1},\label{supp:localization}\\
 a_0(r)&=\chi(r)\theta((r-N)/\Delta)\ee^{-\lambda y(\sqrt r-\sqrt N)},\nonumber\\
 a_1(r)&=\chi(r)\{1-\theta((r-N)/\Delta)\}\ee^{sr-\lambda y\sqrt r}.\nonumber
\end{align}
Write $B_{j,u}(v)=\int a_j(r)\ee^{-\ii uy\sqrt r-\ii vr}\dd r$,
$z=u-\ii\lambda$, $s_0=0$ and $s_1=s$.
Use the characteristic functions $\varphi_{+,s_j}(u,v)$
and $\varphi_{-,s_j}(-u,v)$ defined above.
\begin{lemma}\label{supp:contour}
Let $b_0=\ee^{-\lambda y\sqrt N}$ and $b_1=\ee^{-sq}$.
\textup{(a)} The first derivative at $w=0$ and the second derivative at $w$ have the representation
\begin{equation}\label{supp:shiftedidentity}
 f_{n-i}^{(a)}(w_a,q)=-\frac1{(2\pi)^2}\sum_{j=0}^1 b_j
 \iint \ee^{\ii vq}B_{j,u}(v)\Xi(z,v+\ii s_j)\mathcal R_{a,j}\dd u\dd v,
\end{equation}
where $a=1,2$, $w_1=0$, $w_2=w$, and
\begin{align*}
 \mathcal R_{1,j}&=\mathcal Z_{+,s_j}^{n-i}\varphi_{+,s_j}(u,v)^{n-i}
                   -\mathcal Z_{-,s_j}^{n-i}\varphi_{-,s_j}(-u,v)^{n-i},\\
 \mathcal R_{2,j}&=\ii z\{\ee^{\ii zw}\mathcal Z_{+,s_j}^{n-i}\varphi_{+,s_j}(u,v)^{n-i}
                   +\ee^{-\ii zw}\mathcal Z_{-,s_j}^{n-i}\varphi_{-,s_j}(-u,v)^{n-i}\}.
\end{align*}
\textup{(b)} The amplitudes satisfy
\begin{equation}\label{supp:amplitudederivative}
 |a_j^{(k)}(r)|\le C_k\Delta^{-k}\ee^{-c|r-N|/\Delta},\qquad j=0,1,
\end{equation}
so Lemma~\ref{supp:marknorms}\textup{(b)} applies.

\textup{(c)} For $|w|\le h$, the factors $b_j\mathcal Z_{\pm,s_j}^{n-i}$,
the noise factors and $\ee^{\lambda|w|}$ are bounded together by
\begin{equation}\label{supp:Emajorant}
 \mathcal E_\lambda=\exp\left\{-\lambda y\sqrt N+\frac{T\lambda^2}{2}
 +(n-i)\lambda|\mu|+L_4(n-i)\lambda^3m_3\ee^{L_4{\lambda h}}+L_4{\lambda h}\right\},
\end{equation}
where $L_4$ is a fixed sufficiently large absolute constant.
\end{lemma}
The proof is given in Section~\Sref{supp:auxiliary-proofs}.

At large $|u|$, we use the following bound for $K$.
\begin{lemma}\label{supp:complexchord}
\textup{(a)} On the real chord cone, for $u\ne0$,
\begin{equation}\label{supp:complexchordbound}
 |K(u-\ii\lambda,v+\ii s)|
 \le C\ee^{C{\lambda h}}\left\{
 \frac{1-|\varphi_{+,s}(u,v)|^2}{\tau^2u^2}+\frac{\lambda^2}{u^2}\right\}.
\end{equation}
\textup{(b)} The corresponding matching product satisfies
\begin{equation}\label{supp:complexhigh}
 |\varphi_{+,s}(u,v)|^{n-i}|K(u-\ii\lambda,v+\ii s)|^J
 \le C\ee^{C{\lambda h}}(T^{-1}+\lambda^2)^J|u|^{-2J}.
\end{equation}
The negative branch is controlled by the reflected matching factor in $\Xi$.
\end{lemma}
The proof is given in Section~\Sref{supp:auxiliary-proofs}.

We now bound the derivative integrals in \eqref{supp:shiftedidentity}.
Assume \eqref{supp:tiltedcondition}.
Take $c>0$ small enough for \eqref{supp:tiltedlow}
and divide the real variables on each shifted contour into
\[
 \begin{array}{c|l}
 \mathrm{I}&|u|\le c/\Gamma,\quad |v|\le c|u|/\Gamma,\\
 \mathrm{II}&|u|\le c/\Gamma,\quad |v|>c|u|/\Gamma,\\
 \mathrm{III}&|u|>c/\Gamma,\quad (u,v)\in\calC,\\
 \mathrm{IV}&|u|>c/\Gamma,\quad (u,v)\notin\calC.
 \end{array}
\]
We bound the four contributions in the next lemma.

\begin{lemma}\label{supp:fourregions}
Assume \eqref{supp:tiltedcondition}. Let $J=8$, and, for $k=0,1,2,3$, put
\begin{gather*}
 \mathcal J_k=T^{-(k+1)/2}+\sqrt y\,\Delta T^{-(2k+3)/4},
 \qquad U=\mathcal J_1+\lambda\mathcal J_0,\\
 V=\mathcal J_3+3\lambda\mathcal J_2+3\lambda^2\mathcal J_1+\lambda^3\mathcal J_0 +(s+\Delta^{-1})U+y(\mathcal J_2+\lambda\mathcal J_1),\\
 P_J=(T^{-1}+\lambda^2)^J,\qquad
 R_m=\Gamma^{2J-m}+\sqrt y\,\Delta\Gamma^{2J-m-1/2}\quad(m=1,2).
\end{gather*}
For the four regions defined above, let $(\mathcal A_r,\mathcal B_r)$ be the quantities in the following table.
\begin{center}
\begin{tabular}{@{}cll@{}}
\toprule
Region & $\mathcal A_r$ (first derivative) & $\mathcal B_r$ (second derivative)\\\midrule
I & $T({n|\mu|} U+nm_3 V)$ & $U$\\
II & $\Gamma\Delta^{-1}$ & $\Gamma^2\Delta^{-2}+\lambda\Gamma\Delta^{-1}$\\
III & $P_JR_1$ & $P_J(R_2+\lambda R_1)$\\
IV & $(h/\Delta)^{31}\Gamma^{30}$ & $(h/\Delta)^{31}(\Gamma^{29}+\lambda\Gamma^{30})$\\\bottomrule
\end{tabular}
\end{center}
With $L_4$ fixed as in \eqref{supp:Emajorant},
\begin{equation}\label{supp:regiontablebound}
 |f_{n-i}'(0,q)|\le C\mathcal E_\lambda\sum_{r=1}^4\mathcal A_r,
 \qquad
 \sup_{|w|\le h}|f_{n-i}''(w,q)|\le C\mathcal E_\lambda\sum_{r=1}^4\mathcal B_r.
\end{equation}
\end{lemma}
The proof is given in Section~\Sref{supp:auxiliary-proofs}.

Use the fixed smooth functions, $J$, and $L_2,L_3,L_4$ above. Put
\[
 \alpha=\min\left\{\frac1{6L_3},\frac1{2L_4}\right\}.
\]
The fixed absolute constant $\varepsilon_0$ satisfies
\[
 \varepsilon_0\le
 \min\left\{\frac1{64\sqrt{L_2}},\frac{\alpha}{8L}\right\}.
\]
Fix $0<t_0<\min\{3/8,\ee^{-1}\}$ such that
\[
 \frac{3\alpha}{\varepsilon_0}t_0\log(1/t_0)\le\frac14.
\]
For fixed absolute constants $C\ge1$ and $c\in(0,1)$,
the fixed absolute constant $\eta_0$ satisfies
\[
 \eta_0\le\min\left\{\frac1{16},\frac{\sqrt{t_0}}{2\sqrt{L_2}}\right\},
 \qquad C\eta_0^{2/3}(1+|\log\eta_0|)\le c.
\]
The constants $C,c$ in this choice are fixed after $\alpha,\varepsilon_0,t_0$.

By Lemma~\Mref{lem:noisemoments}, $|Y|\le1/8$.
Since $\chi$ is supported in $(1/16,4)$ and $Q_{n-i}\ge0$,
$f_{n-i}$ and its derivatives vanish for $q>5$.
Thus all three parts are immediate for $q>5$.

\begin{proof}[Proof of Lemma~\Mref{lem:growingkernel}~\textup{(i)}]
We first obtain bounds for all $T\ge t_0$, $q\ge0$ and $N\in[1/16,7]$.
Lemma~\Mref{lem:conditionalmoments} and the choice of $\eta(\varepsilon)$ give
\[
 \Lambda\le20\varepsilon^{-2}R_n(x)<\eta_0/2,\qquad
 \gamma\le2nm_3\le2\Lambda<\eta_0.
\]
The fixed bound $\eta_0\le\sqrt{t_0}/(2\sqrt{L_2})$ implies $L_2\gamma^2\le t_0$.
Thus the preceding Fourier bounds apply throughout this range.
Choose $\lambda=y\sqrt N/T$ and $s=y^2/T$.
Since $N\ge T$, this choice satisfies $\lambda\ge y/\sqrt N$ and $s\le\lambda^2$.
It minimizes $-\lambda y\sqrt N+T\lambda^2/2$.
The parameter bounds are
\begin{equation}\label{supp:optimaltilt}
 \lambda\asymp_{t_0}(1+y),\quad {\lambda h}\le C\varepsilon_0,\quad
 \Gamma\asymp_{t_0}\gamma\asymp_{t_0} nm_3,\quad \Delta^{-1}\asymp_{t_0}(1+y)^2.
\end{equation}
The bounds in (\ref{supp:optimaltilt}) give
\[
 (y+\lambda)\Gamma\le C(1+y)nm_3\le\frac{C\Lambda}{(1+y)^2}.
\]
The fixed choice of $\eta_0$ therefore gives (\ref{supp:tiltedcondition}).
The moment errors in \eqref{supp:Emajorant} are $O_{t_0}(\Lambda)$,
and
\begin{equation}\label{supp:macroexponent}
 -\lambda y\sqrt N+T\lambda^2/2=-y^2N/(2T).
\end{equation}
Consequently, $\mathcal E_\lambda\le C\ee^{-y^2N/(2T)}$.
Apply Lemma~\ref{supp:fourregions}.
In Region~I, $T\ge t_0$ makes each $\mathcal J_k$ bounded.
Hence $\mathcal A_1\le C((1+y)n|\mu|+(1+y)^3 nm_3)=C\Lambda$
and $\mathcal B_1\le C(1+y)$.
For Regions~II--IV, put $z_*=nm_3(1+y)^3\le\Lambda\le1$.
In Region~II, substitute $\Gamma\asymp nm_3$ and
$\Delta^{-1}\asymp(1+y)^2$ into the second row:
\[
 \mathcal A_2\le C\frac{z_*}{1+y}\le C\Lambda,\qquad
 \mathcal B_2\le C\left\{\frac{z_*^2}{(1+y)^2}+z_*\right\}\le C(1+y).
\]
In Region~III, we have $P_J\le C(1+y)^{16}$ and
$\sqrt y\,\Delta\le C(1+y)^{-3/2}$.
Use $nm_3=z_*(1+y)^{-3}$ in the third row to obtain
\[
 \mathcal A_3\le C(1+y)^{-29}(z_*^{15}+z_*^{29/2}),\qquad \mathcal B_3\le C(1+y)^{-26}(z_*^{14}+z_*^{27/2})
                    +C(1+y)\mathcal A_3.
\]
Since $z_*\le\Lambda\le1$, these bounds are at most $C\Lambda$ and $C(1+y)$.
In Region~IV, use $h/\Delta\le C(1+y)$ and $\Gamma\asymp nm_3$:
\[
 \mathcal A_4\le C(1+y)^{-59}z_*^{30}\le C\Lambda,\qquad \mathcal B_4\le C\{(1+y)^{-56}z_*^{29}+(1+y)^{-58}z_*^{30}\}
                 \le C(1+y).
\]
Adding the four regions gives
\[
 \sum_r\mathcal A_r\le C\Lambda,\qquad
 \sum_r\mathcal B_r\le C(1+y).
\]
Substitute these estimates and \eqref{supp:macroexponent}
into \eqref{supp:regiontablebound}.
We have proved, throughout this general range,
\[
 |f_{n-i}^{\prime}(0,q)|\le C\Lambda\ee^{-y^2N/(2T)},\qquad
 \sup_{|w|\le h}|f_{n-i}^{\prime\prime}(w,q)|\le C(1+y)\ee^{-y^2N/(2T)}.
\]
For $i\le n/2$, Lemma~\Mref{lem:conditionalmoments} gives $3/8\le T\le4/3$.
Since $t_0<3/8$ and $q\le5$, we have $T>t_0$ and $N\in[3/8,19/3]\subset[1/16,7]$.
The general bounds therefore prove part~\textup{(i)}.
\end{proof}

\begin{proof}[Proof of Lemma~\Mref{lem:growingkernel}~\textup{(ii)}]
Let $i>n/2$, $T\ge t_0$ and $q\ge1/16$.
For $q\le5$, $T\le4/3$ gives $N=q+T\in[1/16,19/3]\subset[1/16,7]$.
The general bounds established in part~\textup{(i)} apply.
This proves (\Mref{eq:macroodd})--(\Mref{eq:macroeven}) in part~\textup{(ii)}.
\end{proof}

\begin{proof}[Proof of Lemma~\Mref{lem:growingkernel}~\textup{(iii)}]
Let $i>n/2$, $L_2\gamma^2\le T<t_0$
and $q\in[1/16,5]$. Since $t_0<3/8$, we have $N=q+T\in[1/16,7]$.
Take the tilt

\begin{equation}\label{supp:shorttilt}
 \lambda=\frac\alpha h\log(1/T),\qquad s=\lambda y/\sqrt N.
\end{equation}
Since $T<t_0\le\ee^{-1}$ and $\alpha/\varepsilon_0>4$,
we have $\lambda\ge4y\ge y/\sqrt N$ and $s\le\lambda^2$.
Thus $\lambda h=\alpha\log(1/T)$, so
$\Gamma=L_3T^{-L_3\alpha}\gamma$ by \eqref{supp:tiltedmoments}.
Since $T\ge L_2\gamma^2$, we have
$T^{-L_3\alpha}\le C\gamma^{-2L_3\alpha}$
and $\log(1/T)\le C(1+|\log\gamma|)$.
Since $1-2L_3\alpha\ge2/3$ and $\gamma\le1$, it follows that
\begin{equation}\label{supp:shortcondition}
 (y+\lambda)\Gamma
 \le C(1+y)\gamma^{2/3}(1+|\log\gamma|).
\end{equation}
For small $t>0$, the function $t^{2/3}(1+|\log t|)$ is increasing.
Since $\gamma\le2nm_3\le2\eta_0(1+y)^{-3}$, the right side is at most
\[
 C\eta_0^{2/3}(1+y)^{-1}
 \{1+|\log\eta_0|+3\log(1+y)\}.
\]
The fixed choice of $\eta_0$ therefore gives \eqref{supp:tiltedcondition}.

Put $B_*=1+y+\lambda+T^{-1}$.
Lemma~\ref{supp:fourregions} now gives
\begin{equation}\label{supp:generalcomplexbound}
 |f_{n-i}'(0,q)|\le C({n|\mu|}+nm_3)B_*^{128}\mathcal E_\lambda,
 \quad \sup_{|w|\le h}|f_{n-i}''(w,q)|\le CB_*^{128}\mathcal E_\lambda.
\end{equation}
To bound these factors, use
\[
 s\le\lambda^2,\qquad \Delta^{-1}\le B_*^2,\qquad
 \Gamma\le c,\qquad \Gamma\le C\ee^{C\lambda h} nm_3.
\]
Keep one factor $nm_3$ in the odd contributions from Regions II--IV.
After retaining this factor, the powers of $B_*$ in the four rows are at most
$12,4,20,63$.
For example, the largest term in Region~IV is bounded by
\[
 (h/\Delta)^{31}(\Gamma^{29}+\lambda\Gamma^{30})
 \le CB_*^{62}(1+\lambda)\le CB_*^{63}.
\]
The fixed $L_4$ includes the factors $\ee^{C\lambda h}$.
Bounding the remaining powers by $CB_*^{128}$ proves \eqref{supp:generalcomplexbound}.

We next bound $\mathcal E_\lambda$.
Divide the positive terms in \eqref{supp:Emajorant},
apart from $L_4\lambda h$, by $\lambda y\sqrt N$.
The ratios satisfy
\begin{align}
 \frac{T\lambda^2/2}{\lambda y\sqrt N}&\le\frac{3\alpha}{\varepsilon_0}T\log(1/T),\nonumber\\
 \frac{(n-i)\lambda|\mu|}{\lambda y\sqrt N}&\le CT{n|\mu|}/y,\label{supp:shortratios}\\
 \frac{L_4(n-i)\lambda^3m_3\ee^{L_4{\lambda h}}}{\lambda y\sqrt N}
 &\le C\Lambda T^{1-L_4\alpha}\log^2(1/T).\nonumber
\end{align}
The first ratio is at most $1/4$ by the fixed bound on $t_0$.
Since $L_4\alpha<1$, the other two ratios are at most $1/8$ each
by the fixed bound on $\eta_0$. Their sum is therefore at most $1/2$.
Since $\sqrt N\ge1/4$ and $\lambda=\alpha(1+y)\log(1/T)/\varepsilon_0$,
\[
 \frac{\lambda y\sqrt N}{2}
 \ge\frac{\alpha y(1+y)}{8\varepsilon_0}\log(1/T)
 \ge Ly(1+y)\log(1/T).
\]
The last inequality uses $\varepsilon_0\le\alpha/(8L)$. Therefore,
\begin{equation}\label{supp:shortmajorant}
 \mathcal E_\lambda\le\ee^{L_4{\lambda h}}\ee^{-\lambda y\sqrt N/2}
 \le\ee^{L_4{\lambda h}}T^{L y(1+y)}.
\end{equation}
It remains to combine this decay with the polynomial factors.
For $T<1$, we have $\lambda\le C(1+y)T^{-1}$.
Hence
\[
 B_*^{128}\ee^{L_4{\lambda h}}\le C(1+y)^{128}T^{-512}.
\]
Since $L\ge1600$ and $t_0\le\ee^{-1}$, for $y\ge2$,
\begin{equation}\label{supp:polynomialpaid}
 (1+y)^{128}T^{-512}T^{L y(1+y)}
 \le C\Phibar(y)T^2.
\end{equation}
Indeed, $Ly(1+y)-512-2\ge(L/2)y(1+y)$.
Divide by $T^2\Phibar(y)$ and apply \eqref{supp:mills}.
The ratio is bounded by
$C(1+y)^{129}\exp\{-(L/2)y(1+y)+y^2/2\}$,
which is bounded for $y\ge2$.

Substitute \eqref{supp:shortmajorant} and \eqref{supp:polynomialpaid}
into \eqref{supp:generalcomplexbound}.
Since $n|\mu|+nm_3\le\Lambda$, we obtain
\[
 |f_{n-i}'(0,q)|\le C\Lambda\Phibar(y)T^2,\qquad
 \sup_{|w|\le h}|f_{n-i}''(w,q)|\le C\Phibar(y)T^2,
\]
which proves part~\textup{(iii)}.

\end{proof}

\subsection{Proof of Lemma~\Mref{lem:reflection}}\label{supp:reflection-proof}
\begin{proof}
Put $F_i=f_{n-i}^\psi(W_i,Q_i)$ for $0\le i\le n$.
For $E\in\calF_k$ and $i>k$, $\ind_E$ is $\calF_{i-1}$-measurable.
Using (\Mref{eq:reflection-step}) and the tower property, we have
\[
 \begin{aligned}
 \E_x^*[\ind_E\Delta_i(\psi)]
 &=\E_x^*\!\left[\ind_E\E_x^*(\Delta_i(\psi)\mid\calF_{i-1})\right]\\
 &=\E_x^*\!\left[\ind_E\{\E_x^*(F_i\mid\calF_{i-1})-F_{i-1}\}\right]\\
 &=\E_x^*[\ind_E(F_i-F_{i-1})].
 \end{aligned}
\]
Summing over $i=k+1,\ldots,n$ gives
\[
 \E_x^*\!\left[\ind_E\sum_{i=k+1}^n\Delta_i(\psi)\right]
 =\E_x^*[\ind_E(F_n-F_k)].
\]
Here $F_n=\psi(W_n,Q_n)$ and $F_k=f_{n-k}^\psi(W_k,Q_k)$.
Since $0\le\psi\le1$, we have $0\le F_i\le1$.
Therefore,
\[
 \left|\E_x^*\!\left[\ind_E\sum_{i=k+1}^n\Delta_i(\psi)\right]\right|
 \le\E_x^*[\ind_E|F_n-F_k|]\le\PP_x^*(E),
\]
which proves (\Mref{eq:signedblock}).
\end{proof}

\subsection{Proof of Lemma~\Mref{lem:clock}}\label{supp:clock-proof}
\begin{proof}
Put $k=\lfloor n/2\rfloor$.
We first compare the two events in (\Mref{eq:clock}) with deviations from the mean.
The threshold condition gives $(1+y)h_x\le2\varepsilon\le\varepsilon_0\le1$.
Since $|A|\le h_x$ and $y\ge2$, we have $m_2\le h_x^2\le1/9$.
Lemma~\Mref{lem:conditionalmoments} gives $3/4\le nm_2\le4/3$.
Hence $n\ge27/4>4$ and $k\ge n/3$. Thus
\[
 \E_x^*Q_k=km_2\ge\frac13nm_2\ge\frac14,
 \qquad \frac34\le\E_x^*Q_n=nm_2\le\frac43.
\]
It follows that
\[
 \begin{aligned}
 \{Q_k<1/16\}&\subset\{|Q_k-km_2|\ge1/8\},\\
 \{Q_n\notin[1/4,2]\}&\subset\{|Q_n-nm_2|\ge1/8\}.
 \end{aligned}
\]

To bound these deviations, put $\theta=8(1+y)^2$.
As $(1+y)h_x\le\varepsilon_0$, we have $\theta h_x^2\le8\varepsilon_0^2\le1$.
For $B=A^2-m_2$, we have
\[
 \E_x^*B=0,\qquad |B|\le h_x^2,\qquad
 \E_x^*B^2=m_4-m_2^2\le m_4.
\]
For $|t|\le1$, Taylor's formula gives $\ee^t\le1+t+2t^2$.
Since $|\theta B|\le1$, we obtain
\[
 \E_x^*\ee^{\pm\theta B}
 \le1\pm\theta\E_x^*B+2\theta^2\E_x^*B^2
 \le1+2\theta^2m_4.
\]
Independence therefore gives, for $1\le j\le n$,
\[
 \E_x^*\exp\{\pm\theta(Q_j-jm_2)\}
 \le(1+2\theta^2m_4)^j\le\exp(2\theta^2nm_4).
\]
Lemma~\Mref{lem:conditionalmoments} gives $nm_4\le h_xnm_3$.
Its bound for $\Lambda_y$ and the choice of $\eta(\varepsilon)$ give
$\Lambda_y\le20\varepsilon^{-2}R_n(x)<\eta_0/2$.
Therefore,
\[
 \theta^2nm_4
 \le64\varepsilon_0(1+y)^3nm_3\le64\varepsilon_0\Lambda_y
 \le32\varepsilon_0\eta_0\le\frac12.
\]
The last inequality uses $64\varepsilon_0\eta_0\le1$.
Averaging the two exponential bounds and using $\ee^u-1\le2u$ for $0\le u\le1$, we get
\[
 \E_x^*\!\left[\cosh\{\theta(Q_j-jm_2)\}-1\right]
 \le\exp(2\theta^2nm_4)-1\le4\theta^2nm_4.
\]
Since $\theta/8=(1+y)^2\ge1$, we have
$\cosh(\theta/8)-1\ge c\ee^{\theta/8}$.
Markov's inequality applied to the nonnegative variable
$\cosh\{\theta(Q_j-jm_2)\}-1$ yields
\[
 \PP_x^*(|Q_j-jm_2|\ge1/8)
 \le\frac{\E_x^*[\cosh\{\theta(Q_j-jm_2)\}-1]}{\cosh(\theta/8)-1}
 \le C\theta^2nm_4\ee^{-\theta/8}.
\]
Applying this bound at $j=k,n$ and using the event inclusions gives
\[
 \PP_x^*(Q_k<1/16)+\PP_x^*(Q_n\notin[1/4,2])
 \le C\Lambda_y\ee^{-(1+y)^2}.
\]
Finally, Mills' inequalities give
\[
 \ee^{-(1+y)^2}\le\frac{C\ee^{-y^2/2}}{1+y}\le C\Phibar(y),
\]
which proves (\Mref{eq:clock}).
\end{proof}

\subsection{Proofs of auxiliary lemmas}\label{supp:auxiliary-proofs}

\begin{proof}[Proof of Lemma~\ref{supp:lowfreq}]
Use the conditional law $P_A$ and the moment notation in (\Mref{eq:conditionalmoments}).

Fix a small absolute $c>0$ and then take $\kappa$ sufficiently small.
For $u\ne0$, put $d_0=\min\{c/|u|,|u|/(4|v|)\}$,
where the second term is infinite when $v=0$.
The frequency restrictions then give $d_0\ge C\gamma$,
with $C$ sufficiently large.
Let $A,A'$ be independent copies under the conditional law.
Write $D=A-A'$ and $\mathcal A=\{|A|\le d_0,|A'|\le d_0\}$.
Lyapunov's inequality and $\tau^2\le m_2$ give
\[
 \sqrt{m_2}\le\frac{m_3}{m_2}\le\frac{m_3}{\tau^2}=\gamma.
\]
Since $D^2\le2(A^2+A'^2)$, independence gives
\begin{align*}
 \E_x^*[D^2\ind_{\mathcal A^c}]
 &\le C\{\E_x^*[A^2\ind_{\{|A|>d_0\}}]+{m_2}\PP_x^*(|A|>d_0)\}\\
 &\le C\{m_3/d_0+m_2m_3/d_0^3\}\\
 &\le C\tau^2\{\gamma/d_0+(\gamma/d_0)^3\}\le\tau^2.
\end{align*}
The choice of $\kappa$ gives the last inequality.
Since $\E_x^* D^2=2\tau^2$, we obtain
$\E_x^*[D^2\ind_{\mathcal A}]\ge\tau^2$.

On $\mathcal A$, we have $|v(A+A')|\le|u|/2$
and $|D(u+v(A+A'))|\le3c$.
By the choice of $c$, the identity below and $1-\cos t\ge c't^2$
for small $|t|$ give
\[
 1-|\varphi(u,v)|^2
 =\E_x^*[1-\cos\{D(u+v(A+A'))\}]
 \ge c'u^2\E_x^*[D^2\ind_{\mathcal A}]\ge c'\tau^2u^2.
\]
Apply $\sqrt{1-t}\le\ee^{-t/2}$ to obtain \eqref{supp:lowfreqbound}.
The calculation also applies to $-u$.
When $u=0$, the assumptions force $v=0$.
\end{proof}

\begin{proof}[Proof of Lemma~\ref{supp:oddcommon}]
Recall $h=\varepsilon_0/(1+y)$ and $|A|\le h$; use $z,\zeta$ from the statement.

\textup{(a)} Write $d=2\ii\E_x^*[\sin(zA)\ee^{\ii vA^2-sA^2}]$.
Subtract the linear term to write
\[
 d-2\ii z\mu
 =2\ii\E_x^*[(\sin(zA)-zA)\ee^{\ii vA^2-sA^2}] +2\ii z\E_x^*[A(\ee^{\ii vA^2-sA^2}-1)].
\]
Use $|\ee^{\ii vA^2-sA^2}-1|\le(|v|+s)A^2$
and $|\sin(zA)-zA|\le C|zA|^3\ee^{\lambda h}$.
Taking expectations gives \eqref{supp:oddcommon0}.

\textup{(b)} On the real plane,
\[
 \partial_u d(u,v)=2\ii\E_x^*[A\cos(uA)\ee^{\ii vA^2}].
\]
Subtract $2\ii\mu$ and use $|\cos(uA)-1|\le Cu^2A^2$
and $|\ee^{\ii vA^2}-1|\le|v|A^2$.
This gives (\ref{supp:oddcommon1}).

\textup{(c)} The second derivative is
\[
 \partial_u^2d(u,v)=-2\ii\E_x^*[A^2\sin(uA)\ee^{\ii vA^2}].
\]
The bound $|\sin(uA)|\le|uA|$ gives (\ref{supp:oddcommon2}).
\end{proof}

\begin{proof}[Proof of Lemma~\ref{supp:marknorms}]
Use the amplitude $a$ and the transforms $B_u,B^{\sin}_{u,w}$ in the statement.

We first prove
\begin{equation}\label{supp:sobolev}
 \|\widehat b\|_1\le C\|b\|_2^{1/2}(\|b\|_2+\|b'\|_2)^{1/2}.
\end{equation}
Split the Fourier integral at $|v|=R$.
Cauchy--Schwarz and Plancherel give
\[
 \|\widehat b\|_1
 \le C\bigl(\sqrt R\|b\|_2+R^{-1/2}\|b'\|_2\bigr).
\]
For $b\ne0$, substitute $R=\max\{1,\|b'\|_2/\|b\|_2\}$
to obtain \eqref{supp:sobolev}.
For $b=0$, the inequality is immediate.

\textup{(a)} For $0<y\le2$, put $b(r)=a(r)\ee^{\ii uy\sqrt r}$.
The support of $a$ is bounded away from zero, so differentiation gives
\[
 \|b\|_2\le C,\qquad \|b'\|_2\le C(1+y|u|),\qquad
 \|b''\|_2\le C(1+y|u|)^2.
\]
Apply \eqref{supp:sobolev} to $b$ and $b'$.
Since $\int|v||B_u(v)|\dd v=\|\widehat{b'}\|_1$,
these bounds prove \eqref{supp:markexp}.
When $y|u|\le1$, each fixed derivative of the sine amplitude is $O(y|u|)$.
When $y|u|>1$, split the sine into two exponentials.
The resulting bounds give \eqref{supp:marksine}.

For $|v|>c|u|/a_*$, the smallness of $a_*y$ gives
\[
 \frac{|v|}{2}\le\left|\frac{uy}{2\sqrt r}-v\right|\le\frac{3|v|}{2}
 \qquad(r\in[a_0,a_1]).
\]
Integration by parts $N_0$ times gives $C_{N_0}|v|^{-N_0}$.
Compact support removes the boundary terms.
For $y|u|\le1$, integration against $\ee^{-\ii vr}$ retains the factor $y|u|$.
For $c|u|/a_*\le1$, use the bounds already proved.
For $c|u|/a_*>1$, the last integration is
\[
 \int_{c|u|/a_*}^{\infty}v^{j-N_0}\dd v
 =\frac{(c|u|/a_*)^{j+1-N_0}}{N_0-j-1}.
\]
This proves \eqref{supp:markoutside}.

\textup{(b)} For $y\ge1$, substitute $r=N+\Delta z$.
On the support of $a(N+\Delta z)$, define
\[
 b_u(z)=a(N+\Delta z)
 \exp\left\{\ii uy\left(\sqrt{N+\Delta z}-\sqrt N-\frac{\Delta z}{2\sqrt N}\right)\right\},
\]
and extend $b_u$ by zero elsewhere.
The Fourier transform is
\[
 B_u(v)=\Delta\ee^{\ii uy\sqrt N-\ii vN}
 \widehat b_u\left(\Delta\left[v-\frac{uy}{2\sqrt N}\right]\right).
\]
Changing variables in the frequency integral gives
\begin{align*}
 \int|B_u(v)|\dd v&=\|\widehat b_u\|_1,\\
 \int|v||B_u(v)|\dd v
 &\le\frac{y|u|}{2\sqrt N}\|\widehat b_u\|_1
       +\Delta^{-1}\|\widehat{b_u'}\|_1.
\end{align*}

The remaining phase has first derivative bounded by $Cy|u|\Delta^2|z|$
and second derivative bounded by $Cy|u|\Delta^2$.
Equation \eqref{supp:localampcondition} bounds the derivatives of
$a(N+\Delta z)$ by $C_k\ee^{-c|z|}$.
It follows that
\[
 \|b_u\|_2\le C,\qquad \|b_u'\|_2\le C(1+y|u|\Delta^2),\qquad
 \|b_u''\|_2\le C(1+y|u|\Delta^2)^2.
\]
Apply \eqref{supp:sobolev} to $b_u$ to obtain \eqref{supp:localizedmark}.
Applying it to $b_u'$ also gives
$\|\widehat{b_u'}\|_1\le C(1+y|u|\Delta^2)^{3/2}$.
Hence
\begin{align*}
 \int|v||B_u(v)|\dd v
 &\le C\{\Delta^{-1}(1+y|u|\Delta^2)^{3/2}
           +y|u|(1+\Delta\sqrt{y|u|})\}\\
 &\le C(\Delta^{-1}+y|u|)(1+\Delta\sqrt{y|u|}).
\end{align*}
For the last step, use $(1+t)^{3/2}\le C(1+t)(1+\sqrt t)$ for $t\ge0$
and $\Delta\le1$.
This proves \eqref{supp:localizedweighted}.

To prove \eqref{supp:localizedoutside}, first note that
$\|a^{(k)}\|_1\le C_k\Delta^{1-k}$.
For $|v|>c|u|/a_*$, integration by parts gives
$|B_u(v)|\le C_{K_0}\Delta^{1-K_0}|v|^{-K_0}$.
If $\Delta|u|/a_*\ge1$, then
\[
 \Delta^{1-K_0}\int_{c|u|/a_*}^{\infty}v^{-K_0}\dd v
 \le C_{K_0}(\Delta|u|/a_*)^{1-K_0}.
\]
If $\Delta|u|/a_*<1$, then $y|u|\Delta^2\le a_*y\Delta$ is bounded.
Equation \eqref{supp:localizedmark} therefore bounds the full integral by $C$.
These two cases prove \eqref{supp:localizedoutside}.
\end{proof}

\begin{proof}[Proof of Lemma~\ref{supp:highderivatives}]
Use $T=(n-i)m_2$, $h=\varepsilon_0/(1+y)$ and the cone $\calC$ in (\Sref{supp:real}).
Let $A,A'$ be independent conditional copies, and put $D=A-A'$.

\textup{(a)} We first prove (\ref{supp:Kderivatives}).
Write $\theta=u+v(A+A')$.
Differentiating (\Mref{eq:Kexact}) gives
\begin{align*}
 K_u&=\tau^{-2}\E_x^*\left\{\frac{D\sin(D\theta)}{\theta^2}-\frac{2(1-\cos(D\theta))}{\theta^3}\right\},\\
 K_{uu}&=\tau^{-2}\E_x^*\left\{\frac{D^2\cos(D\theta)}{\theta^2}-\frac{4D\sin(D\theta)}{\theta^3}
 +\frac{6(1-\cos(D\theta))}{\theta^4}\right\}.
\end{align*}
On $\calC$, $|\theta|\asymp|u|$.
Equation (\Mref{eq:gap}) gives
$1-|\varphi(u,v)|^2\le C\tau^2u^2$.
Since $\sin^2 t\le2(1-\cos t)$, Cauchy--Schwarz gives
\[
 \bigl(\E_x^*|D\sin(D\theta)|\bigr)^2
 \le(\E_x^*D^2)(\E_x^*\sin^2(D\theta))
 \le4\tau^2(1-|\varphi(u,v)|^2).
\]
Substitute these bounds into the two derivative formulas to obtain
the assertions for $K$.

For $a=0,1,2$, we have
$|\partial_u^a\widehat\rho(\gamma u)|\le C_a|u|^{-a}$.
Use bounded moments when $\gamma|u|\le1$
and integration by parts when $\gamma|u|>1$.
The product rule now gives the bounds for $\Xi$ in \eqref{supp:Kderivatives}.

\textup{(b)} Also, $|\Xi_u|\le\E_{Z,Y}|Z|$ and $|\Xi_{uu}|\le\E_{Z,Y} Z^2$.
Lemma~\Mref{lem:noisemoments} therefore gives \eqref{supp:Xiglobal}.

\textup{(c)} To prove (\ref{supp:highderivativebound}), write
$\varphi_c=\ee^{-\ii\mu u}\varphi$.
Then
\[
 |\varphi_{c,u}|\le\tau\sqrt{1-|\varphi|^2},\qquad |\varphi_{c,uu}|\le\tau^2.
\]
For the first inequality, use Cauchy--Schwarz on the covariance of
$A-\mu$ and $\ee^{\ii uA+\ii vA^2}$, with
\[
 \E_x^*\left|\ee^{\ii uA+\ii vA^2}-\varphi(u,v)\right|^2=1-|\varphi(u,v)|^2.
\]
In the product rule, at most two derivatives fall on the noise factors.
Each power $\varphi(\pm u,v)^{n-i}$ is paired with $K(\pm u,v)^J$.
With $J=8$, at least six of these factors remain undifferentiated.
By (\Mref{eq:gap}), they satisfy
\[
 K(\pm u,v)^6\le\frac{C(1-|\varphi(\pm u,v)|^2)^6}{\tau^{12}|u|^{12}}.
\]
For the undifferentiated power, \eqref{supp:gapopt} gives
$((n-i)\tau^2)^{-6}|u|^{-12}$.
One derivative of $\varphi_c^{n-i}$ contributes at most
$(n-i)\tau\sqrt{1-|\varphi|^2}|\varphi|^{n-i-1}$.
Two derivatives contribute
$(n-i)\tau^2|\varphi|^{n-i-1}$ and
$(n-i)(n-i-1)\tau^2(1-|\varphi|^2)|\varphi|^{n-i-2}$.
Applying \eqref{supp:gapopt} to these three terms gives
\begin{align*}
 \frac{(n-i)\tau}{\tau^{12}|u|^{12}}
 \sup_{0\le z\le1}z^{n-i-1}(1-z^2)^{13/2}
 &\le\frac{C\sqrt{(n-i)\tau^2}}{((n-i)\tau^2)^6|u|^{12}},\\
 \frac{(n-i)\tau^2}{\tau^{12}|u|^{12}}
 \sup_{0\le z\le1}z^{n-i-1}(1-z^2)^6
 &\le\frac{C(n-i)\tau^2}{((n-i)\tau^2)^6|u|^{12}},\\
 \frac{(n-i)(n-i-1)\tau^2}{\tau^{12}|u|^{12}}
 \sup_{0\le z\le1}z^{n-i-2}(1-z^2)^7
 &\le\frac{C(n-i)\tau^2}{((n-i)\tau^2)^6|u|^{12}}.
\end{align*}
Derivatives of the remaining noise factors contribute
$|u|^{-1}$ or $|u|^{-2}$.
On this range, these are bounded by $C\sqrt T$ and $CT$.
Lemma~\Mref{lem:conditionalmoments} and the choice of $\eta(\varepsilon)$ give
$n|\mu|\le\Lambda_y<\eta_0/2$.
The mean phase contributes $(n-i)|\mu|\le CTn|\mu|\le C\sqrt T$.
Since $T\asymp(n-i)\tau^2$, these bounds prove \eqref{supp:highderivativebound}.

\textup{(d)} Independence gives
\[
 \E_x^*|S_{n-i}-(n-i)\mu|\le\tau\sqrt{n-i},\qquad
 \E_x^*\{S_{n-i}-(n-i)\mu\}^2=(n-i)\tau^2.
\]
Differentiate the characteristic function of $(S_{n-i},Q_{n-i})$.
The two moment bounds give
\begin{align*}
 |\partial_u\varphi(u,v)^{n-i}|
 &\le\tau\sqrt{n-i}+(n-i)|\mu|\le C\sqrt T,\\
 |\partial_u^2\varphi(u,v)^{n-i}|
 &\le(n-i)\tau^2+(n-i)^2\mu^2\le CT.
\end{align*}
Here $(n-i)\tau^2\asymp T$ and $(n-i)|\mu|\le CTn|\mu|\le C\sqrt T$.
The same calculation applies to $-u$, proving \eqref{supp:powerglobal}.
\end{proof}

\begin{proof}[Proof of Lemma~\ref{supp:shortregions}]
Let $\gamma\le y\le2$, $T=(n-i)m_2$ with $L_2\gamma^2\le T<y^2$, and $|w|\le y/16$.
Use $h=\varepsilon_0/(1+y)$ and the sets $\calL,\calC$ from (\Sref{supp:real}).
We first estimate derivatives of the powers of $\varphi(u,v)$ and $\varphi(-u,v)$.
Put
\[
 a_{n-i}(u,v)=(n-i)|\mu||u|+(n-i) m_3(|u|^3+|uv|),\qquad
 b_{n-i}(u,v)=(n-i)|\mu|+(n-i) m_3(u^2+|v|).
\]
Here $(n-i)|\mu|\le CT{n|\mu|}$ and $(n-i) m_3\le CT nm_3$.
On $\calL$, subtract $\ii\mu$ from
$\varphi_u=\ii\E_x^*[A\ee^{\ii uA+\ii vA^2}]$ to obtain
\[
 |\varphi_u|\le|\mu|+|u|m_2+|v|m_3\le|\mu|+Cm_2|u|.
\]
Here $|v|m_3\le\kappa|u|m_3/\gamma=\kappa|u|\tau^2\le\kappa|u|m_2$.
Also $|\varphi_{uu}|\le m_2$.
Since $T\ge L_2\gamma^2$ and $m_2\le\gamma^2$, we have $n-i\ge L_2\ge2J+4$.
Use Lemmas~\ref{supp:oddcommon} and~\ref{supp:lowfreq}
and \eqref{supp:powerdifference}.
Each bound in the following table is multiplied by $C\ee^{-cTu^2}$:
\begin{equation}\label{supp:oddpowerderivatives}
 \begin{array}{c|l}
 j&\text{bound for }|\partial_u^j(\varphi(u,v)^{n-i}-\varphi(-u,v)^{n-i})|\\\hline
 0&a_{n-i}\\
 1&b_{n-i}+\sqrt T(1+\sqrt T|u|)a_{n-i}\\
 2&(n-i) m_3|u|+\sqrt T(1+\sqrt T|u|)b_{n-i}+T(1+Tu^2)a_{n-i}.
 \end{array}
\end{equation}
To obtain the table, differentiate the sum in \eqref{supp:powerdifference}.
One derivative gives at most $(n-i)(n-i-1)$ terms.
Two derivatives give $(n-i)(n-i-1)$ terms with one second derivative
and $(n-i)(n-i-1)(n-i-2)$ terms with two first derivatives.
Lemma~\Mref{lem:conditionalmoments} and the choice of $\eta(\varepsilon)$ give
$n|\mu|\le\Lambda_y<\eta_0/2$.
For the differentiated factors, use
\[
 (n-i)(|\mu|+Cm_2|u|)\le CTn|\mu|+CT|u|\le C\sqrt T(1+\sqrt T|u|).
\]
The remaining powers give $\ee^{-cTu^2}$ by Lemma~\ref{supp:lowfreq}.

For Region~1, on $\calL\cap\calC$, combine \eqref{supp:oddpowerderivatives} and \eqref{supp:Kderivatives}. The Gaussian decay gives
\begin{align}
 \ee^{cTu^2}|D_{n-i}''|&\le C\bigl[(n-i) m_3|u|
  +\{\sqrt T(1+\sqrt T|u|)+|u|^{-1}\}b_{n-i}\nonumber\\
 &\hspace{21mm}+\{T(1+Tu^2)+\sqrt T(1+\sqrt T|u|)/|u|+u^{-2}\}a_{n-i}\bigr],\label{supp:region1odd}\\
 \ee^{cTu^2}|(uE_{n-i})''|&\le C\{T|u|(1+Tu^2)+\sqrt T(1+\sqrt T|u|)+|u|^{-1}\}.
 \label{supp:region1even}
\end{align}
To integrate these expressions, put
$J_k=\int\ee^{-cTu^2}|u|^kF_y(u)\dd u$ for $k\ge-1$.
We have
\[
 J_k\le C_kT^{-(k+1)/2}R(T,y)\quad(k\ge0),\qquad J_{-1}\le CR(T,y).
\]
For $J_{-1}$, substitute $z=\sqrt T|u|$ and put $a=y/\sqrt T>1$.
Use $F_y\le az$ on $z\le1/a$ and $F_y\le1+\sqrt{az}$ on $1/a<z\le1$.
These integrals give $\log a+C\sqrt a\le C\sqrt a$.
For $z>1$, use the Gaussian factor.

The weighted bound in \eqref{supp:marksine} controls the $|v|$ terms
by $(1+y|u|)F_y(u)$.
Thus the integral of \eqref{supp:region1odd} is bounded by
$C\{(n-i)|\mu|+(n-i) m_3/T+(n-i) m_3(1+y/\sqrt T)\}R(T,y)$.
For example, the mixed inverse term is $(n-i) m_3(J_{-1}+yJ_0)$.
Use $(n-i)|\mu|\le CTn|\mu|$, $(n-i) m_3\le CT nm_3$
and $(n-i) m_3(1+y/\sqrt T)\le C nm_3$.
The integral is therefore at most $C(Tn|\mu|+nm_3)R(T,y)$.
Equation \eqref{supp:region1even} gives $CR(T,y)$ in the same way.
Restoring $y^{-2}$ proves \eqref{supp:region1bounds}.

For Region~2, on $\calL\setminus\calC$, use the global noise derivatives \eqref{supp:Xiglobal}. The odd product is bounded, apart from the Gaussian factor, by
\[
 (n-i) m_3|u|+\{\sqrt T(1+\sqrt T|u|)+\gamma\}b_{n-i} +\{T(1+Tu^2)+\gamma\sqrt T(1+\sqrt T|u|)+h\gamma\}a_{n-i}.
\]
The powers of $|u|$ are nonnegative and at most five.
The powers of $|v|$ are at most one.
Since $T\le4/3$ and $h,\gamma$ are bounded, their coefficients are at most
$C((n-i)|\mu|+(n-i)m_3)$.
Here $|v|>|u|/(4h)$.
Apply \eqref{supp:outsideintegral} with $a_*=4h$ and $N_0=32$,
using $\ee^{-cTu^2}\le1$.
The odd contribution is at most $Cy^{-2}T(n|\mu|+nm_3)$.
The even integrand is bounded by $C\{\sqrt T+\gamma+|u|(T+h\gamma)\}$.
Its integral is at most $Cy^{-2}(\sqrt T+\gamma)$.
This proves \eqref{supp:region2bounds}.

In Region~3, where $|u|\le\kappa/\gamma$ outside $\calL$, use \eqref{supp:powerglobal} and \eqref{supp:Xiglobal} to obtain
\[
 |D_{n-i}''|\le C(T+h\gamma),\qquad
 |(uE_{n-i})''|\le C\{\sqrt T+\gamma+|u|(T+h\gamma)\}.
\]
Since $\gamma\le2h$, we have
$\sqrt T\gamma\le(T+\gamma^2)/2\le(T+2h\gamma)/2$.
Apply \eqref{supp:outsideintegral} with $a_*\asymp\gamma$.
The integrals of the constant and $|u|$ terms are bounded by
$C\gamma$ and $C\gamma^2$, respectively.
This gives \eqref{supp:region3bounds}.

For Region~4, with $|u|>\kappa/\gamma$ inside $\calC$,
Lemma~\ref{supp:highderivatives} gives
\[
 |D_{n-i}''|\le CT^{-5}|u|^{-12},\qquad
 |(uE_{n-i})''|\le CT^{-5}|u|^{-11}.
\]
Since $|u|\sqrt T\ge\kappa\sqrt{L_2}\ge1$,
the second bound also includes the term $2E_{n-i}'$.
By \eqref{supp:marksine}, it remains to integrate against
$F_y(u)\le1+\sqrt{y|u|}$.
For the two derivatives, the power integrals are
\[
 \int_{\kappa/\gamma}^{\infty}u^{-12}(1+\sqrt{yu})\dd u
 \le C(\gamma^{11}+\sqrt y\,\gamma^{21/2}),\qquad \int_{\kappa/\gamma}^{\infty}u^{-11}(1+\sqrt{yu})\dd u
 \le C(\gamma^{10}+\sqrt y\,\gamma^{19/2}).
\]
Multiply by $y^{-2}T^{-5}$ to obtain
\eqref{supp:region4odd}--\eqref{supp:region4even}.

In Region~5, where $|u|>\kappa/\gamma$ outside $\calC$,
the bounds from \eqref{supp:powerglobal} and \eqref{supp:Xiglobal} also apply.
Apply \eqref{supp:highoutsideintegral} with $a_*=4h$ and $N_0=32$.
The constant and $|u|$ terms give, respectively,
\[
 h^{31}\int_{\kappa/\gamma}^{\infty}u^{-31}\dd u\le Ch^{31}\gamma^{30},\qquad h^{31}\int_{\kappa/\gamma}^{\infty}u^{-30}\dd u\le Ch^{31}\gamma^{29}.
\]
Multiplying by the coefficients in the Region~3 bounds and by $y^{-2}$
proves \eqref{supp:region5odd}--\eqref{supp:region5even}.

\end{proof}

\begin{proof}[Proof of Lemma~\ref{supp:tilt}]
Let $y\ge2$, $N\in[1/16,7]$, $h=\varepsilon_0/(1+y)$,
$\lambda\ge y/\sqrt N$ and $s\in\{0,\lambda y/\sqrt N\}$.
Use the laws $\PP_{\pm,s}$ in (\ref{supp:jointtilt}).

\textup{(a)} Since $N\ge1/16$ and $yh$ is small,
$sh^2\le4\lambda(yh)h\le C\lambda h$. Hence the density ratio relative to the conditional law
is between $\ee^{-C\lambda h}$ and $\ee^{C\lambda h}$.
For independent copies under either tilted law, the variance identity gives
\[
 \tau_{\pm,s}^2=\tfrac12\E_{\pm,s}(A-A')^2.
\]
The density ratio for this pair is between $\ee^{-C\lambda h}$ and $\ee^{C\lambda h}$.
Comparing with $\tfrac12\E_x^*(A-A')^2=\tau^2$ gives the variance bounds.
Also, $\E_{\pm,s}|A|^3\le\ee^{C\lambda h}m_3$.
Divide by the variance lower bound to obtain \eqref{supp:tiltedmoments}
with the fixed $L_3$.
For $W=\ee^{\pm\lambda A-sA^2}$, we have
$|W-1|\le\ee^{C\lambda h}(\lambda|A|+sA^2)$
and $\mathcal Z_{\pm,s}\ge\ee^{-C\lambda h}$.
Hence
\[
 \left|\E_{\pm,s}A\right|
 =\left|\frac{\mu+\E_x^*[A(W-1)]}{\mathcal Z_{\pm,s}}\right|
 \le C\ee^{C\lambda h}(|\mu|+\lambda m_2+sm_3),
\]
which proves \eqref{supp:tiltedmean}.

\textup{(b)} For the normalizing constant, let $F(a)=\ee^{\pm\lambda a-sa^2}$. Then
\[
 F^{(3)}(a)=F(a)\{(\pm\lambda-2sa)^3-6s(\pm\lambda-2sa)\}.
\]
Since $\lambda\ge y/\sqrt N$, we have $s\le\lambda^2$.
For $|a|\le h$, $|sa|\le sh\le4\lambda yh\le C\lambda$ and $F(a)\le\ee^{\lambda h}$.
Thus $|F^{(3)}(a)|\le C\lambda^3\ee^{C\lambda h}$.
Taylor's formula gives
\[
 \left|F(a)-1\mp\lambda a-(\lambda^2/2-s)a^2\right|
 \le C\lambda^3|a|^3\ee^{C\lambda h}.
\]
Taking expectations yields
\[
 \mathcal Z_{\pm,s}
 \le1\pm\lambda\mu+(\lambda^2/2-s)m_2
       +C\lambda^3m_3\ee^{C\lambda h}.
\]
Apply $\log z\le z-1$ to obtain \eqref{supp:normalizer}.

\textup{(c)} Now assume (\ref{supp:tiltedcondition}).
Lemma~\Mref{lem:conditionalmoments} and the choice of $\eta(\varepsilon)$ give
$n|\mu|\le\Lambda_y<\eta_0/2$.
The same lemma gives $n\tau^2\asymp1$ and $m_2\asymp\tau^2$;
also $\tau\le\sqrt{m_2}\le\gamma$.
Consequently,
\[
 |\mu|/\tau\le Cn|\mu|\gamma,\qquad m_2/\tau\le C\gamma,
 \qquad m_3/\tau=\gamma\tau\le C\gamma^2.
\]
Equation (\ref{supp:tiltedmean}) and the variance lower bound now give
\[
 \frac{|\E_{\pm,s}A|}{\tau_{\pm,s}}
 \le C\ee^{C\lambda h}\{n|\mu|\gamma+\lambda\gamma+s\gamma^2\}\le1.
\]
The last inequality uses $s\le\lambda^2$, (\ref{supp:tiltedcondition}) and the fixed $L_3$.
Consequently,
\[
 \tau_{\pm,s}^2\le\E_{\pm,s}A^2
 =\tau_{\pm,s}^2+(\E_{\pm,s}A)^2\le2\tau_{\pm,s}^2.
\]
Lyapunov's inequality also gives
\[
 \E_{\pm,s}A^2\le\left(\frac{\E_{\pm,s}|A|^3}{\tau_{\pm,s}^2}\right)^2\le\Gamma^2.
\]
Repeat the argument for (\ref{supp:lowfreqbound}) under each tilted law.
Equation \eqref{supp:tiltedmoments} then gives \eqref{supp:tiltedlow}.
\end{proof}

\begin{proof}[Proof of Lemma~\ref{supp:contour}]
Let $y\ge2$, $T=(n-i)m_2$ and $N=q+T\in[1/16,7]$.
Use $h=\varepsilon_0/(1+y)$, $\lambda\ge y/\sqrt N$, $s=\lambda y/\sqrt N$
and $\Delta=(1+\lambda y)^{-1}$ in (\ref{supp:localization}).

\textup{(a)} We prove (\ref{supp:shiftedidentity}).
Apply Fourier inversion in \eqref{supp:transform} and change $u$ to $-u$
in one of the two exponential terms.
This replaces the sine by $\ee^{-\ii uy\sqrt r}$.
At $w=0$, the first derivative contains
$\varphi(u,v)^{n-i}-\varphi(-u,v)^{n-i}$.
The second derivative at $w$ contains
$\ii u(\ee^{\ii uw}\varphi(u,v)^{n-i}+\ee^{-\ii uw}\varphi(-u,v)^{n-i})$.
Both have coefficient $-(2\pi)^{-2}$.
Split $\chi$ at $N$ on width $\Delta$ as in \eqref{supp:localization}.

The differentiated integrands are entire in the two frequencies,
since the increments and noise are bounded and the boundary integral has compact support.
On the strip $z=u+\ii\eta$, $\zeta=v+\ii\xi$,
with $|\eta|\le\lambda$ and $|\xi|\le s$,
integrate by parts in the boundary integral and the smooth scalar bump.
Let $I$ denote either differentiated integrand.
Set $b=0$ for the first derivative and $b=1$ for the second.
Integrating the boundary integral by parts three times,
and the scalar bump integral eight times, gives a finite constant $B$
for the fixed laws and strip parameters such that
\begin{equation}\label{supp:stripenvelope}
 |I(u+\ii\eta,v+\ii\xi)|\le B
 (1+|u|)^{b-5}(1+|v|)^{-3},\qquad b=0,1.
\end{equation}
The bound is integrable, and the vertical-boundary integrals tend to zero.
Cauchy's theorem and dominated convergence allow
$u\mapsto u-\ii\lambda$ in both pieces.
In the lower piece, also shift $v\mapsto v+\ii s$.

The identities
\[
 \varphi(z,v+\ii s_j)=\mathcal Z_{+,s_j}\varphi_{+,s_j}(u,v),\qquad
 \varphi(-z,v+\ii s_j)=\mathcal Z_{-,s_j}\varphi_{-,s_j}(-u,v)
\]
give \eqref{supp:shiftedidentity}.

\textup{(b)} We prove (\ref{supp:amplitudederivative}) on
$r\in\supp\chi\subset[1/16,4]$, with $N\in[1/16,7]$.
The square roots and their denominators below are therefore bounded away from zero.
Moreover,
\[
 \lambda y\ge\frac{y^2}{\sqrt N}\ge\frac4{\sqrt7}>1,\qquad
 \lambda y\asymp\Delta^{-1}.
\]
For $r\ge N$, use
$\sqrt r-\sqrt N=(r-N)/(\sqrt r+\sqrt N)\ge c(r-N)$.
The exponential factor in $a_0$ is then at most $\ee^{-c|r-N|/\Delta}$.
For $r\le N$,
\[
 \lambda y\sqrt r-sr
 =\frac{\lambda y\sqrt r}{\sqrt N(\sqrt N+\sqrt r)}(N-r)
 \ge c\lambda y(N-r),
\]
which gives the same decay for $a_1$.
On $|r-N|\le\Delta$, the two exponential factors are bounded
because $\lambda y\Delta\le1$.
Each derivative of either exponent is bounded by $C_k\lambda y\le C_k\Delta^{-1}$.
The derivatives of the smooth cutoff contribute $C_k\Delta^{-k}$.
The product rule now proves \eqref{supp:amplitudederivative}.

\textup{(c)} It remains to prove (\ref{supp:Emajorant}).
The two prefactors are $\ee^{-\lambda y\sqrt N}\mathcal Z_{\pm,0}^{n-i}$
and $\ee^{-sq}\mathcal Z_{\pm,s}^{n-i}$.
Equation \eqref{supp:normalizer} and
\[
 -sq+T(\lambda^2/2-s)=-\lambda y\sqrt N+T\lambda^2/2
\]
bound the normalizing terms in (\ref{supp:Emajorant}).
The noise factor satisfies
\[
 |\Xi(z,v+\ii s_j)|\le\E_{Z,Y}\ee^{\lambda|Z|+s_j|Y|}\le\ee^{C\lambda h}.
\]
Here we used the support bounds in Lemma~\Mref{lem:noisemoments}
and $s_jh^2\le C\lambda h$.
For the additional spatial factors, use $|z|\le|u|+\lambda$
and $|\ee^{\pm\ii zw}|\le\ee^{\lambda h}$.
Together, these bounds prove (\ref{supp:Emajorant}) with the fixed $L_4$.
\end{proof}

\begin{proof}[Proof of Lemma~\ref{supp:complexchord}]
Use $h=\varepsilon_0/(1+y)$, $\lambda\ge y/\sqrt N$ and $s\in\{0,\lambda y/\sqrt N\}$,
with $y\ge2$ and $N\in[1/16,7]$.
Let $A,A'$ be independent conditional copies, and put $D=A-A'$.

\textup{(a)} For real $a,b$, write
\[
 1-\cos(a+\ii b)=1-\cos a\cosh b+\ii\sin a\sinh b.
\]
Use $\sin^2a\le2(1-\cos a)$,
$\cosh b-1\le Cb^2\ee^{|b|}$ and $|\sinh b|\le|b|\ee^{|b|}$.
The inequality $2rs\le r^2+s^2$ then gives
\[
 |1-\cos(a+\ii b)|\le C\ee^{|b|}(1-\cos a+b^2).
\]
In the chord formula, write the denominator as $t+\ii B$,
where $t=u+v(A+A')$ and $B=-\lambda+s(A+A')$.
On $\calC$, we have $|t|\ge|u|/2$, $|B|\le C\lambda$
and $|DB|\le C\lambda h$.
Hence
\[
 |K(u-\ii\lambda,v+\ii s)|
 \le\frac{C\ee^{C{\lambda h}}}{\tau^2u^2}
 \E_x^*[1-\cos(Dt)+\lambda^2D^2].
\]
The density ratios of the untilted and tilted pair laws
are bounded by $\ee^{C\lambda h}$.
Since $1-\cos(Dt)\ge0$, its untilted expectation is at most
$\ee^{C\lambda h}$ times its tilted expectation,
which is $1-|\varphi_{+,s}(u,v)|^2$.
Also, $\E_x^* D^2=2\tau^2$.
This proves \eqref{supp:complexchordbound}.

\textup{(b)} Let $T=(n-i)m_2\ge L_2\gamma^2$.
Since $m_2\le\gamma^2$, $n-i\ge L_2\ge2J+4$.
Raise the bound to the $J$th power and multiply by $|\varphi_{+,s}(u,v)|^{n-i}$.
In the binomial expansion, apply \eqref{supp:gapopt} to
$z^{n-i}(1-z^2)^j$, with $z=|\varphi_{+,s}(u,v)|$.
This gives
\[
 C\ee^{C{\lambda h}}|u|^{-2J}
 \sum_{j=0}^J\binom Jj((n-i)\tau^2)^{-j}\lambda^{2(J-j)}.
\]
Since $(n-i)\tau^2\asymp T$, the sum is at most $C(T^{-1}+\lambda^2)^J$.
This proves \eqref{supp:complexhigh}.
The same calculation with $-A$ gives the negative branch;
its variance and absolute moments are unchanged.
Bounded support bounds the other noise factors by $\ee^{C\lambda h}$.
\end{proof}

\begin{proof}[Proof of Lemma~\ref{supp:fourregions}]
Let $y\ge2$, $q\ge0$, $T=(n-i)m_2\ge L_2\gamma^2$ and $N=q+T\in[1/16,7]$.
Use $h=\varepsilon_0/(1+y)$ and $\lambda,s,\Delta,\Gamma$ from the statement,
and assume (\ref{supp:tiltedcondition}).

Apply (\ref{supp:shiftedidentity}) to each boundary piece $j=0,1$, with $s_j\in\{0,s\}$.
The bounds below use $s_j\le s$.
The choice of $L_4$ in \eqref{supp:Emajorant} includes
the factors $\ee^{C\lambda h}$ below.

For Region~I, use \eqref{supp:tiltedlow}.
Put $a=\varphi(z,v+\ii s_j)$, $b=\varphi(-z,v+\ii s_j)$
and $Z_*=\max(\mathcal Z_{+,s_j},\mathcal Z_{-,s_j})$.
By \eqref{supp:tiltedlow},
\[
 \sum_{r=0}^{n-i-1}|a|^{n-i-1-r}|b|^r
 \le(n-i)Z_*^{n-i-1}
 \ee^{-c\ee^{-C\lambda h}(n-i-1)\tau^2u^2}.
\]
Here $n-i\ge L_2$, $(n-i)\tau^2\asymp T$ and $Z_*^{-1}\le\ee^{C\lambda h}$.
Equation \eqref{supp:powerdifference} therefore gives
\begin{equation}\label{supp:complexpowers}
 |a^{n-i}-b^{n-i}|\le C\ee^{C{\lambda h}}Z_*^{n-i}(n-i)|a-b|\ee^{-c\ee^{-C{\lambda h}}Tu^2}.
\end{equation}
The factor $Z_*^{n-i}$ is included in $\mathcal E_\lambda$.
Use Lemma~\ref{supp:oddcommon} for $a-b$.
Lemma~\ref{supp:marknorms}\textup{(b)} bounds the integrals in $v$
with weights $1$ and $|v|$.
The integrals in $u$ satisfy
\[
 \int\ee^{-c\ee^{-C{\lambda h}}Tu^2}|u|^k(1+\sqrt y\,\Delta|u|^{1/2})\dd u
 \le C_k\ee^{C_k{\lambda h}}\mathcal J_k.
\]
Use $(n-i)|\mu|\le CTn|\mu|$, $(n-i) m_3\le CT nm_3$
and $|z|\le|u|+\lambda$.
The first-derivative bound is $T(n|\mu|U+nm_3 V)$.
The second derivative has the multiplier $|u|+\lambda$,
so its bound is $U$.

In Region~II, bound the normalized powers by one.
Condition (\ref{supp:tiltedcondition}) makes $y\Gamma$ sufficiently small.
We may therefore apply \eqref{supp:localizedoutside} with $a_*\asymp\Gamma$ and $K_0=32$.
Substitute $u=(\Gamma/\Delta)t$ in the two integrals:
\[
 \begin{aligned}
 \int_0^\infty(1+\Delta u/\Gamma)^{-31}\dd u
 &=\frac\Gamma\Delta\int_0^\infty(1+t)^{-31}\dd t
 \le C\Gamma\Delta^{-1},\\
 \int_0^\infty u(1+\Delta u/\Gamma)^{-31}\dd u
 &=\frac{\Gamma^2}{\Delta^2}\int_0^\infty t(1+t)^{-31}\dd t
 \le C\Gamma^2\Delta^{-2}.
 \end{aligned}
\]
The second derivative also contains $\lambda$ times the first integral.
These are the two bounds in Region~II.

For Region~III,
Lemma~\ref{supp:complexchord} gives $P_J|u|^{-2J}$.
For $m=1,2$, the required power integrals are
\[
 \int_{c/\Gamma}^{\infty}u^{m-1-2J}
       (1+\sqrt y\,\Delta u^{1/2})\dd u
 \le C\{\Gamma^{2J-m}+\sqrt y\,\Delta\Gamma^{2J-m-1/2}\}
 =CR_m.
\]
Thus the first-derivative bound is $P_JR_1$.
The second derivative contains the additional factor $|u|+\lambda$,
which gives $P_J(R_2+\lambda R_1)$ in Region~III.

In Region~IV, $4hy\le4\varepsilon_0$ satisfies the smallness condition
in \eqref{supp:localizedoutside}.
Apply that bound with $a_*=4h$ and $K_0=32$.
Use $(1+t)^{-31}\le t^{-31}$ to obtain
\[
 \begin{aligned}
 \int_{c/\Gamma}^{\infty}(1+\Delta u/(4h))^{-31}\dd u
 &\le C(h/\Delta)^{31}\int_{c/\Gamma}^{\infty}u^{-31}\dd u
 \le C(h/\Delta)^{31}\Gamma^{30},\\
 \int_{c/\Gamma}^{\infty}u(1+\Delta u/(4h))^{-31}\dd u
 &\le C(h/\Delta)^{31}\Gamma^{29}.
 \end{aligned}
\]
Adding $\lambda$ times the first integral gives the second-derivative bound
in Region~IV.

Add the four contributions to obtain \eqref{supp:regiontablebound}.
\end{proof}

\section{Proof of Lemma~\Mref{lem:sandwich}}\label{supp:sandwich-proof}
Use the conditional sample and independent noise fixed at $x$
in Step~4 of Subsection~\Mref{sec:proofpropositions}.
The conditional functions are zero on nonpositive arguments.

\begin{proof}
Let $y>0$ satisfy $|y-x|\le[2(1+x)]^{-1}$.
On $Q_n\in[1/4,2]$, the bound $|Y|\le1/8$ gives
$Q_n+Y\in[1/8,3]$ and $\chi(Q_n+Y)=1$.
Also,
\[
 |\sqrt{Q_n+Y}-\sqrt{Q_n}|
 =\frac{|Y|}{\sqrt{Q_n+Y}+\sqrt{Q_n}}\le2|Y|.
\]
Since $y\le1+x$ and $\sqrt{Q_n}\ge1/2$, the definition of $\delta$ gives
\[
 \delta\sqrt{Q_n}\ge |Z|+y|\sqrt{Q_n+Y}-\sqrt{Q_n}|.
\]
As $W_n\ge0$, the triangle inequality gives
\[
 \ind_{\{W_n<(y-\delta)\sqrt{Q_n}\}}
 \le H_y(W_n+Z,Q_n+Y)
 \le\ind_{\{W_n<(y+\delta)\sqrt{Q_n}\}}.
\]
The lower event is empty if $y-\delta\le0$.
For the absolute sum, replace $W_n$ by $|S_n|$ in the outer indicators
and use $\big||S_n+Z|-|S_n|\big|\le|Z|$.
Average over the sample and the independent noise.
The complement of $Q_n\in[1/4,2]$ costs at most $b_x$,
so the two comparisons prove (\Mref{eq:sandwich}).

By (\Mref{eq:smoothing-size}) and Lemma~\Mref{lem:conditionalmoments},
\[
 \delta\le\frac{400\varepsilon}{1+x},\qquad
 \E_{Z,Y}\delta\le400\gamma\le1600r_n.
\]
The same moment bounds give $\gamma\le2 nm_3$
and $nm_3\le2R_n(x)/(1+x)^3$.
Mills' inequalities therefore imply
\[
 \phi(x)\E_{Z,Y}\delta\le800\phi(x)nm_3
 \le1600\frac{R_n(x)\phi(x)}{(1+x)^3}\le6400w_n(x).
\]

It remains to bound $b_x$.
For $x\ge2$, Lemma~\Mref{lem:conditionalmoments} gives $\gamma\le2h_x<2\le x$.
Thus $y=x$ satisfies (\Mref{eq:comparisonthresholds}).
Lemma~\Mref{lem:clock} gives
$b_x\le C\Lambda_x\Phibar(x)\le Cw_n(x)$.

For $0<x<2$, the mean of $Q_n$ lies in $[3/4,4/3]$
and its variance is at most $nm_4=n\E_x^* A^4$.
Chebyshev's inequality and Lemma~\Mref{lem:conditionalmoments} give
\[
 b_x\le4nm_4\le4h_x nm_3\le CR_n(x)\le Cw_n(x),
\]
where the last step uses $\Phibar(x)\ge\Phibar(2)$.
Together these bounds prove (\Mref{eq:sandwichcost}).
\end{proof}

\end{document}